\documentclass[10pt]{amsart}
\usepackage[dvipdfmx]{graphicx}
\usepackage{amsmath,amscd,amsthm,amsxtra,yhmath,stackrel}
\usepackage{epsfig,graphics,color,colortbl}
\usepackage{amssymb,latexsym}
\usepackage{mathrsfs}
\usepackage{bm}
\usepackage[poly,all]{xy}
\usepackage{hyperref}
\usepackage{tikz-cd}
\usepackage[draft]{todonotes}
\usepackage{upgreek}

\usepackage{ytableau}
\newcommand{\tl}[1]{\substack{\scalebox{0.8}{#1}}}

\newtheorem{thm}{\bf Theorem}[section]
\newtheorem{df}[thm]{\bf Definition}
\newtheorem{prop}[thm]{\bf Proposition}
\newtheorem{cor}[thm]{\bf Corollary}
\newtheorem{lem}[thm]{\bf Lemma}
\newtheorem{rem}[thm]{\bf Remark}
\newtheorem{ex}[thm]{\bf Example}

\numberwithin{equation}{section}

\newcommand{\mc}{\mathcal}
\newcommand{\mf}{\mathfrak}

\newcommand{\mb}{\bm}

\newcommand{\pf}{\noindent{\bfseries Proof. }}
\newcommand{\ov}{\overline}

\newcommand{\U}{{\mc U}}
\newcommand{\V}{\mc{V}}
\newcommand{\W}{\mc{W}}
\newcommand{\cP}{\mathscr{P}}

\newcommand{\I}{\mathbb{I}}
\newcommand{\be}{{\bf e}}

\newcommand{\Z}{\mathbb{Z}}
\newcommand{\Q}{\mathbb{Q}}
\newcommand{\C}{\mathbb{C}}

\newcommand{\e}{\epsilon}
\newcommand{\de}{\delta}
\newcommand{\ude}{\updelta}

\newcommand{\td}{\widetilde}

\newcommand{\lang}{\langle}
\newcommand{\rang}{\rangle}
\newcommand{\gl}{\mf{gl}}
\newcommand{\agl}{\widehat{\mf{gl}}}

\newcommand{\La}{\Lambda}

\newcommand{\la}{\lambda}

\newcommand{\ot}{\otimes}

\newcommand{\si}{(-1)^{\e_i}}

\newcommand{\bq}{{\bf q}}

\usepackage{stmaryrd}
\usepackage{verbatim}
\usepackage{mathtools}
\usepackage[title]{appendix}

\theoremstyle{definition}

\theoremstyle{remark}

\newcommand{\ket}[1]{\left|#1\right\rangle}

\newcommand{\abs}[1]{\left|#1\right|}
\begin{document}
\title
{Affinization of $q$-oscillator representations of $U_q(\mf{gl}_n)$}

\author{JAE-HOON KWON}

\address{Department of Mathematical Sciences and RIM, Seoul National University, Seoul 08826, Korea}
\email{jaehoonkw@snu.ac.kr}

\author{SIN-MYUNG LEE}

\address{Department of Mathematical Sciences, Seoul National University, Seoul 08826, Korea}
\email{luckydark@snu.ac.kr}

\thanks{This work is supported by the National Research Foundation of Korea(NRF) grant funded by the Korea government(MSIT) (No.\,2019R1A2C1084833 and 2020R1A5A1016126).}

\begin{abstract}
We introduce a category $\widehat{\mc{O}}_{\rm osc}$ of $q$-oscillator representations of the quantum affine algebra $U_q(\widehat{\gl}_n)$. We show that $\widehat{\mc{O}}_{\rm osc}$ has a family of irreducible representations, which naturally corresponds to finite-dimensional irreducible representations of quantum affine algebra of untwisted affine type $A$. It is done by constructing a category of $q$-oscillator representations of the quantum affine superalgebra of type $A$, which interpolates these two family of irreducible representations.
The category $\widehat{\mc{O}}_{\rm osc}$ can be viewed as a quantum affine analogue of the semisimple tensor category generated by unitarizable highest weight representations of $\gl_{u+v}$ ($n=u+v$) appearing in the $(\gl_{u+v},\gl_\ell)$-duality on a bosonic Fock space.
\end{abstract}

\maketitle
\setcounter{tocdepth}{1}
\tableofcontents

\section{Introduction}

Let $\mf{g}$ be a simple Lie algebra and let $U_q(\widehat{\mf g})$ be the quantized enveloping algebra (without derivation) associated to the affine Lie algebra $\widehat{\mf g}$. The theory of finite-dimensional $U_q(\widehat{\mf g})$-modules has a rich structure in connection with various topics in mathematics  and mathematical physics, especially including the quiver Hecke algebras \cite{KKK,KKK2} and the cluster algebras \cite{HL10,KKOP1,KKOP2}, where the notion of $R$ matrix underlies these connections in a fundamental way.

Beyond the finite-dimensional $U_q(\widehat{\mf g})$-modules, there are two important categories whose irreducible objects are in general infinite-dimensional.
In \cite{H05}, Hernandez introduces a category $\widehat{\mc O}$ of $U_q(\widehat{\mf{g}})$-modules which can be viewed as a quantum affine analogue of the category $\mc{O}$ for a symmetrizable Kac-Moody algebra $\mf{g}$.
In \cite{MY}, when $\mf{g}$ is finite-dimensional, Mukhin and Young extend many fundamental results in finite-dimensional representations to the case of $\widehat{\mc{O}}$ including the classification of irreducible objects and the theory of $q$-characters, and then study minimal affinization of (parabolic) Verma $U_q(\mf{g})$-modules.
In \cite{HJ}, Hernandez and Jimbo introduce a category of $U_q(\mf{b})$-modules, where $\mf{b}$ is a Borel subalgebra of $\widehat{\mf g}$. This category plays a  crucial role in generalizing the Baxter’s relation and hence giving a new interpretation of $q$-characters of finite-dimensional $U_q(\widehat{\mf g})$-modules \cite{FH}.
When the modules are restricted as $U_q(\mf{b})$-modules, $\widehat{\mc O}$ is a subcategory of the Hernandez-Jimbo's category.

The purpose of this paper is to introduce a new subcategory of $\widehat{\mc O}$, which consists of infinite-dimensional $U_q(\widehat{\gl}_n)$-modules, but still has properties very similar to those of finite-dimensional $U_q(\widehat{\gl}_n)$-modules.

We assume $n\ge 4$ and fix $2\le r\le n-2$. We consider a family of infinite-dimensional irreducible representations of $\gl_n$ called {\em oscillator representations}, which appear in the $(\gl_{u+v},\gl_\ell)$-duality on a bosonic Fock space $S({\C^u}^*\ot {\C^\ell}^*\oplus \C^v \ot \C^\ell)$, where $\ell\ge 1$, $u=r$ and $v=n-r$ \cite{KV}. These representations are unitarizable highest weight representations, which are of importance in themselves in the study of real Lie groups (cf.~\cite{HTW} and references therein).
From a viewpoint of the theory of reductive dual pairs \cite{H}, they form a semisimple tensor category $O_{\rm osc}$, whose Grothedieck ring is a homomorphic image of the dual of the coalgebra $K(GL)=\bigoplus_{\ell\ge 0}K(GL_\ell)$. Here $K(GL_\ell)$ is the Grothendieck group of the category of the finite-dimensional representations of $GL_\ell$ over the complex numbers with integral weights and the comultiplication is induced from the branching rule.

In this paper, we construct a quantum affine analogue of oscillator representations of $\gl_n$ in a systematic way.
The main idea is to replace the finite-dimensional $U_q(\gl_n)$-modules in the category of finite-dimensional $U_q(\widehat{\gl}_n)$-modules by the $q$-analogues of oscillator representations of $\gl_n$ in $O_{\rm osc}$, and then as in \cite{Kas02} to apply the fusion construction using $R$ matrix to the tensor products of fundamental type $q$-oscillator representations. We remark that the fundamental type $q$-oscillator representation also appears in \cite{Ku18}, which is given in a more general way, in the study of the tetrahedron equation and its solutions. There are several other works on special $q$-oscillator representations of type $A_{n-1}^{(1)}$ \cite{BGKNR16,BGKNR17,KuO}, but they do not seem to be the ones constructed here.
This work is also partly motivated by \cite{KO21}, where a special case of higher level oscillator representation for the quantum affine (super)algebras of other type is constructed. 

Let us explain our results in more details.
The first main result in this paper is the construction of $q$-oscillator representations $\V^\la$ of $U_q(\gl_n)$, which are parametrized by a set of generalized partitions $\la$ of length $\ell\ge 1$. We also show that they form a semisimple tensor category $\mc{O}_{\rm osc}$ with the irreducible objects $\V^\la$, which can be viewed as a $q$-analogue of $O_{\rm osc}$. 
We then define a monoidal or tensor subcategory $\widehat{\mc O}_{\rm osc}$ of $\widehat{\mc O}$, where each object belongs to $\mc{O}_{\rm osc}$ as a $U_q(\gl_n)$-module.

Next, we consider a construction of irreducible $q$-oscillator representations of $U_q(\widehat{\gl}_n)$.
The level one $q$-oscillator representation $\V^{(l)}\in \mc{O}_{\rm osc}$ for $l\in \Z$ is naturally extended to a $U_q(\widehat{\gl}_n)$-module, which we call a fundamental representation in $\widehat{\mc O}_{\rm osc}$, and denote by $\W_l$. We define a normalized $R$ matrix on $\W_l(z_1)\ot \W_m(z_2)$ for $l,m\in\Z$, where $z_1, z_2$ denote the generic spectral parameters. 
The tensor product $\W_l(z_1)\ot \W_m(z_2)$ is semisimple as a $U_q(\gl_n)$-module, and decomposes into an infinite direct sum of level two $q$-oscillator representations:
\begin{equation}\label{eq:decomp of two tensor}
\W_l(z_1)\ot \W_m(z_2) \cong \W_m(z_2)\ot \W_l(z_1) = 
\bigoplus_{t\in\Z_+} \mc{V}^{(l_1+t,l_2-t)},
\end{equation}
where $l_1=\max\{l,m\}$ and $l_2=\min\{l,m\}$.
As the second main result in this paper, we give an explicit formula for the spectral decomposition of the normalized $R$ matrix $\mc{R}^{\rm norm}_{(l,m)}$ on $\W_l(z_1)\ot \W_m(z_2)$ with respect to the decomposition \eqref{eq:decomp of two tensor}:
\begin{equation}\label{eq:spectral decomp}
\begin{split}
\mc{R}^{\rm norm}_{(l,m)} 
&= \mc P^{l,m}_0 + \sum_{t=1}^\infty\prod_{i=1}^{t}\dfrac{1-q^{d+2i}z}{z-q^{d+2i}} \mc P^{l,m}_t\\
\end{split}
\end{equation}
(up to a constant multiple) where $z=z_1/z_2$ and $d=|l-m|$, and $\mc{P}^{l,m}_t$ ($t\ge 0$) is a suitably normalized projection onto $\mc{V}^{(l_1+t,l_2-t)}$.
Then we apply the fusion construction to have a family of irreducible $U_q(\widehat{\gl}_n)$-modules as the non-zero image of the composition of the normalized $R$ matrices:
\begin{equation}\label{eq:fusion}
 \xymatrixcolsep{2.5pc}\xymatrixrowsep{3pc}\xymatrix{ 
\mathcal{W}_{l_{1}}(c_{1})\otimes\cdots\otimes\mathcal{W}_{l_{\ell}}(c_{\ell})\ \ar@{->}[r] 
& \ \mathcal{W}_{l_{\ell}}(c_{\ell})\otimes\cdots\otimes\mathcal{W}_{l_{1}}(c_{1})},
\end{equation}  
for any pair of ${\mb l}=(l_1,\dots,l_\ell)\in \Z^\ell$ and ${\mb c}=(c_1,\dots,c_\ell)\in (\Q(q)^\times)^\ell$ such that $c_i/c_j\not\in q^{|l_i-l_j|+2\mathbb{N}}$ for $i<j$.
We remark that the coefficients in the spectral decomposition \eqref{eq:spectral decomp} are the same as the ones in the case of the tensor product of finite-dimensional fundamental representations of $U_q(\widehat{\mathfrak{gl}}_m)$ %type $A_m^{(1)}$
for all sufficiently large $m$, under a suitable correspondence between fundamental representations. 
Hence the pairs $({\mb l},{\mb c})$ in \eqref{eq:fusion} are in  natural correspondence with those for the irreducible finite-dimensional representations of $U_q(\widehat{\mathfrak{gl}}_m)$.

To prove the formula \eqref{eq:spectral decomp}, we define more generally $q$-oscillator representations of quantum superalgebras of type $A$, and their affinizations over the quantum affine superalgebras. We use the generalized quantum group $\U(\e)$ of type $A$ introduced by Kuniba-Okado-Sergeev \cite{KOS}, which appears in the study of three-dimensional Yang-Baxter equation and which is also isomorphic to the quantum affine superalgebra under mild extension. Here $\e$ is a sequence of $0$'s and $1$'s corresponding to a non-conjugate simple root system associated to the generators of $\U(\e)$. We define the associated categories ${\mc O}_{{\rm osc},\e}$ and $\widehat{\mc O}_{{\rm osc},\e}$ in the same way as in ${\mc O}_{\rm osc}$ and $\widehat{\mc O}_{\rm osc}$.

Suppose that the length of $\e$ is sufficiently large, and let $\underline{\e}$ and $\ov{\e}$ be subsequences of $\e$ consisting of only $0$'s and $1$'s respectively. 
We note that $\widehat{\mc{O}}_{\rm osc,\underline{\e}}$ is the category $\widehat{\mc{O}}_{\rm osc}$ for $U_q(\widehat{\gl}_M)$-modules, and $\widehat{\mc{O}}_{\rm osc,\ov{\e}}$ is the category of finite-dimensional $U_q(\widehat{\gl}_N)$-modules, where $M$ and $N$ are the lengths of $\underline{\e}$ and $\ov{\e}$, respectively.
We define exact monoidal functors called truncation functors:
\begin{equation}\label{eq:super duality?}
\xymatrixcolsep{3pc}\xymatrixrowsep{0.3pc}\xymatrix{
 & \widehat{\mc{O}}_{\rm osc,\e}  \ar@{->}_{\mf{tr}^{\e}_{\underline{\e}}}[dl]\ar@{->}^{\mf{tr}^{\e}_{\ov{\e}}}[dr] &  \\
 \widehat{\mc{O}}_{\rm osc,\underline{\e}} & &  \widehat{\mc{O}}_{\rm osc,\ov{\e}}
}
\end{equation}
which preserve the normalized $R$ matrices and fundamental representations, and sends the irreducible representations in ${\mc O}_{{\rm osc},\e}$ to those in ${\mc O}_{{\rm osc},\underline{\e}}$ and ${\mc O}_{{\rm osc},\ov{\e}}$ or zero. 
This shows that the coefficients in the spectral decompositions are preserved whenever they are non-zero, and hence implies the formula \eqref{eq:spectral decomp}. 

Roughly speaking, a $q$-oscillator representation of $\U(\e)$ interpolates a $q$-oscillator representation of $U_q(\widehat{\gl}_M)$ and a finite-dimensional representation of  $U_q(\widehat{\gl}_N)$. The diagram \eqref{eq:super duality?} is reminiscent of the super duality between oscillator representations and integrable representations of classical types \cite{CW}, or the quantum super duality of affine type $A$ \cite{KL}, where an equivalence is shown to exist after taking a suitable limit of the associated categories. Our results provide some evidence for the existence of such an equivalence related to \eqref{eq:super duality?} (see also \cite{KO21}).

The paper is organized as follows. In Section \ref{sec:notations}, we introduce necessary notations and the fundamental type $q$-oscillator representations $\W_l$.
In Section \ref{sec:cat of q-osc}, we define the category $\mc{O}_{\rm osc}$ of $q$-oscillator representations of $U_q(\gl_n)$. In Section \ref{sec:cat of aff q-osc}, we define the category $\widehat{\mc{O}}_{\rm osc}$ and give the spectral decomposition of the normalized $R$ matrix \eqref{eq:spectral decomp}. We also compute the $\ell$-highest weight of the fundamental representation $\W_l$. In Section \ref{sec:irreducible of affine q-osc}, we construct irreducible representations in $\widehat{\mc{O}}_{\rm osc}$ including the Kirillov-Reshetikhin type modules and compute the characters in case of Kirillov-Reshetikhin type modules. In Section \ref{sec:spectral decomposition}, we study $q$-oscillator representations of the quantum affine superalgebras and prove the spectral decomposition of the normalized $R$ matrix in this case, which implies \eqref{eq:spectral decomp} under the truncation functor.
\smallskip 

{\bf Acknowledgement}
The authors would like to thank M. Okado for his interest in this work and letting us know the reference \cite{Ku18}, and the referees for very careful reading of the manuscript and helpful comments.

\section{Preliminary}\label{sec:notations}

\subsection{Notations}

We assume that $n$ is a positive integer such that $n\geq 4$.
We denote by $\mathbb{Z}_+$ (resp. $\mathbb{N}$) the set of non-negative (resp. positive) integers. Let $q$ be an indeterminate and $$[m]=\frac{q^m-q^{-m}}{q-q^{-1}}\quad (m\in \mathbb{Z}_+).$$
Let $\Bbbk=\mathbb{Q}(q)$ and $\Bbbk^{\times}=\Bbbk \setminus \left\{0\right\}$.
We assume the following notations throughout the paper except in Section \ref{sec:spectral decomposition}: 
\begin{itemize}

\item[$\bullet$] $I=\{\,0,1,\ldots,n-1\,\}$,

\item[$\bullet$] $r\in I\setminus\{0,1,n-1\}$,

\item[$\bullet$] $P = \Z\Lambda_r \oplus \Z\de_1\oplus\cdots\oplus \Z\de_n$,

\item[$\bullet$] $P^\vee={\rm Hom}_\mathbb{Z}(P,\mathbb{Z})=\Z c\oplus \Z\de^\vee_1\oplus\cdots\oplus\Z\de^\vee_n$ such that the pairing $\langle\,,\, \rangle$ on $P\times P^\vee$ is given by
$$
\lang \de_i, \de^\vee_j \rang =\de_{ij},\quad 
\lang \de_i,c \rang=\lang \La_r,\de^\vee_i \rang=0,\quad
\lang \La_r,c\rang =1,
$$

\item[$\bullet$] $\alpha_i=\de_i-\de_{i+1}\in P$, 
$\alpha_i^\vee = \de^\vee_i-\de^\vee_{i+1} + (\de_{ir}-\de_{i0})c \in P^\vee$  ($i\in I$),

\end{itemize}
where we understand the subscript $i\in I$ modulo $n$.

Let $U_q(\widehat{\mf{gl}}_n)$ be the quantum affine algebra associated to $\mf{gl}_n$, which is the associative $\Bbbk$-algebra with $1$ 
generated by $q^h, e_i, f_i$ for $h\in P^\vee$ and $i\in I$ 
satisfying
\begin{gather*}
q^0=1, \quad q^{h +h'}=q^{h}q^{h'} \quad (h, h' \in P^\vee),\\  
q^h e_i q^{-h}=q^{\lang\alpha_i,h\rang}e_i,\quad 
q^h f_i q^{-h}=q^{-\lang\alpha_i,h\rang}f_i \\ 
e_if_j - f_je_i =\delta_{ij}\frac{k_{i} - k_{i}^{-1}}{q-q^{-1}} \\
e_i e_j -  e_j e_i = f_i f_j -  f_j f_i =0
 \quad \text{($i-j\not\equiv \pm 1\!\!\pmod n$)},\\ 
\begin{array}{ll}
e_i^2 e_j- [2] e_i e_j e_i + e_j e_i^2= 0,\ \ 
f_i^2 f_j- [2] f_i f_j f_i+f_j f_i^2= 0
\end{array}
\quad \text{($i-j\equiv \pm 1\!\!\pmod n$)}, 
\end{gather*} 
where $k_i=q^{\alpha^\vee_i}$ for $i\in I$.
There is a Hopf algebra structure on $U_q(\agl_n)$, where the comultiplication $\Delta$ is given by 
\begin{equation*}
\begin{array}{l}
\Delta(q^h)=q^h\otimes q^h, \\ 
\Delta(e_i)= 1\ot e_i + e_i\ot k_i^{-1}, \\
\Delta(f_i)= f_i\ot 1 + k_i\ot f_i , \\  
\end{array}\quad (h\in P^\vee, i\in I),
\end{equation*}
and the antipode $S$ is given by 
\begin{equation*}
S(q^h)=q^{-h}, \quad  S(e_i)=-e_i k_i, \quad S(f_i)=-k_i^{-1} f_i.
\end{equation*}

We let $U_q(\gl_n)$ be the subalgebra of $U_q(\agl_n)$ generated by
$q^h, e_i, f_i$ for $h\in P^\vee$ and $i\in I\setminus\{0\}$.
%Let $\mf{l}$ be the maximal Levi subalgebra of $\gl_n$ corresponding to the set of simple roots $\{\,\alpha_i\,|\,i\in J\,\}$.
Let $U_q(\mf{l})$ be the subalgebra of $U_q(\gl_n)$ generated by $q^h, e_i, f_i$ for $h\in P^\vee$ and $i\in I\setminus\{0,r\}$, where $\mf{l}$ is the maximal Levi subalgebra of $\gl_n$ corresponding to the set of simple roots $\{\,\alpha_i\,|\,i\in I\setminus\{0,r\}\,\}$.

Let $V$ be a $U_q(\agl_n)$-module. For $\la\in P$, we define the $\la$-weight space of $V$ to be
\begin{equation*}
V_\la 
= \{\,u\in V\,|\,q^h u= q^{\lang\la,h\rang} u \ \ (h\in P^\vee) \,\},
\end{equation*}
and write ${\rm wt}(u)=\la$ for $u\in V_\la$.

\begin{rem}\label{rem:convention for P}
{\rm We may regard $P$ as the weight lattice for $\agl_n$ modulo the imaginary root $\de$ of type $A_{n-1}^{(1)}$. We remark that we take $P^\vee$ as an extension of the coweight lattice $\bigoplus_{i=1}^n\Z\de_i^\vee$ for $\gl_n$ with respect to the simple root $\alpha_r$ unlike the usual convention (with respect to $\alpha_0$). This convention enables us to describe the main results in this paper more easily.
}
\end{rem}

\subsection{Fundamental $q$-oscillator representation $\W_l$}

Let us introduce a $U_q(\agl_n)$-module, which plays an important role in this paper. 
Let 
$$
\mathcal{W}=\bigoplus_{\mathbf{m}\in\mathbb{Z}_{+}^{n}}\Bbbk\ket{\mathbf{m}}.
$$ 
Let $\{\be_1,\dots,\be_n\}$ denote the standard basis of $\Z^n$ and ${\bf 0}=(0,\dots,0)$.

\begin{prop}\label{prop:osc module-1}
For $x\in\Bbbk^\times$, the following formula defines a $U_q(\agl_n)$-module structure on $\mathcal{W}$, which we denote by $\W(x)$:
{\allowdisplaybreaks
\begin{align*}
&
\begin{array}{l}
q^c \ket{\bf m}=q^{-1}\ket{\bf m}, \\
q^{\de_i^\vee}\ket{\bf m}=q^{-m_{i}}\ket{\bf m} \quad (1\le i\le r),\\
q^{\de_j^\vee}\ket{\bf m}=q^{m_{j}}\ket{\bf m} \quad (r+1\le j\le n),\\
\end{array}\\
&
\begin{array}{l}
e_{0}\ket{\bf m} =x\ket{{\bf m}+\be_{1}+\be_{n}},\\
f_{0}\ket{{\bf m}} =-x^{-1}[m_{1}][m_{n}]\ket{{\bf m}-\be_{1}-\be_{n}},\\
\end{array}\\
&
\begin{array}{l}
e_{r}\ket{\bf m} =-[m_{r}][m_{r+1}]\ket{{\bf m}-\be_{r}-\be_{r+1}},\\
f_{r}\ket{\bf m} =\ket{{\bf m}+\be_{r}+\be_{r+1}},\\
\end{array}\\
&
\begin{array}{l}
e_{i}\ket{\bf m} =[m_{i}]\ket{{\bf m}-\be_{i}+\be_{i+1}}\\
f_{i}\ket{\bf m} =[m_{i+1}]\ket{{\bf m}+\be_{i}-\be_{i+1}}\\
\end{array}\quad(1\leq i< r), \\
&
\begin{array}{l}
e_{j}\ket{\bf m} =[m_{j+1}]\ket{{\bf m}+\be_{j}-\be_{j+1}}\\
f_{j}\ket{\bf m} =[m_{j}]\ket{{\bf m}-\be_{j}+\be_{j+1}}\\
\end{array}\quad (r< j\le n-1),
\end{align*}
where ${\bf m}=(m_1,\dots,m_n)\in\Z_+^n$ and we understand $\ket{\bf m'}=0$ unless $\ket{\bf m'}$ on the right-hand side in the above formula belongs to $\Z_+^n$.} 
\end{prop}
\pf We can check directly that the operators defined above satisfy the defining relations of $U_q(\agl_n)$.
\qed\vskip 2mm

Note that $\W$ decomposes into a direct sum of weight spaces, where we have for ${\bf m}=(m_1,\dots,m_n)\in \Z_+^n$
\begin{equation}\label{eq:weight of m}
{\rm wt}(\ket{\bf m})= 
-\La_r - \sum_{i=1}^rm_i\delta_i  + \sum_{j=r+1}^nm_j\delta_j.
\end{equation}
For ${\bf m}\in \Z_+^n$, put 
\begin{equation}\label{eq:notation for state vector}
\begin{split}
|{\bf m}|_-&=m_1+\dots + m_r,\quad |{\bf m}|_+=m_{r+1}+\dots + m_n,\\
|{\bf m}| &= |{\bf m}|_++|{\bf m}|_-,\ \
l( {\bf m} ) =|{\bf m}|_+-|{\bf m}|_-.
\end{split}
\end{equation}
For $l\in \Z$, let
\begin{equation*}
\W_l = \bigoplus_{l( {\bf m} )=l}\Bbbk\ket{\mathbf{m}},
\end{equation*}
and let $v_l\in \W_l$ given by
\begin{equation}\label{eq:v_l}
v_l=
\begin{cases}
\ket{l\be_{r+1}} & \text{if $l\ge 0$},\\
\ket{-l\be_{r}} & \text{if $l< 0$}.
\end{cases}
\end{equation}
It is clear from definition that $\W_l$ is invariant under the action of $U_q(\agl_n)$, and hence becomes a $U_q(\agl_n)$-submodule, which we denote by $\W_l(x)$.

\begin{prop}\label{prop:osc module-2}
For $l\in \Z$ and $x\in\Bbbk^\times$, $\W_l(x)$ is an irreducible $U_q(\agl_n)$-module and $\W(x)$ decomposes as a $U_q(\agl_n)$-module as follows:
\begin{equation*}
\W(x) = \bigoplus_{l\in\Z}\W_l(x).
\end{equation*}
Moreover, $\W_l(x)$ is an irreducible highest weight $U_q(\gl_n)$-module with highest weight vector $v_l$.  
\end{prop}
\pf One can check directly that $e_iv_l=0$ for all $i\in I\setminus\{0\}$ and $\W_l(x)$ is generated by $v_l$ as a $U_q(\mf{gl}_n)$-module. 
Moreover, given any non-zero vector $v$, we may obtain $v_l$ (up to scalar multiplication) by applying $e_{i_1}\dots e_{i_t}$ to $v$ for some $i_1,\dots,i_t\in I\setminus\{0\}$. This proves the assertion.
\qed\vskip 2mm

Let us call $\W_l(x)$ the {\em $l$-th fundamental $q$-oscillator representation of $U_q(\mf{gl}_n)$}.

\section{The semisimple tensor category $\mc{O}_{\rm osc}$ of $q$-oscillator $U_q(\gl_n)$-modules}\label{sec:cat of q-osc}
 
\subsection{Oscillator representations of $\gl_n$}
We take
\begin{itemize}

\item[$\bullet$] $P_{\rm fin}^\vee= \bigoplus_{i=1}^n\Z\ude^\vee_i\subset P^\vee$,  where $\ude^\vee_n=\de^\vee_n$ and $\ude^\vee_{i}=\ude^\vee_{i+1}+\alpha^\vee_i$ for $1\le i\le n-1$,

\item[$\bullet$] $P_{\rm fin}=P/\Z(\La_r-\de_1-\dots-\de_r)$,

\end{itemize}
as the coweight and weight lattices for $\gl_n$, respectively. We denote by $\ude_i$ the image of $\de_i$ in $P_{\rm fin}$ so that $\lang \ude_i,\ude^\vee_j \rang=\de_{ij}$. The image of $\La_r$ in $P_{\rm fin}$ is the $r$-th fundamental weight for $\gl_n$, which we denote by $\varpi_r$. For $\varpi\in P_{\rm fin}$, let $V(\varpi)$ be the irreducible highest weight representation of $\gl_n$ with highest weight $\varpi$.

For a positive integer $\ell$, let 
\begin{equation*}
\mc{P}(GL_\ell) =\{\,\la=(\la_1,\dots,\la_\ell)\in \Z^\ell\,|\, \la_1\ge \dots\ge \la_\ell\,\} 
\end{equation*}
be the set of generalized partitions of length $\ell$, while a partition is a finite non-increasing sequence of positive integers. It is well-known that $\mc{P}(GL_\ell)$ parametrizes the finite-dimensional irreducible representations of $GL_\ell=GL_\ell(\C)$ with integral weights, say $V_{GL_\ell}(\la)$ for $\la\in \mc{P}(GL_\ell)$.
We may identify $\la$ with a generalized Young diagram, for example, when $\la=(4,2,2,0,0,-1,-3)\in \mc{P}(GL_7)$, we have
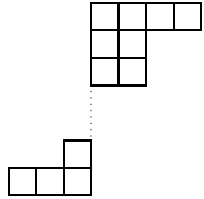
\[ \begin{tikzpicture}
\node [above right] at (0, 0) 
{$ 
\ytableausetup{mathmode, boxsize=1em} 
\begin{ytableau}
\none & \none &  \none & {} & {} & {} & {} \\
\none & \none &  \none & {} & {} & \none & \none \\
\none & \none &  \none & {} & {} & \none & \none \\
\none & \none & \none  & \none & \none & \none & \none \\
\none & \none &  \none & \none & \none & \none & \none \\
\none & \none &  {} & \none & \none & \none & \none \\
{} & {} &  {} & \none & \none & \none & \none \\
\end{ytableau}
$};
\draw [dotted] (1.225, 2.36em) -- (1.225, 4.6em);
\end{tikzpicture} \]

Let $\la\in \mc{P}(GL_\ell)$ be given. When $\lambda$ can be written as
\begin{equation*}\label{eq:s and t}
\la=(\la_1\ge \dots \ge \la_s>\la_{s+1}=\dots =\la_{t}=0>\la_{t+1}\ge \dots\ge \la_\ell) 
\end{equation*} 
for some $s < t$, we put
\begin{equation}\label{eq:lambda pm}
\la^+=(\la_1\ge \dots \ge \la_s),\quad \la^-= (-\la_\ell \ge \dots\ge -\la_{t+1}).
\end{equation}
Note that $s$ (resp. $t$) does not exist if $\la_i<0$ (resp. $\la_i>0$) for all $1\le i\le n$.

Let  
\begin{equation*}
\mc{P}(GL_\ell)_{(r,n-r)}=\{\,\la\in \mc{P}(GL_\ell)\,|\, \text{$\ell(\la^-)\le r$, $\ell(\la^+)\le n-r$}\,\},
\end{equation*}
where $\ell(\mu)$ denotes the length of a partition $\mu$, and put
\begin{equation*}
\mc{P}(GL)_{(r,n-r)} = \bigsqcup_{\ell\ge 1}\mc{P}(GL_\ell)_{(r,n-r)}.
\end{equation*}

For $\la\in \mc{P}(GL_\ell)_{(r,n-r)}$, we define the weight $\varpi_\la\in P_{\rm fin}$ by
\begin{equation}\label{eq:Lambda highest weight}
\begin{split}
\varpi_{\la}
&= -\ell \varpi_r +\sum_{i=1}^s\la_i\ude_{r+i} + \sum_{j=t+1}^\ell\la_j\ude_{r-\ell +j}\\
&= -\ell \varpi_r +\la_1\ude_{r+1}+\dots+\la_s\ude_{r+s} + \la_{t+1}\ude_{r-\ell+t+1} +\dots +\la_\ell\ude_r.
\end{split}
\end{equation}

\begin{rem}\label{rem:comb rule for h.w.}
{\rm 
There is a simple combinatorial rule to describe $\varpi_\la$. 
Let us first fill $\la^+$ with $r+1, r+2, \dots$ row by row from the top, and then fill $\la^-$ with $r, r-1, \dots$ row by row from the top (or fill each row corresponding to the negative parts of $\la$ with $r, r-1, \dots$ from the bottom). Then each coefficient of $\ude_i$ (up to sign) is given by counting the number of occurrences of $i$'s in the diagram.

Let $\la\in \mc{P}(GL_7)$ be as above with $s=3$ and $t=5$. Suppose that $n=8$ and $r=3$ so that $\la\in \mc{P}(GL_7)_{(3,5)}$. Then we have
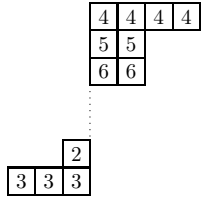
\[ \begin{tikzpicture}
\node [above right] at (0, 0) 
{$ 
\ytableausetup{mathmode, boxsize=1em} 
\begin{ytableau}
\none & \none &  \none & {\tl 4} & {\tl 4} & {\tl 4} & {\tl 4} \\
\none & \none &  \none & {\tl 5} & {\tl 5} & \none & \none \\
\none & \none &  \none & {\tl 6} & {\tl 6} & \none & \none \\
\none & \none & \none  & \none & \none & \none & \none \\
\none & \none &  \none & \none & \none & \none & \none \\
\none & \none & {\tl 2} & \none & \none & \none & \none \\
{\tl 3} & {\tl 3} &  {\tl 3} &  \none & \none & \none & \none \\
\end{ytableau}
$};
\draw [dotted] (1.225, 2.36em) -- (1.225, 4.6em);
\end{tikzpicture} \]
and hence
\[
\varpi_\la = -7\varpi_3 + 4\ude_4 +2\ude_5 +2\ude_6 - \ude_2- 3\ude_3.
\] }
\end{rem}

Let  
$$V^{\la} = V(\varpi_\la)$$
denote the irreducible highest weight $\gl_n$-module with highest weight $\varpi_\la$. For example, we have $V^{(l)}=V(\varpi_{(l)})$ for $l\in \Z$, where 
\begin{equation*}
 \varpi_{(l)}=
\begin{cases}
 -\varpi_r + l\ude_{r+1} & \text{if $l\ge 0$},\\
 -\varpi_r - l\ude_{r} & \text{if $l< 0$}.
\end{cases}
\end{equation*}

Note that $V^\la$ is not finite-dimensional. Instead, $V^\la$ belongs to the parabolic BGG category for $\gl_n$ with respect to the maximal Levi sublagebra $\mf{l}$ since $\varpi_\la$ is $\mf{l}$-dominant, that is, $\lang \varpi_\la,\alpha_i^\vee\rang\in\Z_+$ for $i\in I\setminus\{0,r\}$, but not $\gl_n$-dominant since $\lang \varpi_\la, \alpha^\vee_r\rang <0$.

The irreducible representation $V^\la$ is of special importance, which has been studied in connection with unitary representations of Lie groups.
They also satisfy a nice duality with the finite-dimensional representations of $GL_\ell$ as follows (see \cite{EHW,H,HTW,KV} and references therein for more details).
\begin{thm}[{\cite{H}}]\label{thm:Howe duality}
Let $W=\bigoplus_{l\in\Z}V^{(l)}$.
There exists an action of $GL_\ell$ on $W^{\ot \ell}$ $(\ell\ge 1)$ which commutes with the action of $\mathfrak{gl}_n$ and gives the following multiplicity-free decomposition as a $(\gl_n,GL_\ell)$-module:
\begin{equation*}
W^{\ot \ell} = \bigoplus_{\la\in \mc{P}(GL_\ell)_{(r,n-r)}}V^\la\ot V_{GL_\ell}(\la).
\end{equation*}
\end{thm}\vskip 2mm

Suppose that $\mu\in \mc{P}(GL_{\ell})_{(r,n-r)}$ and $\nu\in \mc{P}(GL_{\ell'})_{(r,n-r)}$ are given. Then $V^\mu\ot V^\nu$ is semisimple  since it is a $\mathfrak{gl}_n$-submodule of $W^{\ot{(\ell+\ell')}}$ by Theorem \ref{thm:Howe duality}, and  
\begin{equation}\label{eq:LR for osc}
V^\mu\ot V^\nu = \bigoplus_{\la\in \mc{P}(GL_{\ell+\ell'})_{(r,n-r)}}(V^\la)^{\oplus c_{\mu\nu}^\la},
\end{equation}  
where  $c_{\mu\nu}^\la$ is given by the branching multiplicity 
\begin{equation}\label{eq:LR}
\dim{\rm Hom}_{GL_\ell\times GL_{\ell'}}\left( V_{GL_\ell}(\mu)\ot V_{GL_{\ell'}}(\nu), V_{GL_{\ell+\ell'}}(\la) \right).
\end{equation}
The multiplicity \eqref{eq:LR} is equal to the usual Littlewood-Richardson coefficients, while $V^\mu\ot V^\nu$ has infinitely many distinct irreducible components.  

\begin{ex}{\rm
For $l,m\in\Z$, we have
\begin{equation}\label{eq:LR for fundamentals}
V^{(l)}\ot V^{(m)} = 
V^{(l_1,l_2)}\oplus V^{(l_1+1,l_2-1)}\oplus \dots = \bigoplus_{t\in \Z_+}V^{(l_1+t,l_2-t)},
\end{equation}
where $l_1=\max\{l,m\}$ and $l_2=\min\{l,m\}$.
}
\end{ex}

\begin{df}\label{def:cat-osc-cl}
{\rm
Let ${O}_{\rm osc}$ be the category of $\gl_n$-modules 
$V$ such that
\begin{itemize}
 \item[(1)] $V=\bigoplus_{\mu\in P_{\rm fin}}V_\mu$, where $\dim V_\mu<\infty$, and the set of weights of $V$ is finitely dominated,

 \item[(2)] $V=\bigoplus_{\ell\ge 1}V_{\ell}$, where $V_\ell$ is a direct sum of $V^\la$'s for $\la\in \mc{P}(GL_\ell)_{(r,n-r)}$, and $V_\ell=0$ for all sufficiently large $\ell$. 
\end{itemize}
} 
\end{df}

The category ${O}_{\rm osc}$ is a full subcategory of (a version of) the parabolic BGG category of $\gl_n$ with respect to $\mf l$, where an object $V$ has a weight space decomposition with finite-dimensional weight spaces, the set of weights is finitely dominated, and it decomposes into irreducible finite-dimensional representations of $\mf{l}$. It is closed under tensor product by \eqref{eq:LR for osc}.

\begin{prop}
The category ${O}_{\rm osc}$ is a semisimple tensor category. 
\end{prop} 
 
\begin{rem}{\rm
The $\mathfrak{gl}_n$-modules satisfying the conditions in Definition~\ref{def:cat-osc-cl} except the finite dominance
also form a semisimple tensor category.
%The set of weights of a $\gl_n$-module $V$ satisfying the condition (2) only is not necessarily finitely dominated, but the category of such $\gl_n$-modules $V$ is still a semisimple tensor category.
} 
\end{rem}

\subsection{The character of $V^\la$}
Let us recall the character of ${V}^\la$ for $\la\in \mc{P}(GL_\ell)_{(r,n-r)}$.
Let $\Z[P_{\rm fin}]$ be the group algebra of $P_{\rm fin}$ with basis $\{\,e^\mu\,|\,\mu\in P_{\rm fin}\,\}$. Put 
\begin{equation*}
 t=e^{-\varpi_r},\quad 
 x_i=e^{\ude_{r+i}} \quad (1\le i\le n-r),\quad 
 y_j=e^{-\ude_{j}} \quad (1\le j\le r).
\end{equation*}
For partitions $\mu$ and $\nu$ such that $\ell(\mu)=s\le n-r$ and $\ell(\nu)=t\le r$, let $s_\mu(x)$ and $s_\nu(y)$ be the Schur polynomials in $x_i$'s and $y_j$'s corresponding to $\mu$ and $\nu$, respectively. We may also define the Schur polynomials in $x_i$'s and $y_j$'s associated to skew Young diagrams.

For $\la,\eta,\xi\in \mc{P}(GL_\ell)$,
let
\begin{equation*}\label{eq:LR*}
\hat{c}^\la_{\eta\xi^*}=\dim {\rm Hom}_{GL_\ell}\left(V_{GL_{\ell}}(\la), V_{GL_\ell}(\eta)\ot V_{GL_{\ell}}(\xi^*)\right),
\end{equation*}
where $\xi^*=(-\xi_\ell,\dots,-\xi_{1})\in \mc{P}(GL_\ell)$. 
Then we have the following by \cite[Proposition 3.14]{K08} (by putting $\mc{A}=\{\,r+1,\dots,n\,\}$ and $\mc{B}=\{\,1,\dots,r\,\}$ in loc.~cit.).
\begin{prop}[{\cite[Theorem 5.3]{CLZ}}]\label{prop:ch of classical osc}
 For $\la\in \mc{P}(GL_\ell)_{(r,n-r)}$, we have
\begin{equation*}
\begin{split}
 {\rm ch}{V}^\la
 &=t^\ell\sum_{\mu,\nu}\hat{c}^\la_{\mu\nu^*}s_\mu(x)s_\nu(y)\\
 &=t^\ell\sum_{\eta,d}s_{(\la+(d^\ell))/\eta}(x)s_{(d^\ell)/\eta}(y),
\end{split}
\end{equation*}
where the sum is over the partitions $\mu$ and $\nu$ such that $\ell(\mu)\le \min\{n-r,\ell\}$ and $\ell(\nu)\le \min\{r,\ell\}$, or over the partitions $\eta$ and $d\in\Z_+$ such that $\la+(d^\ell)$ is a partition and $\eta\subset \la+(d^\ell)$, $\eta\subset (d^\ell)$ as Young diagrams.
Here we regard $\mu=(\mu_1,\dots,\mu_{\ell(\mu)},0,\dots,0)$ and $\nu^*=(0,\dots,0,-\nu_{\ell(\nu)},\dots,-\nu_1)$.
\end{prop}

Let $\mu=(\mu_1,\dots,\mu_s)$ and $\nu=(\nu_1,\dots,\nu_t)$ be partitions with $\ell(\mu)=s\le n-r$ and $\ell(\nu)=t\le r$. For $\ell \ge n$, let
\begin{equation}\label{eq:generalized partition}
[\mu,\nu]_\ell=
(\mu_1\ge \dots \ge \mu_s>\underbrace{0=\dots =0}_{\ell-s-t}>-\nu_{t}\ge \dots\ge -\nu_1)\in \mc{P}(GL_\ell)_{(r,n-r)}.
\end{equation}

\begin{cor}[{\cite[Theorem 3.2]{HTW}}]\label{cor:ch of classical osc}
 Under the above hypothesis, if $\ell\ge n$, then we have
\begin{equation*}
 {\rm ch}V^{[\mu,\nu]_\ell}
 =t^\ell s_\mu(x)s_{\nu}(y)\dfrac{1}{\prod_{i,j}(1-x_iy_j)}.
\end{equation*}
Hence the normalized characters $t^{-\ell}{\rm ch}V^{[\mu,\nu]_\ell}$ coincide for all $\ell\ge n$.
\end{cor}
\pf It follows from the second formula in Proposition \ref{prop:ch of classical osc} and the skew Cauchy identity \cite{SS}.
\qed

\begin{rem}{\rm
 The stabilized limit of the normalized character $t^{-\ell}{\rm ch}V^{[\mu,\nu]_\ell}$ is equal to the character of the parabolic Verma module with highest weight $\varpi_{[\mu,\nu]_\ell}+\ell\varpi_r$  %$\varpi_\la+\ell\varpi_r$
 with respect to $\mf{l}$ (cf.~\eqref{eq:Lambda highest weight}).
} 
\end{rem}

\subsection{$q$-oscillator representations of $U_q(\gl_n)$ and the category ${\mc O}_{\rm osc}$} 
In this subsection, we introduce a $q$-analogue of ${O}_{\rm osc}$.

Let us first prove that $\W_l(x)$ is a $q$-analogue of $V^{(l)}$ as a $U_q(\gl_n)$-module for $l\in\Z$  and $x\in\Bbbk^\times$, and $\W(x)^{\ot \ell}$ is semisimple for $\ell\ge 1$.
Let ${\bf A}=\Z[q,q^{-1}]$. For $x\in\Bbbk^\times$, let 
\begin{equation*}
\W(x)_{{\bf A}} =\bigoplus_{\mathbf{m}\in\mathbb{Z}_{+}^{n}}{\bf A}\ket{\mathbf{m}}.
\end{equation*}
The ${\bf A}$-module $\W(x)_{{\bf A}}$ is invariant under $e_i$, $f_i$, $q^h$ and $\{q^h\}:=\frac{q^h-q^{-h}}{q-q^{-1}}$ for $i\in I\setminus\{0\}$ and $h\in P^\vee_{\rm fin}$, and hence the $\C$-vector space 
\begin{equation*}
\ov{\W(x)}=\W(x)_{{\bf A}}\otimes_{{\bf A}}\mathbb{\C}
\end{equation*}
is invariant under the $\mathbb{C}$-linear endomorphisms ${\rm E}_i$, ${\rm F}_i$ and ${\rm D}_a$ induced from $e_i$, $f_i$ and $\{q^{\ude_a}\}$, respectively for $i\in I\setminus\{0\}$ and $1\le a\le n$. Here $\C$ is an ${\bf A}$-module such that $f(q)\cdot z = f(1)z$ for $f(q)\in {\bf A}$ and $z\in \C$.
The endomorphisms ${\rm E}_i$, ${\rm F}_i$ and ${\rm D}_a$ satisfy the defining relations for the  enveloping algebra $U(\gl_n)$ (cf.~\cite[Chapter 5]{Ja}). 
Hence $\ov{\W(x)}$ becomes a $U(\gl_n)$-module and so does $\ov{\W_l(x)}:=(\W_l(x)\cap \W(x)_{\bf A})\ot_{\bf A}\C\subset \ov{\W(x)}$ for $l\in\Z$.
Similarly, we can check that for $\ell\ge 1$, $\ov{\W(x)^{\ot\ell}}:=(\W(x)_{\bf A})^{\ot\ell}\ot_{\bf A}\C$ is a $U(\gl_n)$-module and it is isomorphic to $\ov{\W(x)}^{\ot\ell}$.

\begin{lem}\label{lem:classical limit-1}
Let $x\in \Bbbk^\times$ be given. We have the following as $\gl_n$-modules:
\begin{itemize}
\item[(1)] $\ov{\W_l(x)} \cong V^{(l)}$ for $l\in\Z$, and hence $\ov{\W(x)}\cong W$,

\item[(2)] $\ov{\W(x)^{\ot\ell}}\cong W^{\ot \ell}$ for $\ell\ge 1$.
\end{itemize}
\end{lem}
\pf (1) We may assume that $l=0$ since the proof for the other cases is almost the same. It is clear that ${\rm E}_i(v_0\ot 1)={\rm E}_i(\ket{{\bf 0}} \ot 1)=0$ for all $i\in I\setminus\{0\}$. 
Since
\begin{equation*}
\begin{split}
{\rm D}_a(|{\bf 0}\rangle \ot 1) 
& = \left(\frac{q^{\ude_a}-q^{-\ude_a}}{q-q^{-1}}|{\bf 0}\rangle\right) \ot 1  = \left(\frac{q^{-1}-q}{q-q^{-1}}|{\bf 0}\rangle\right) \ot 1=- |{\bf 0}\rangle \ot 1,
\end{split}
\end{equation*}
for $1\le a\le r$ and ${\rm D}_a(|{\bf 0}\rangle \ot 1)= |{\bf 0}\rangle \ot 1$ otherwise, the vector $\ket{{\bf 0}} \ot 1$ is a highest weight vector with highest weight $-\varpi_r$.

We see from the actions of ${\rm E}_i$ for $i\in I\setminus\{0\}$ that any submodule of $\ov{\W_0(x)}$ contains $|{\bf 0}\rangle \ot 1$. This implies that $\ov{\W_0({x})}$ is an irreducible highest weight module with highest weight $-\varpi_r$, that is, $V^{(0)}$. 
(2) It follows from (1).
\qed \newline

For $\gamma\in P$, let $V_q(\gamma)$ be the irreducible highest weight $U_q(\gl_n)$-module with highest weight $\gamma$. 
For $\la\in \mc{P}(GL_\ell)_{(r,n-r)}$, let 
\begin{equation}\label{eq:Lambda highest weight-2}
\La_{\la}
= -\ell \La_r +\sum_{i=1}^s\la_i\de_{r+i} + \sum_{j=t+1}^\ell\la_j\de_{r-\ell +j}\in P,
\end{equation}
where $s=\ell(\la^+)$ and $t=\ell(\la^-)$ (cf.~\eqref{eq:lambda pm}), and let 
\begin{equation*}
\mc{V}^{\la}=V_q(\La_\la).
\end{equation*}
Note that the image of $\La_\la$ in $P_{\rm fin}$ is equal to $\varpi_\la$.

\begin{prop}\label{prop:hwsub-W^l}
Let $\ell\ge 1$ and $x\in \Bbbk^\times$ be given. Any highest weight $U_q(\gl_n)$-submodule in $\W(x)^{\ot \ell}$ is isomorphic to $\mc{V}^{\la}$ for some $\la\in \mc{P}(GL_\ell)_{(r,n-r)}$.
\end{prop}

%\todo{It seems that the semisimplicity of $\mathcal{W}^{\otimes \ell}$ is not explicitly referred to except in Prop.3.14, so I just remove that part and briefly mention in the proof of Prop.3.14. I also relabel this proposition and accordingly Prop.6.15. (Lee)}

\pf We may write $\W=\W(x)$ and $\W_l=\W_l(x)$ since $x$ is irrelevant for the action of $U_q(\gl_n)$. 

Let $V$ be a submodule of $\W^{\ot \ell}$, which is generated by a highest weight vector $v$ with highest weight $\La$. We may assume that $v \in \W_{\bf A}^{\ot \ell}$ by multiplying a scalar.

Let $V_{\bf A}$ be the ${\bf A}$-span of $f_{i_1}\ldots f_{i_s}v$ for $s\geq 0$ and $i_1,\ldots,i_s\in I\setminus\{0\}$.  
We have $V_{{\bf A}}=\bigoplus_{\mu\in P_{\rm fin}}V_{\mu,{\bf A}}$, where $V_{\mu,{\bf A}}=V_{\bf A}\cap V_\mu$ and ${\rm rank}_{\bf A}V_{\mu,{\bf A}}=\dim_{\Bbbk}V_\mu$ (cf.~\cite[Chapter 5]{Ja}).
We can check that $V_{\bf A}$ is invariant under $e_i$, $f_i$, $q^h$ and $\{q^h\}$ for $i\in I\setminus\{0\}$ and $h\in P^\vee_{\rm fin}$, and hence 
\begin{equation}\label{eq:classical limit}
\ov{V}:=V_{\bf A}\ot_{\bf A}\C,
\end{equation}
becomes a $U(\gl_n)$-module, where $\ov{V}=\bigoplus_{\mu}\ov{V}_\mu$ with $\dim_{\mathbb{C}}\ov{V}_\mu={\rm rank}_{\bf A}V_{\mu,{\bf A}}$.

Since $\ov{V}$ is a $U(\gl_n)$-submodule of $W^{\ot\ell}$ by Lemma \ref{lem:classical limit-1} and it is also a highest weight module generated by $v\ot 1$ with highest weight $\La$, it follows from Theorem \ref{thm:Howe duality} that $\ov{V}$ is isomorphic to $V^\la$ for some $\la\in \mc{P}(GL_\ell)_{(r,n-r)}$. This implies that $V$ is irreducible and  $V\cong \mc{V}^\la$.\qed
%(2) 
%Let $v$ be a non-zero weight vector of $\W^{\ot \ell}$ and $V=U_q(\gl_n)v$.
%Since $\W=\bigoplus_{l\in\Z}\mc{V}^{(l)}$ and $\mc{V}^{(l)}$ is a highest weight $U_q(\gl_n)$-module for $l\in\Z$, ${\rm wt}(v)$ belongs to $\sum_{i=1}^t\La_{(l_k)}-\sum_{i\in I\setminus\{0\}}\Z_+\alpha_i$ for some $l_1,\dots,l_t\in \Z$. Hence the set of weights of $V$ is finitely dominated. In particular, there exists only finitely many weights of $\mathcal{W}^{\otimes \ell}$ which is greater than $\mathrm{wt}(v)$.

%Now let $U$ be the sum of all irreducible submodules of $\mathcal{W}^{\otimes \ell}$ and suppose $U$ is proper. Then we can take a non-zero maximal weight vector $u\in V\setminus U$. By (1), $u$ generates an irreducible submodule not intersecting with $U$, which is absurd. Hence $\mathcal{W}^{\otimes \ell}=U$ and is semisimple.\qed

We also have seen in the proof of Proposition \ref{prop:hwsub-W^l} that
\begin{cor}\label{cor:classical limit of V lambda}
For $\la\in \mc{P}(GL_\ell)_{(r,n-r)}$, we have as a $\gl_n$-module
\begin{equation*}
\ov{\mc{V}^{\la}}\cong V^\la.
\end{equation*}
\end{cor}

%\begin{rem}{\rm
%Another proof of Proposition \ref{prop:semisimplicity of osc} using polarization will be given in Section \ref{sec:spectral decomposition}. 
%} 
%\end{rem}

\begin{df}\label{def:cat osc}
{\rm Let $\mc{O}_{\rm osc}$ be the category of $U_q(\gl_n)$-modules 
$V$ such that 
\begin{itemize}
 \item[(1)] $V=\bigoplus_{\mu\in P_{\rm fin}}V_\mu$, where $\dim V_\mu<\infty$, and the set of weights of $V$ is finitely dominated,

 \item[(2)] $V=\bigoplus_{\ell\ge 1}V_{\ell}$, where $V_\ell$ is a direct sum of $\V^\la$'s for $\la\in \mc{P}(GL_\ell)_{(r,n-r)}$ with finite multiplicity for each $\la$, and $V_\ell=0$ for all sufficiently large $\ell$.
 
\end{itemize}
}
\end{df}

\begin{prop}\label{prop:semisimplicity}
The category $\mc{O}_{\rm osc}$ is closed under tensor product, and hence it is a semisimple tensor category.
\end{prop}
\pf {It follows directly from the semisimplicity of $\mathcal{W}^{\otimes \ell}$, which will be proved later in Lemma~\ref{lem:pol and semisimple} in a more general setting.}
\qed

\begin{rem}\label{rem:LR for q osc}
{\rm The decomposition of the tensor product $\mc{V}^\la\ot \mc{V}^\mu$ is the same as in the case of $V^\la \ot V^\mu$ \eqref{eq:LR for osc} since the characters of $\mc{V}^\la$ and $V^\la$ are the same.
}
\end{rem}

\section{Category $\widehat{\mc O}_{\rm osc}$ and normalized $R$ matrix}\label{sec:cat of aff q-osc}

\subsection{The category $\widehat{\mc O}_{\rm osc}$}\label{subsec:aff osc cat type 0}
Let us introduce our main object in this paper.
\begin{df}
{\rm
Let $\widehat{\mc O}_{\rm osc}$ be the category of $U_q(\agl_n)$-modules $V$ such that
\begin{itemize}
\item[(1)] $V=\bigoplus_{\mu\in P}V_\mu$ with $\dim V_\mu<\infty$, and the set of weights of $V$ is finitely dominated,

\item[(2)] $V$ belongs to ${\mc O}_{\rm osc}$ when restricted as a $U_q(\gl_n)$-module.
\end{itemize}}
\end{df}

The category $\widehat{\mc O}_{\rm osc}$ is a full subcategory of the category $\widehat{\mathcal{O}}$ introduced by Hernandez \cite{H05}, and it is closed under taking  submodules, quotients, and tensor product.
Note that $\W_l(x)\in \widehat{\mc O}_{\rm osc}$ for $l\in \Z$ and $x\in \Bbbk^\times$ since $\W_l(x)\cong \V^{(l)}$ as a $U_q(\gl_n)$-module.

Let us introduce the affinization of a module in $\widehat{\mc O}_{\rm osc}$, which is defined in the same way as in the case of finite-dimensional $U_q(\agl_n)$-module.
First, let
\begin{itemize}
\item[$\bullet$] $P_{\rm af} = \Z\Lambda_r \oplus \Z\de_1\oplus\cdots\oplus \Z\de_n\oplus\Z\delta \supset P$,

\item[$\bullet$]  $\mb{\alpha}_i=\de_i-\de_{i+1} +\de_{i0}\de\in P_{\rm af}$ ($i\in I$),

\item[$\bullet$] $P_{\rm af}^\vee={\rm Hom}_\mathbb{Z}(P_{\rm af},\mathbb{Z})=\Z c\oplus \Z\de^\vee_1\oplus\cdots\oplus\Z\de^\vee_n\oplus\Z d \supset P^\vee$, where the pairing on $P_{\rm af}\times P^\vee_{\rm af}$ is extended such that 
$$
\lang \de, \de^\vee_j \rang = 
\lang \de,c \rang=\lang \La_r,d \rang=\lang \de_i,d \rang=0,\quad
\lang \de,d\rang =1,
$$

\item[$\bullet$] ${\rm cl} : P_{\rm af} \longrightarrow P$ : the linear map given by ${\rm cl}\vert_{P}={\rm id}_P$ and ${\rm cl}(\de)=0$,

\item[$\bullet$] $\iota : P \longrightarrow P_{\rm af}$ : a section of ${\rm cl}$ given by $\iota(\delta_i)=\de_i$ for $1\le i\le n$.

\end{itemize}

Let $U_q(\widetilde{\mathfrak{gl}}_n)$ be the associative $\Bbbk$-algebra which is defined by the same relations as in $U_q(\agl_n)$, where $P$ and $P^\vee$ are replaced by $P_{\rm af}$ and $P_{\rm af}^\vee$.
Let $z$ be an indeterminate. For a $U_q(\agl_n)$-module $V$ in $\widehat{\mc O}_{\rm osc}$, we define a $U_q(\widetilde{\mathfrak{gl}}_n)$-module $V_{\rm aff}$ by 
\begin{equation*}\label{eq:affinization}
V_{\rm aff} = \Bbbk[z,z^{-1}]\ot V, 
\end{equation*}
where $q^h$, $e_i$, and $f_i$ acts as $1\ot q^h$, $z^{\delta_{i,0}}\ot e_i$, and $z^{-\delta_{i,0}}\ot f_i$, respectively ($h\in P^\vee$, $i\in I$), and $q^d(z^k\ot v)=q^k z^k\ot v$ ($k\in \Z$, $v\in V$).
For $\la\in P_{\rm af}$, we define the $\la$-weight space of $V_{\rm aff}$ by
\begin{equation*}
(V_{\rm aff})_{\la}= z^k\ot V_{{\rm cl}(\la)},
\end{equation*}
where $\la- \iota \circ {\rm cl}(\la)=k\de$. 
We have
$V_{\rm aff}=\bigoplus_{\la\in P_{\rm af}}(V_{\rm aff})_\la$, 
and $e_i(V_{\rm aff})_\la \subset (V_{\rm aff})_{\la+\mb{\alpha}_i}$, 
$f_i(V_{\rm aff})_\la \subset (V_{\rm aff})_{\la-\mb{\alpha}_i}$ for $i\in I$ and $\la\in P_{\rm af}$.
The map sending $g(z)\ot m$ to $zg(z)\ot m$ for $m\in V$ and $g(z)\in \Bbbk[z,z^{-1}]$ gives an isomorphism of $U_q(\agl_n)$-modules 
\begin{equation*}
z : V_{\rm aff} \longrightarrow V_{\rm aff},
\end{equation*}
such that $z(V_{\rm aff})_\la \subset (V_{\rm aff})_{\la+\de}$ for $\la\in P_{\rm af}$.  
For $x\in \Bbbk^\times$, we define a $U_q(\agl_n)$-module
\begin{equation*}
V_x = V_{\rm aff}/(z-x)V_{\rm aff}.
\end{equation*}
For example, if $V=\W_l(1)$ ($l\in \Z$), then $V_x=\W_l(x)$ for $x\in \Bbbk^\times$.

\subsection{Irreducibility of $\W_l (x)\ot \W_m (y)$}\label{subsec:Irreducibility}
Suppose that $l,m\in \Z$ and $x,y\in \Bbbk^\times$ are given.
%We regard $\W_l=\W_l(1)$ and $\W_m=\W_m(1)$ as $U_q(\agl_n)$-modules.
Let $v_l\in \W_l(x)$ and $v_m\in \W_m(y)$ be given as in \eqref{eq:v_l}, and write $v_{l}=\ket{\mathbf{m}}$ and $v_{m}=\ket{\mathbf{m}^{\prime}}$.

Let $l_1=\max(l,m)$ and $l_2=\min(l,m)$ so that $l_1\ge l_2$.
Let
\begin{equation}\label{eq:L}
 L=\max\{-l_1,l_2,0\}=
\begin{cases}
l_2 & \text{if }l_2\geq0,\\
-l_1 & \text{if } l_1 \leq 0,\\
0 & \text{if } l_1 \geq 0\geq l_2.
\end{cases}
\end{equation}
Equivalently, $L$ is the smallest non-negative integer such that $l_1+L\geq0\geq l_2-L$. 

By Propositions \ref{prop:hwsub-W^l} and \ref{prop:semisimplicity}, $\W_l(x)\ot \W_m(y)$ is a direct sum of $\mc{V}^{\nu}$ for $\nu\in \mc{P}(GL_2)_{(r,n-r)}$.
By considering its classical limit and \eqref{eq:LR for fundamentals}, we obtain the following decomposition of $\W_l(x)\ot \W_m(y)$ as a $U_q(\gl_n)$-module:
\begin{equation}\label{eq:classical decomp of fundamentals}
 \W_l(x)\ot\W_m(y) = \bigoplus_{i=-L}^\infty \V^{(l_1+L+i,l_2-L-i)}.
\end{equation}

Set 
\begin{equation*}
v_{a,b}^{+}=\ket{\mathbf{m}+a(\be_{r}+\be_{r+1})}\otimes\ket{\mathbf{m}^{\prime}+b(\be_{r}+\be_{r+1})},
\end{equation*}
for $a,b\ge 0$ and
\[
v_{a,b}^{-}=
\begin{cases}
\ket{\mathbf{m}+a(-\be_{r+1}+\be_{r+2})}\otimes\ket{\mathbf{m}^{\prime}+b(-\be_{r+1}+\be_{r+2})} & \text{if } l_2\geq0,\\
\ket{\mathbf{m}+a(\be_{r-1}-\be_{r})}\otimes\ket{\mathbf{m}^{\prime}+b(\be_{r-1}-\be_{r})} & \text{if } l_1 \leq0,
\end{cases}
\]
for $a,\,b\ge 0$ with $a+b\leq L$. Note that $v_{0,0}^\pm=v_l\ot v_m$.

\begin{lem}\label{lem:sing-vector}
The $U_{q}(\mathfrak{gl}_{n})$-highest weight vector of $\V^{(l_1+L+i,l_2-L-i)}$ in $\W_l(x) \ot \W_m(y)$ for $i\ge -L$ is equal (up to non-zero scalar multiplication) to the vector $u_i$ given as follows:
\begin{enumerate}
\item for $i\geq0$,
\[
u_{i}=\sum_{j=0}^{i}\left[(-1)^{j}\prod_{k=1}^{j}\left(q^{-(\abs m+2i-2k+1)}\frac{[\abs m+i+1-k][i+1-k]}{[\abs l+k][k]}\right)\right]v_{j,i-j}^{+},
\]
where the coefficient of $v_{0,i}^+$ is understood to be $1$,
\item for $-L\le i\le -1$,
\[
u_{i}=\sum_{j=0}^{-i}\left[(-1)^{j}\prod_{k=1}^{j}\left(q^{|m|+2i+2k}\frac{[-i+1-k]}{[k]}\right)\right]v_{j,-i-j}^{-}.
\]
\end{enumerate}
\end{lem}
\pf For $a,b\ge 0$, we have
\begin{equation}\label{eq:formula e v_{a,b}}
\begin{split}
 e_r v^+_{a,b} &= - q^{2b+|m|+1}[|l|+a][a]v^+_{a-1,b} - [|m|+b][b]v^+_{a,b-1},\\
 e_{r+1} v^-_{a,b} &= q^{-m+2b}[a]v^-_{a-1,b} + [b]v^-_{a,b-1}\quad (l_2\ge 0),\\
 e_{r-1} v^-_{a,b} &= q^{m+2b}[a]v^-_{a-1,b} + [b]v^-_{a,b-1}\quad (l_1\le 0),
\end{split}
\end{equation}
and $e_k v_{a,b}^{\pm}=0$ for all $k\neq 0,r,r\pm 1$.
Using \eqref{eq:formula e v_{a,b}}, it is straightforward to check that $e_ku_i=0$ for all $k\in I\setminus\{0\}$. Indeed, we only need to verify $e_ku_i=0$ for $k=r, r\pm 1$ by definition of $u_i$, which we leave it to the reader.
\qed\newline

Next, we describe how the $U_q(\gl_n)$-highest weight vectors are related under the action of $U_q(\agl_n)$.
Put
\begin{align*}
\boldsymbol{F}^{+} & =(e_{r+1}\cdots e_{n-2}e_{n-1})(e_{r-1}\cdots e_{2}e_{1})e_{0},\\
\boldsymbol{F}^{-} & =
\begin{cases}
e_{r}(e_{r+2}\cdots e_{n-2}e_{n-1})(e_{r-1}\cdots e_{2}e_{1})e_{0} & \text{if }l_2\geq0,\\
e_{r}(e_{r+1}\cdots e_{n-2}e_{n-1})(e_{r-2}\cdots e_{2}e_{1})e_{0} & \text{if }l_1\leq0,
\end{cases}\\
\boldsymbol{E}^{+} & =f_{0}(f_{1}\cdots f_{r-2}f_{r-1})(f_{n-1}\cdots f_{r+2}f_{r+1}),\\
\boldsymbol{E}^{-} & =
\begin{cases}
f_{0}(f_{n-1}\cdots f_{r+3}f_{r+2})(f_{1}\cdots f_{r-2}f_{r-1})f_{r} & \text{if }l_2\geq0,\\
f_{0}(f_{1}\cdots f_{r-3}f_{r-2})(f_{n-1}\cdots f_{r+2}f_{r+1})f_{r} & \text{if }l_1\leq0.
\end{cases}
\end{align*}

\begin{lem}\label{lem:EF action on v_{a,b}}
We have the following.
\begin{enumerate}
\item For $a,b\ge 0$, 
\begin{align*}
\boldsymbol{F}^+ v_{a,b}^{+} & =yv_{a,b+1}^{+}+xq^{-\abs m-2b-1}v_{a+1,b}^{+},\\
f_{r}v_{a,b}^{+} & =q^{-\abs l-2a-1}v_{a,b+1}^{+}+v_{a+1,b}^{+},\\
\boldsymbol{E}^{+} v_{a,b}^{+} & = - x^{-1}[a+\abs l][a]v_{a-1,b}^{+}-y^{-1}q^{2a+\abs l+1}[b+\abs m][b]v_{a,b-1}^{+}.
\end{align*}

\item For $a,b\ge 0$ with $a+b\leq L$ when $l_2\geq0$,
\begin{align*}
\boldsymbol{F}^{-} v_{a,b}^{-} & = - y[m-b]v_{a,b+1}^{-} - xq^{m-2b}[l-a]v_{a+1,b}^{-},\\
f_{r+1}v_{a,b}^{-} & =q^{l-2a}[m-b]v_{a,b+1}^{-}+[l-a]v_{a+1,b}^{-},\\
\boldsymbol{E}^{-} v_{a,b}^{-} & =-x^{-1}[a]v_{a-1,b}^{-}-y^{-1}q^{-l+2a}[b]v_{a,b-1}^{-}.
\end{align*}

\item For $a,b\ge 0$ with $a+b\leq L$ when $l_1\leq0$,
\begin{align*}
\boldsymbol{F}^{-} v_{a,b}^{-} & =-y[-m-b]v_{a,b+1}^{-}-xq^{-m-2b}[-l-a]v_{a+1,b}^{-},\\
f_{r-1}v_{a,b}^{-} & =q^{-l-2a}[-m-b]v_{a,b+1}^{-}+[-l-a]v_{a+1,b}^{-},\\
\boldsymbol{E}^{-} v_{a,b}^{-} & =-x^{-1}[a]v_{a-1,b}^{-}-y^{-1}q^{l+2a}[b]v_{a,b-1}^{-}.
\end{align*}
\end{enumerate}
\end{lem}

\begin{lem}\label{lem:F action on u_i}
For generic $x,y\in\Bbbk^\times$, we have
\begin{enumerate}
\item for $i\geq 0$, $\boldsymbol{F}^+ u_i \in \Bbbk^\times u_{i+1}+\Bbbk f_r u_i + \Bbbk (1-\delta_{i0})f^{(2)}_r u_{i-1}$,
\item for $-L \leq i \leq 0$ when $l_2 \geq 0$, $\boldsymbol{F}^- u_i \in \Bbbk^\times u_{i-1} + \Bbbk f_{r+1}u_{i}+\Bbbk (1-\delta_{i0})f^{(2)}_{r+1}u_{i+1}$,
\item for $-L \leq i \leq 0$ when $l_1 \leq 0$, $\boldsymbol{F}^- u_i \in \Bbbk^\times u_{i-1} + \Bbbk f_{r-1}u_{i}+\Bbbk (1-\delta_{i0})f^{(2)}_{r-1}u_{i+1}$.
\end{enumerate} 
\end{lem}

\pf
Let us prove only (1), leaving the other two to the reader. Considering the classical decomposition \eqref{eq:classical decomp of fundamentals} and weights, we can write
\begin{equation*}
  \boldsymbol{F}^+ u_i =\sum_{k\geq 0} C^{k}_i f^{(k)}_{r}u_{i+1-k}
\end{equation*}
for some $C^{k}_{i}\in\Bbbk$. First, a direct computation gives us the following identity
\begin{equation*}
e^{2}_r \boldsymbol{F}^+ u_i = [2][\left| m\right|+i][i](yq^{\left|m\right|-2i+1}-xq^{\left|l\right|+1})u_{i-1},
\end{equation*}
which implies $C^{k}_i =0$ whenever $k>2$. Moreover, we obtain
\begin{equation*}
  C^{2}_{i}=\frac{[2][\left| m\right|+i][i]}{[\left|l\right|+\left| m\right|+2i+1][\left| l\right|+\left| m\right|+2i]}(yq^{\left|m\right|-2i+1}-xq^{\left|l\right|+1}).
\end{equation*}
Next, by comparing the coefficients of $v^{+}_{i,0}$ of both sides of 
\begin{equation*}
  e_r \boldsymbol{F}^+ u_i = -C^{1}_{i} [\left| l\right|+\left| m\right|+2i] u_i - C^{2}_{i} [\left| l\right|+\left| m\right|+2i+1]f_r u_{i-1},
\end{equation*}
we can compute $C^{1}_{i}$. Now substituting $C^{1}_{i},C^{2}_{i}$ in 
\begin{equation*}
  \boldsymbol{F}^+ u_i = C^{0}_{i}u_{i+1} +C^{1}_{i}f_r u_i + C^{2}_{i} f^{(2)}_{r} u_{i-1},
\end{equation*}
and comparing the coefficients of $v^{+}_{i+1,0}$, it can be easily seen that $C^{0}_{i}$ is a polynomial in $x,y$, if we regard them as indeterminates. The assertion follows.
\qed

\begin{thm}\label{thm:irr of tensor product of fund 0}
For generic $x,y\in\Bbbk^\times$, the tensor product $\mathcal{W}_{l}(x)\otimes\mathcal{W}_{m}(y)$ is an irreducible $U_{q}(\widehat{\mathfrak{gl}}_{n})$-module in $\widehat{\mc O}_{\rm osc}$.
\end{thm}
\pf
Let $K$ be a non-zero $U_{q}(\widehat{\mathfrak{gl}}_{n})$-submodule
of $\mathcal{W}_{l}(x)\otimes\mathcal{W}_{m}(y)$. 
Since $\mathcal{W}_{l}(x)\otimes\mathcal{W}_{m}(y)$ is semisimple over $U_{q}(\mathfrak{gl}_{n})$, $K$ contains at least
one $u_{i}$. 
We claim that
\begin{enumerate}
\item $u_{0}\in K$ if $u_{i}\in K$ for some $i\geq-L$,
\item $u_{0}$ generates $\mathcal{W}_{l}(x)\otimes\mathcal{W}_{m}(y)$.
\end{enumerate}
Let us first prove (1). Since $\mathcal{W}_{l}(x)\otimes\mathcal{W}_{m}(y)$ is semisimple over $U_q (\mathfrak{gl}_n)$, $K$ contains at least one highest weight vector, say $u_i$.

Suppose that $-L \leq i\leq -1$. We shall prove that $\boldsymbol{E}^- u_i$ is a non-zero multiple of $u_{i+1}$ so that $u_{i+1}\in K$, and then recursively we obtain $u_0 \in K$.
First, observe that for generic $x,y$, $\boldsymbol{E}^- u_i$ is not zero, which can be seen from a direct computation using Lemma~\ref{lem:EF action on v_{a,b}}(1).
Then it remains to check that $\boldsymbol{E}^- u_i$ is a $U_q(\mathfrak{gl}_n)$-highest weight vector.
Indeed, the only nontrivial thing to check is $e_j \boldsymbol{E}^- u_i =0$ for $j=r,r+1$ if $l_2 \geq 0$ or $j=r,r-1$ if $l_1 \leq 0$.
If $l_1 \leq 0$, then
\begin{align*}
  e_r \boldsymbol{E}^- u_i & =f_0 (f_{1}\cdots f_{r-2})(f_{n-1}\cdots f_{r+1})\frac{k_r -{k_r}^{-1}}{q-q^{-1}}u_i\\
  & =[l+m-i-2]f_0 (f_{1}\cdots f_{r-2})(f_{n-1}\cdots f_{r+1})u_i =0,\\
  e_{r-1}\boldsymbol{E}^- u_i &= e_{r+1} (f_{1}\cdots f_{r-2})(f_{n-1}\cdots f_{r+1})f_r u_i \boldsymbol{E}^- e_{r-1}u_i=0.
\end{align*}
The case if $l_2 \geq 0$ is similar.

Thanks to Lemma~\ref{lem:F action on u_i}(1) (resp. (2)), it is an easy induction to see $u_i \in K$ for $i>0$ (resp. $-L\leq i\leq -1$), provided $u_0 \in K$. This completes the proof of the theorem.

\qed\newline

For indeterminates $z_1$ and $z_2$, let 
\begin{equation*}
\begin{split}
(\mathcal{W}_{l})_{\rm aff}=\Bbbk[z_1^{\pm 1}]\ot\mathcal{W}_{l},\quad
(\mathcal{W}_{m})_{\rm aff}=\Bbbk[z_2^{\pm 1}]\ot\mathcal{W}_{m},\\
(\mathcal{W}_{l})^{{}^{\wedge}}_{\rm aff}=\Bbbk(z_1)\ot \mathcal{W}_{l},\quad
(\mathcal{W}_{m})^{{}^{\wedge}}_{\rm aff}=\Bbbk(z_2)\ot \mathcal{W}_{m}.
\end{split}
\end{equation*}

\begin{cor}\label{cor:irreducibility}
The tensor product $(\mathcal{W}_{l})^{{}^{\wedge}}_{\rm aff} \ot (\mathcal{W}_{m})^{{}^{\wedge}}_{\rm aff}$ is an irreducible $\Bbbk(z_1,z_2)\ot U_q(\agl_n)$-module.
\end{cor}
\pf 
If we replace the spectral parameters $x,y$ by $z_1 , z_2$ in the proof of Lemma~\ref{lem:F action on u_i}, 
then it can be easily seen that the coefficient of $u_{i\pm1}$ in $\boldsymbol{F}^{\pm}u_i$ is a polynomial in $z_1$ and $z_2$.
Hence, over the field of rational functions $\Bbbk(z_1 ,z_2)$, the proof of Theorem~\ref{thm:irr of tensor product of fund 0} applies equally well.
\qed 

\subsection{Normalized $R$ matrix}\label{subsec:normalized R matrix}

We keep the notations in Section \ref{subsec:Irreducibility}.
%From now on, let us assume that $\Bbbk$ is the algebraic closure of $\Q(q)$ in $\bigcup_n \C(\!(q^{-1/n})\!)$.
Recall that as a $U_q(\gl_n)$-module, we have the following decomposition:
\begin{equation}\label{eq:LR for q fundamentals}
\W_l\ot \W_m \cong \W_m\ot \W_l = 
\bigoplus_{t\in\Z_+} \mc{V}^{(l_1+t,l_2-t)},
\end{equation}
where $l_1=\max\{l,m\}$ and $l_2=\min\{l,m\}$. (see \eqref{eq:LR for fundamentals} and Remark \ref{rem:LR for q osc}).
For $t\in\Z_+$, choose a highest weight vector  $v(l,m,t)$  of $\mc{V}^{(l_1+t,l_2-t)}$, which is unique up to scalar multiplication by Lemma \ref{lem:sing-vector}, and let $\mc P^{l,m}_{t} : \W_{l} \otimes \W_{m} \longrightarrow \W_{m} \otimes \W_{l}$ be
the $U_q(\gl_n)$-linear map given by $\mc P^{l,m}_{t}(v(l,m,t'))=\delta_{tt'}v(m,l,t')$ for $t'\in\Z_+$.

Let 
\begin{equation*}
(\mathcal{W}_{l})_{\rm aff}\, \widehat{\ot}\,(\mathcal{W}_{m})_{\rm aff},\quad 
(\mathcal{W}_{m})_{\rm aff}\, \widetilde{\ot}\, (\mathcal{W}_{l})_{\rm aff}
\end{equation*} 
be the completions of $(\mathcal{W}_{l})_{\rm aff}\,{\ot}\,(\mathcal{W}_{m})_{\rm aff}$ and 
$(\mathcal{W}_{m})_{\rm aff}\,{\ot}\,(\mathcal{W}_{l})_{\rm aff}$ following \cite[Section 7]{Kas02}.
The universal $R$ matrix for $U_q(\agl_n)$ induces a map
\begin{equation*}
\xymatrixcolsep{4pc}\xymatrixrowsep{3pc}\xymatrix{
(\mathcal{W}_{l})_{\rm aff}\, \widehat{\ot}\,(\mathcal{W}_{m})_{\rm aff} \ \ar@{->}^{\mc{R}^{\rm univ}_{(l,m)}}[r] &\ (\mathcal{W}_{m})_{\rm aff}\, \widetilde{\ot}\, (\mathcal{W}_{l})_{\rm aff} },
\end{equation*}
(see \cite[(7.6)]{Kas02}).
We can check that
\begin{equation*}
\mc{R}^{\rm univ}_{(l,m)}(v_l\ot v_m)= a\left(\tfrac{z_1}{z_2}\right)(v_m\ot v_l),
\end{equation*}
for some non-zero $a(\frac{z_1}{z_2})\in \Bbbk\llbracket \frac{z_1}{z_2}\rrbracket$, which is  invertible. 

Put $z=z_1 /z_2$ and
\begin{equation}\label{eq:c(z)}
 c(z)= \begin{cases}
\prod_{i=1}^{\min\{\left|l\right|,\left|m\right|\}}\frac{1-q^{\left|l-m\right|+2i}z}{z-q^{\left|l-m\right|+2i}} & \text{if } lm>0, \\
1 & \text{if } lm\leq 0.
\end{cases}
\end{equation}
We define the normalized $R$ matrix by
\begin{equation}\label{eq:norm R}
\mc{R}^{\rm norm}_{(l,m)}= c(z)a(z)^{-1}\mc{R}^{\rm univ}_{(l,m)}.
\end{equation}
By \eqref{eq:LR for q fundamentals} and Lemma \ref{lem:sing-vector}, we have
\begin{equation*}\label{eq:normalization}
 \mc{R}^{\rm norm}_{(l,m)}\big\vert_{\V^{(l_1+L,l_2-L)}}=c(z){\rm id}_{\V^{(l_1+L,l_2-L)}}.
\end{equation*}
Since 
%$\mc{R}^{\rm norm}_{l,m}\big\vert_{\mc{V}^{(l,m)}}= c(z){\rm id}_{\mc{V}^{(l,m)}}$, 
$\mc{R}^{\rm norm}_{(l,m)}(v_l \otimes v_m) = c(z) v_m \otimes v_l$,
there exists a $\Bbbk[z_1^{\pm 1},z_2^{\pm 1}]\ot U_q(\agl_n)$-linear map
\begin{equation}\label{eq:R matrix on W_l tensor W_m}
\xymatrixcolsep{3pc}\xymatrixrowsep{3pc}\xymatrix{
(\mathcal{W}_{l})_{\rm aff} \ot (\mathcal{W}_{m})_{\rm aff}\ \ar@{->}^{\hskip -2cm \mc{R}^{\rm norm}_{(l,m)}}[r] &  \
\Bbbk(z_1,z_2)\ot_{\Bbbk[z_1^{\pm 1},z_2^{\pm 1}]}\left((\mathcal{W}_{m})_{\rm aff}\,{\ot}\, (\mathcal{W}_{l})_{\rm aff}\right)
}
\end{equation}
by Corollary \ref{cor:irreducibility}, and $\mc{R}^{\rm norm}_{(l,m)}$ is a unique $\Bbbk(z_1,z_2)\ot U_q(\agl_n)$-linear map \eqref{eq:R matrix on W_l tensor W_m}
satisfying $\mc{R}^{\rm norm}_{(l,m)}(v_l \otimes v_m) = c(z) v_m \otimes v_l$.
%$\mc{R}^{\rm norm}_{l,m}\big\vert_{\mc{V}^{(l,m)}}= {\rm id}_{\mc{V}^{(l,m)}}$.

The following spectral decomposition of $\mc{R}^{\rm norm}_{(l,m)}$ is the main result in this section. The proof will be given in Section \ref{sec:spectral decomposition} (Corollary \ref{cor:spectral decomp}).

\begin{thm}\label{thm:spectral decomposition}
For $l, m\in\mathbb{Z}$, we have
\begin{equation*}\label{eqn:spectral-decomp}
\begin{split}
\mc{R}^{\rm norm}_{(l,m)} 
&= \mc P^{l,m}_0 + \sum_{t=1}^\infty\prod_{i=1}^{t}\dfrac{1-q^{|l-m|+2i}z}{z-q^{|l-m|+2i}} \mc P^{l,m}_t\\
&= \mc P^{l,m}_0 
+ \dfrac{1-q^{|l-m|+2}z}{z-q^{|l-m|+2}} \mc P^{l,m}_1 
+ \dfrac{1-q^{|l-m|+2}z}{z-q^{|l-m|+2}}\dfrac{1-q^{|l-m|+4}z}{z-q^{|l-m|+4}} \mc P^{l,m}_2 + \dots ,
\end{split}
\end{equation*}
after a suitable choice of $v(l,m,t)$ and $v(m,l,t)$ for $t\in\Z_+$.
\end{thm}

\subsection{The $\ell$-highest weight of $\W_l$}
The subalgebra $U_q(\widehat{\mathfrak{sl}}_{n})$ of $U_{q}(\widehat{\mathfrak{gl}}_{n})$ generated by $k_i^{\pm 1}, e_i, f_i$ ($i\in I$) has another set of generators $k_i^{\pm 1}, h_{i,s}^{\pm}, x_{i,t}^{\pm}$ and $C^{\pm1/2}$ ($i\in I\setminus\{0\}, s\in\mathbb{Z}\setminus\{0\}, t\in \Z$) (see \cite{BCP} for precise presentation).

Let $\psi_{i,k}$ ($i\in I\setminus\{0\}, k\ge 0$) be the element determined by the following identity of formal power series in $z$:
\begin{equation*}\label{eq:psi generators}
\sum_{k=0}^\infty \psi_{i,k} z^{k} = k_i\exp\left((q-q^{-1}) \sum_{s=1}^\infty h_{i,s} z^{s} \right).
\end{equation*}
A $U_{q}(\widehat{\mathfrak{sl}}_{n})$-module $V$ is called an $\ell$-highest weight module if it is generated by an $\ell$-highest weight vector $v$, that is,
$x_{i,k}^{+}v=0$ and $\psi_{i,k}v=\Psi_{i, k}v$ for all $i\in I\setminus\{0\},k\ge 0$ and some scalars $\Psi_{i, k}$. 
The series $\Psi=(\Psi_{i,k})_{i\in I\setminus\{0\},k\ge 0}$ is called the $\ell$-highest weight of $V$.

Recall that every (type 1) finite-dimensional irreducible module is an $\ell$-highest weight module, and an $\ell$-highest weight module is finite-dimensional if $\Psi=(\Psi_{i,k})_{i\in I\setminus\{0\},k\ge 0}$ satisfies
\begin{equation*}\label{eq:Drinfeld polynomial}
\Psi_{i}(z):=\sum_{k\geq 0} \Psi_{i,k}z^k = q^{\mathrm{deg}P_i}\frac{P_i (q^{-2}z)}{P_i(z)},
\end{equation*}
for uniquely determined polynomials $P_i(z)\in\Bbbk[z]$ with constant term $1$(see \cite[Theorem~3.3]{CP}, and also \cite{MY} for the characterization of $\ell$-highest weights of irreducible modules in $\widehat{\mc{O}}$).

Let us verify that $\mathcal{W}_{l}$ is an $\ell$-highest weight $U_{q}(\widehat{\mathfrak{sl}}_{n})$-module for all $l\in \mathbb{Z}$, and compute its $\ell$-highest weight. We refer the reader to \cite{BCP} for notations and construction of root vectors $E_{k\delta -\alpha_i}$ for $k\ge 1$ and $i\in I\setminus\{0\}$.
Fix a map $o:I\setminus\left\{0\right\}\rightarrow \left\{\pm 1\right\}$ such that $o(i+1)=-o(i)$ for all $i\in I\setminus \left\{0\right\}$. We recall the following lemmas.
\begin{lem}[{\cite[Lemma~1.5]{BCP}}]\label{lem:aux-1}
 For $i\in I\setminus\{0\}$ and $k>0$, we have
\[
\psi_{i,k}=o(i)^{k}(q-q^{-1})k_{i}\left(E_{k\delta-\alpha_{i}}e_{i}-q^{-2}e_{i}E_{k\delta-\alpha_{i}}\right).
\]
\end{lem}

\begin{lem}[{\cite[Lemma~4.3]{JKP}}]\label{lem:aux-2}
 For $i\in I\setminus\{0\}$ and $k>0$, we have
\begin{align*}
E_{(k+1)\delta-\alpha_{i}}=-\frac{1}{q+q^{-1}} & \left(E_{\delta-\alpha_{i}}e_{i}E_{k\delta-\alpha_{i}}-q^{-2}e_{i}E_{\delta-\alpha_{i}}E_{k\delta-\alpha_{i}}\right.\\
 & \left.-E_{k\delta-\alpha_{i}}E_{\delta-\alpha_{i}}e_{i}+q^{-2}E_{k\delta-\alpha_{i}}e_{i}E_{\delta-\alpha_{i}}\right).
\end{align*}
\end{lem}

\begin{lem}[{\cite[Lemma~4.7]{JKP}}]\label{lem:aux-3}
 For $i\in I\setminus\{0\}$, we have
\[
E_{\delta-\alpha_{i}}=(-q^{-1})^{n-2}(e_{i+1}\dots e_{n-1})(e_{i-1}\dots e_{2}e_{1})e_{0}+\sum_{j_{1},\dots,j_{n-1}}C_{j_{1},\dots,j_{n-1}}(q)e_{j_{1}}\dots e_{j_{n-1}},
\]
where the sum is over the sequences $(j_{1},\dots,j_{n-1})\in I^{n-1}$ such that
$\sum_{k=1}^{n-1}\alpha_{j_{k}}=\sum_{j\in I\setminus\left\{i\right\}}\alpha_j$ with $j_{n-1}\neq0$
and $C_{j_{1},\dots,j_{n-1}}(q)\in\pm q^{-\mathbb{Z}_{+}}$. 
\end{lem}

\begin{thm}
The fundamental $q$-oscillator representation $\mathcal{W}_{l}$
is an $\ell$-highest weight module with $\ell$-highest weight $\Psi=(\Psi_{i,k})_{i\in I\setminus\{0\},k\ge 0}$
given by
\[
\Psi_{i}(z)=\sum_{k\ge 0}\Psi_{i,k}z^k=\begin{cases}
(q^{-l}+u)(1+q^{-l}u)^{-1} & \text{if $l<0$ and $i=r-1$},\\
(u+q^{-\left|l\right|-1})(1+q^{-\left|l\right|-1}u)^{-1} & \text{if $i=r$},\\
(q^{l}+u)(1+q^{l}u)^{-1} & \text{if $l\geq 0$ and $i=r+1$},\\
1 & \text{otherwise},
\end{cases}
\]
where $u=o(r)(-q^{-1})^{n}z$.
\end{thm}

\pf
We claim that $v_{l}$ is an $\ell$-highest weight vector with the
given $\ell$-highest weight. By weight consideration, $x_{i,t}^{+}v_{l}=0$
for all $i\in I\setminus\left\{0\right\},\,t\in\mathbb{Z}$. 
Let us show that
that $v_{l}$ is a simultaneous eigenvector of $\psi_{i,k}$ with
the eigenvalues $\Psi_{i,k}$ above.

First assume $l\geq0$. 
Since $e_j v_l=0$ for all $j\in I\setminus\left\{0\right\}$, we have by Lemma \ref{lem:aux-3},
\begin{align*}
E_{\delta-\alpha_{r}}v_{l} & =(-q^{-1})^{n-2}(e_{r+1}\cdots e_{n-1})(e_{r-1}\cdots e_{1})e_{0}v_{l}\\
 & =(-q^{-1})^{n-2}\ket{\be_{r}+(l+1)\be_{r+1}},
\end{align*}
and then
\begin{align*}
E_{2\delta-\alpha_{r}}v_{l} 
& =-\frac{1}{q+q^{-1}}\left(E_{\delta-\alpha_{r}}e_{r}E_{\delta-\alpha_{r}}-q^{-2}e_{r}E_{\delta-\alpha_{r}}^{2}+q^{-2}E_{\delta-\alpha_{r}}e_{r}E_{\delta-\alpha_{r}}\right)v_{l}\\
& =-\frac{1}{q+q^{-1}}\\
&\quad \left(-(-q^{-1})^{2n-4}[l+1](1+q^{-2})+q^{-2}(-q^{-1})^{2n-4}[2][l+2]\right)\ket{\be_{r}+(l+1)\be_{r+1}}\\
& =(-q^{-1})^{2n-4}(q^{-1}[l+1]-q^{-2}[l+2])\ket{\be_{r}+(l+1)\be_{r+1}} \\
&=-(-q^{-1})^{2n-4}q^{-l-3}\ket{\be_{r}+(l+1)\be_{r+1}}.
\end{align*}
Repeating similar computation and using Lemma \ref{lem:aux-2}, we have
\[
E_{k\delta-\alpha_{r}}v_{l}=(-1)^{k-1}(-q^{-1})^{(n-2)k}(q^{-l-3})^{k-1}\ket{\mathbf{e}_{r}+(l+1)\mathbf{e}_{r+1}},
\]
and by Lemma \ref{lem:aux-1}
\begin{align*}
\psi_{r,k}v_{l} & =o(r)^{k}(q-q^{-1})k_{r}\left(E_{k\delta-\alpha_{r}}e_{r}-q^{-2}e_{r}E_{k\delta-\alpha_{r}}\right)v_{l}\\
 & =o(r)^{k}(q-q^{-1})q^{-l-1}q^{-2}[l+1](-1)^{k-1}(-q^{-1})^{(n-2)k}(q^{-l-3})^{k-1}v_{l}\\
 & =o(r)^{k}(q^{-l-1}-q^{l+1})(-q^{-1})^{nk}(-q^{-l-1})^{k}v_{l}.
\end{align*}
Thus we obtain
\begin{align*}
\Psi_{r}(z)v_{l} & =\sum_{k\geq0}\psi_{r,k}z^kv_{l} 
 =\left(q^{-l-1}-(q^{l+1}-q^{-l-1})\sum_{k\geq1}\left\{ o(r)(-q^{-1})^{n}(-q^{-l-1})z\right\} ^{k}\right)v_{l}\\
 & =\left(q^{-l-1}-(q^{l+1}-q^{-l-1})\frac{-q^{-l-1}u}{1+q^{-l-1}u}\right)v_{l} =\frac{u+q^{-l-1}}{1+q^{-l-1}u}v_{l},
\end{align*}
where $u=o(r)(-q^{-1})^{n}z$.
The computation of $\Psi_{r+1}(z)$ is similar, where we begin from
\begin{align*}
E_{\delta-\alpha_{r+1}}v_{l} & =(-q^{-1})^{n-2}(e_{r+2}\dots e_{n-1})(e_{r}\dots e_{2}e_{1})e_{0}v_{l}\\
 & =-(-q^{-1})^{n-2}[l]\ket{(l-1)\be_{r+1}+\be_{r+2}}.
\end{align*}
For $i\ne r,\,r+1$, it is obvious since $E_{\delta-\alpha_{i}}v_{l}=0$.
The case when $l<0$ can be dealt with by the same computation.
\qed

\section{Irreducible representations in $\widehat{\mc{O}}_{\rm osc}$}\label{sec:irreducible of affine q-osc}

\subsection{Fusion construction of irreducible representations}
For $i\in I\setminus\{0\}$ and $a\in\Bbbk^\times$, let $V(\varpi_i)_a$ be the $i$-th \emph{fundamental representation} of $U_q(\widehat{\mf{sl}}_n)$ with spectral parameter $a$ (see \cite{Kas02} for more details).
Recall that any (type 1) finite-dimensional irreducible representation $V$ of $U_q(\widehat{\mf{sl}}_n)$ can be realized as an irreducible quotient of the submodule generated by the tensor product of $\ell$-highest weight vectors in $V(\varpi_{i_1})_{a_1}\ot\cdots\ot V(\varpi_{i_\ell})_{a_\ell}$ for some $i_s$ and $a_s$ ($1\le s\le\ell$) \cite[Corollary~3.6]{CP} (over the algebraic closure of $\Bbbk$).
Moreover, if we choose a suitable order of the tensor product so that every normalized $R$ matrices on $V(\varpi_{i_s})_{a_s}\ot V(\varpi_{i_t})_{a_t}$ ($s<t$) are defined, then the tensor product itself is generated by the tensor product of $\ell$-highest weight vectors \cite[Proposition~9.4]{Kas02}, and the  irreducible representation $V$ is obtained as the image of the composition of the normalized $R$ matrices
\begin{equation*}
 \xymatrixcolsep{2.5pc}\xymatrixrowsep{3pc}\xymatrix{
V(\varpi_{i_1})_{a_1}\ot\cdots\ot V(\varpi_{i_\ell})_{a_\ell} \ \ar@{->}^{}[r] & \ V(\varpi_{i_\ell})_{a_\ell}\ot\cdots\ot V(\varpi_{i_1})_{a_1}}.
\end{equation*}
Such a procedure to construct irreducible representations is called \emph{fusion construction}. 

In this section, we construct a family of irreducible representations in $\widehat{\mathcal{O}}_\mathrm{osc}$ by using fusion construction applied to fundamental $q$-oscillator representations.

For $\ell\geq 2$, let $z_1 ,\,\dots,\,z_\ell$ be indeterminates and $l_1 ,\,\dots,\,l_\ell \in \mathbb{Z}_+$ given. 
Let $\mc{W}_{l_i}(z_i)$ denote $(\mathcal{W}_{l_i})_{\rm aff}=\Bbbk[z_i^{\pm 1}]\ot\mathcal{W}_{l_i}$ for $1\le i\le \ell$.

Put ${\mb l}=(l_1,\dots,l_\ell)$.
Since the normalized $R$ matrices \eqref{eq:norm R} satisfy the Yang-Baxter equation, we may define a $\Bbbk[z^{\pm 1}_1,\dots,z^{\pm 1}_\ell]\otimes U_q(\widehat{\mathfrak{gl}}_n)$-linear map 
\begin{equation*}\label{eq:norm R_l}
\xymatrixcolsep{2.5pc}\xymatrixrowsep{3pc}\xymatrix{
\mathcal{W}_{l_1}(z_1)\otimes\cdots\otimes\mathcal{W}_{l_\ell}(z_\ell) \ 
\ar@{->}^{\hskip-2cm\mathcal{R}^{\mathrm{norm}}_{{\mb l}}}[r] & \
\Bbbk(z_1,\dots,z_\ell)\otimes_{\Bbbk[z^{\pm 1}_1,\dots,z^{\pm 1}_\ell]} \mathcal{W}_{l_\ell}(z_\ell)\otimes\cdots\otimes \mathcal{W}_{l_1}(z_1)},
\end{equation*}
by choosing any reduced expression of the longest element of the symmetric group $\mathfrak{S}_\ell$ and taking the composition of the normalized $R$ matrices \eqref{eq:norm R} in that order.

\begin{thm}\label{thm:fusion}
Let $(\boldsymbol{l},\boldsymbol{c})$ be a pair of sequences $\boldsymbol{l}=(l_{1},\dots,l_{\ell})\in\mathbb{Z}^{\ell}$ and 
$\boldsymbol{c}=(c_{1},\dots,c_{\ell})\in(\Bbbk^{\times})^{\ell}$ such that $c_{i}/c_{j}$ is not of the form $q^{\left|l_{i}-l_{j}\right|+2a}$ $(a\in \mathbb{N})$ for all $1\leq i<j\leq \ell$. Then
\begin{itemize}
 \item[(1)] the specialization of $\mathcal{R}_{\boldsymbol{l}}^{\mathrm{norm}}$ at $(z_{1},\dots,z_{\ell})=(c_{1},\dots,c_{\ell})$ gives a $U_q(\widehat{\mathfrak{gl}}_n)$-linear map 
\[
\xymatrixcolsep{2.5pc}\xymatrixrowsep{3pc}\xymatrix{ 
R_{\mb{l}}(\mb{c}):\mathcal{W}_{l_{1}}(c_{1})\otimes\cdots\otimes\mathcal{W}_{l_{\ell}}(c_{\ell})\ \ar@{->}[r] 
& \ \mathcal{W}_{l_{\ell}}(c_{\ell})\otimes\cdots\otimes\mathcal{W}_{l_{1}}(c_{1})},
\] 
 
 \item[(2)] $\mathrm{Im}R_{\mb{l}}(\mb{c})$ is an irreducible representation in $\widehat{\mc O}_{\rm osc}$ if it is not zero.
\end{itemize}
\end{thm}
\pf (1) The well-definedness of $R_{\mb{l}}(\mb{c})$ follows from  Theorem \ref{thm:spectral decomposition}.

(2) The spectral decomposition in Theorem \ref{thm:spectral decomposition} implies that the specialization of $\mathcal{R}_{l,l}^{\mathrm{norm}}$
at $z_{1}=z_{2}=c$ for any $c\in\Bbbk^\times$ is the identity map $\mathrm{id}_{\mathcal{W}_{l}(c)\otimes\mathcal{W}_{l}(c)}$. Suppose that $\mathrm{Im}R_{\mb{l}}(\mb{c})$ is non-zero. Then we can prove that $\mathrm{Im}R_{\mb{l}}(\mb{c})$ is irreducible in the same manner as in \cite[Theorems~4.3, 4.4]{KL} and \cite[Lemma~2.6]{KK} replacing the renormalized $R$ matrices there with the normalized $R$ matrices. 
Note that $\mc{R}^{\rm norm}_{(l,m)}(z)$ is an infinite sum of $\mc P^{l,m}_t$, and hence it is not rationally renormalizable in the sense of \cite{KKKO15}. But we may still apply the arguments there using the normalized $R$ matrix (cf.~\cite[Appendix D]{KO21}).
\qed\newline

Let $\mc{P}^+_\ell$ be the set of $({\mb l},{\mb c})$ such that
\begin{itemize}
\item[(1)] ${\mb l}=(l_1,\dots,l_\ell)\in \Z^\ell$, and $\mb{c}=(c_1,\dots,c_\ell)\in (\Bbbk^\times)^\ell$,

\item[(2)] $c_{i}/c_{j}$ does not belong to $\{\,q^{\left|l_{i}-l_{j}\right|+2a}\,|\, a\in \mathbb{N}\,\}$ for all $1\leq i<j\leq \ell$ when $\ell\geq 2$,
\end{itemize}
and let $\mc{P}^+=\bigsqcup_{\ell\ge 1}\mc{P}^+_\ell$.
For $({\mb l},{\mb c})\in \mc{P}^+$, we put 
\begin{equation*}
 \W({\mb l},{\mb c})=\mathrm{Im}R_{\boldsymbol{l}}(\boldsymbol{c}).
\end{equation*} 

\begin{thm}\label{thm:fusion construction}
For $({\mb l},{\mb c})\in \mc{P}^+_\ell$, $\W({\mb l},{\mb c})$ is non-zero  if ${\mb l}^+\in \mc{P}(GL_\ell)_{(r,n-r)}$, 
where ${\mb l}^+$ denotes the rearrangement of ${\mb l}$ such that ${\mb l}^+\in \mc{P}(GL_\ell)$.
\end{thm}
\pf 
We may assume that $\ell\ge 2$ since it is clear for $\ell=1$. We use the results in Section \ref{sec:spectral decomposition}.
Keeping the notations in Section \ref{sec:spectral decomposition}, we consider the following sequences:
\begin{itemize}
 \item[$\bullet$] $\e=\e^{(a,b)}$ as in \eqref{eq:e(a,b)} for sufficiently large $a,b$,
 
 \item[$\bullet$] $\underline{\e}=(0^{a+b+1})$ with $\underline{\e}_-=(0^a)$ and $\underline{\e}_+=(0^{b+1})$,
 
 \item[$\bullet$] $\ov{\e}=(1^{a+b})$ with $\ov{\e}_-=(1^a)$ and $\ov{\e}_+=(1^{b})$,
 
 \item[$\bullet$] $\underline{\e}'=(0^n)$ with $\underline{\e}'_-=(0^r)$ and $\underline{\e}'_+=(0^{n-r})$,
 
\end{itemize}
where we assume that $a\ge r$ and $b\ge n-r$ so that $\underline{\e}'$ is a subsequence of $\underline{\e}$. 

For each sequence ${\mb \e}$ above, let $\U({\mb\e})$ be the generalized quantum group in Section \ref{subsec:GQG}.
Let $\W_{l,{\mb \e}}$ ($l\in \Z$) denote the $q$-oscillator representation of $\U({\mb\e})$ in Proposition \ref{prop:osc module-3}, $\mc{R}^{\rm norm}_{(l,m),{\mb \e}}$ the normalized $R$ matrix on $(\mathcal{W}_{l,{\mb \e}})_{\rm aff} \ot (\mathcal{W}_{m,{\mb \e}})_{\rm aff}$ in \eqref{eq:R matrix on W_l tensor W_m for e}, and $\V^{\la}_{\mb\e}$ the irreducible representation of $\mathring{\U}({\mb\e})$ appearing in $\W_{l,{\mb \e}}^{\ot \ell}$ ($\ell\ge 1$) (see Section \ref{subsec:q-osc for e} and Corollary \ref{cor:decomp of W for e}). We also consider the category $\mc{C}({\mb \e})$, and then the truncation functors between them (see Section \ref{subsec:truncation functor}).

Let $({\mb l},{\mb c})\in \mc{P}^+_\ell$ be given such that ${\mb l}^+\in \mc{P}(GL_\ell)_{(r,n-r)}$. 
For each sequence ${\mb \e}$ above, we also have 
a $\Bbbk[z^{\pm}_1,\dots,z^{\pm}_\ell]\otimes \mc{U}({\mb \e})$-linear map 
\begin{equation*}\label{eq:norm R_l e}
\xymatrixcolsep{2.5pc}\xymatrixrowsep{3pc}\xymatrix{
\mathcal{W}_{l_1,{\mb \e}}(z_1)\otimes\cdots\otimes\mathcal{W}_{l_\ell,{\mb \e}}(z_\ell) \ 
\ar@{->}^{\hskip-2cm\mathcal{R}^{\mathrm{norm}}_{{\mb l},{\mb \e}}}[r] & \
\Bbbk(z_1,\dots,z_\ell)\otimes_{\Bbbk[z^{\pm}_1,\dots,z^{\pm}_\ell]} \mathcal{W}_{l_\ell,{\mb \e}}(z_\ell)\otimes\cdots\otimes \mathcal{W}_{l_1,{\mb \e}}(z_1)}.
\end{equation*}
Thanks to Theorem \ref{thm:spectral decomposition for e}, we have a $\mc{U}({\mb \e})$-linear map
\begin{equation*}
\xymatrixcolsep{2.5pc}\xymatrixrowsep{3pc}\xymatrix{ 
R_{\mb{l},{\mb \e}}(\mb{c}):\mathcal{W}_{l_{1},{\mb \e}}(c_{1})\otimes\cdots\otimes\mathcal{W}_{l_{\ell},{\mb \e}}(c_{\ell})\ 
\ar@{->}[r] & \
 \mathcal{W}_{l_{\ell},{\mb \e}}(c_{\ell})\otimes\cdots\otimes\mathcal{W}_{l_{1},{\mb \e}}(c_{1})},
\end{equation*}
by specializing at $(z_{1},\dots,z_{\ell})=(c_{1},\dots,c_{\ell})$.
We have
\begin{equation}\label{eq:truncation of fusion}
 \mf{tr}^{\e}_{\underline{\e}}
 \left( R_{\mb{l},{\e}}(\mb{c}) \right)=
 R_{\mb{l},{\underline{\e}}}(\mb{c}),\quad
  \mf{tr}^{\e}_{\ov{\e}}
 \left(\kappa R_{\mb{l},{\e}}(\mb{c}) \right)=
  R_{\mb{l},{\ov{\e}}}(\mb{c}),
\end{equation}
by Lemma \ref{lem:truncation from e to (0)} and Corollary \ref{cor:truncation from e to (1)}, respectively, where $\kappa$ is the constant given in \eqref{eq:kappa}.

First, consider the case when ${\mb \e}=\ov{\e}=(1^{a+b})$. 
Recall that $\U(\ov{\e})$ is isomorphic to $U_{\td{q}}(\agl_{a+b})$, where $\td{q}=-q^{-1}$, and $\W_{l_i,\ov{\e}}(c_i)$ is isomorphic to the $(a+l_i)$-th fundamental representation with spectral parameter $c_i$ (see the proof of Theorem \ref{thm:spectral decomposition for e}).

For $1\le i\le \ell$, let $w_i$ be a $\mathring{\U}(\ov{\e})$-highest weight vector of $\W_{l_i,\ov{\e}}(c_i)$, which is unique up to scalar multiplication.
Then the $\mathring{\mc{U}}(\ov{\e})$-submodule of $\mathcal{W}_{l_{i},{\mb \e}}(c_{i})\otimes \mathcal{W}_{l_{j},{\mb \e}}(c_{j})$ generated by $w_i\ot w_j$ is isomorphic to $\mc{V}_{\ov{\e}}^{{(l_i,l_j)}^+}$ for all $1\le i,j\le \ell$ ($i\neq j$).
Since ${\mb l}^+\in \mc{P}(GL_\ell)_{(r,n-r)}$, we have ${\mb l}^+\in \mc{P}(GL_\ell)_{(\e_-,\e_+)}$.
By Theorem \ref{thm:spectral decomposition for e}, we have
\begin{equation}\label{eq:R sends extremal vector to itself}
R_{\mb{l},\ov{\e}}(\mb{c})(w_1\ot \dots\ot w_\ell)=c (w_\ell\ot\dots\ot w_1), 
\end{equation}
for some $c\in \Bbbk^\times$. In particular, ${\rm Im}R_{\mb{l},\ov{\e}}(\mb{c})$ is non-zero.

Since the $\mathring{\mc{U}}(\ov{\e})$-submodule of $\mathcal{W}_{l_{1},\ov{\e}}(c_{1})\otimes\cdots\otimes\mathcal{W}_{l_{\ell},\ov{\e}}(c_{\ell})$ (resp. $\mathcal{W}_{l_{\ell},\ov{\e}}(c_{\ell})\otimes\cdots\otimes\mathcal{W}_{l_{1},\ov{\e}}(c_{1})$) generated by $w_1\ot \dots\ot w_\ell$ (resp. $w_\ell\ot\dots\ot w_1$) is isomorphic to $\mc{V}_{\ov{\e}}^{{\mb l}^+}$, we have 
$R_{\mb{l},\ov{\e}}(\mb{c})(\mc{V}_{\ov{\e}}^{{\mb l}^+})=\mc{V}_{\ov{\e}}^{{\mb l}^+}$ by \eqref{eq:R sends extremal vector to itself}.

The multiplicity of $\mc{V}_{\ov{\e}}^{{\mb l}^+}$ in $\mathcal{W}_{l_{1},\ov{\e}}(c_{1})\otimes\cdots\otimes\mathcal{W}_{l_{\ell},\ov{\e}}(c_{\ell})$ or $\mathcal{W}_{l_{\ell},\ov{\e}}(c_{\ell})\otimes\cdots\otimes\mathcal{W}_{l_{1},\ov{\e}}(c_{1})$ is equal to one by Pieri rule for the usual Schur polynomials.
By Theorem \ref{thm:trunc sends simple to simple}, it is also true for ${\mb \e}=\e$. Hence we have $R_{\mb{l},{\e}}(\mb{c})(\mc{V}_{{\e}}^{{\mb l}^+})=\mc{V}_{{\e}}^{{\mb l}^+}$ by the second identity in \eqref{eq:truncation of fusion}.

Since ${\mb l}^+\in \mc{P}(GL_\ell)_{(r,n-r)}=\mc{P}(GL_\ell)_{(\underline{\e}'_-,\underline{\e}'_+)}\subset \mc{P}(GL_\ell)_{(\underline{\e}_-,\underline{\e}_+)}$, we have 
$\mf{tr}^\e_{\underline{\e}}\left(\mc{V}^{\la}_\e\right)=\mc{V}^{\la}_{\underline{\e}}$ by Theorem \ref{thm:trunc sends simple to simple}.  
Hence by the first identity in \eqref{eq:truncation of fusion} we conclude that
$R_{\mb{l},\underline{\e}}(\mb{c})(\mc{V}_{\underline{\e}}^{{\mb l}^+})=\mc{V}_{\underline{\e}}^{{\mb l}^+}$, which implies that $R_{\mb{l},\underline{\e}}(\mb{c})$ is non-zero and so is $R_{\mb{l},\underline{\e}'}(\mb{c})=R_{\mb{l}}(\mb{c})$. 
\qed 

\begin{rem}\label{rem:fusion for e}
{\rm %The proof of Theorem \ref{thm:fusion construction} indeed gives a correspondence between irreducible objects in $\widehat{\mc{O}}_{\rm osc, {\mb \e}}$.  
Let us keep the notations in Section \ref{sec:spectral decomposition}.
It is not difficult to see that Theorem \ref{thm:fusion} also holds for ${\mb \e}$ by the same arguments and Theorem \ref{thm:spectral decomposition for e}. 
For $({\mb l},{\mb c})\in \mc{P}^+_\ell$, let $\W_{{\mb \e}}({\mb l},{\mb c})={\rm Im}R_{\mb{l},\mb{\e}}(\mb{c})$.
By \eqref{eq:truncation of fusion}, we have
\begin{equation}\label{eq:tr of simples}
\mf{tr}^{\e}_{\underline{\e}}
 \left(\W_{{\e}}({\mb l},{\mb c}) \right)=
 \W_{{\underline \e}}({\mb l},{\mb c}),\quad
  \mf{tr}^{\e}_{\ov{\e}}
 \left(\W_{\e}({\mb l},{\mb c}) \right)=
   \W_{{\ov \e}}({\mb l},{\mb c}).
\end{equation}
Hence $\mf{tr}^\e_{{\mb\e}}$ yields a correspondence between irreducible representations when $\W_{{\mb \e}}({\mb l},{\mb c})$ are non-zero.

Consider the functors $\mf{tr}^{\e}_{\underline{\e}}$, $\mf{tr}^{\e}_{\underline{\e}'}$, and $\mf{tr}^{\e}_{\ov{\e}}$:
\begin{equation*}\label{eq:triple of truncation}
\xymatrixcolsep{3pc}\xymatrixrowsep{0.3pc}\xymatrix{
& & \mc{C}(\e)  \ar@{->}_{\mf{tr}^{\e}_{\underline{\e}}}[dl]\ar@{->}^{\mf{tr}^{\e}_{\ov{\e}}}[dr] &  \\
& \mc{C}(\underline{\e}) \ar@{->}_{\mf{tr}^{\underline{\e}}_{\underline{\e}'}}[dl] & &  \mc{C}(\ov{\e}) \\
\mc{C}(\underline{\e}') & & &
}
\end{equation*}
By Theorem \ref{thm:trunc sends simple to simple}, this induces the following:
\begin{equation}\label{eq:triple of truncation osc}
\xymatrixcolsep{3pc}\xymatrixrowsep{0.3pc}\xymatrix{
& & \widehat{\mc{O}}_{\rm osc,\e}  \ar@{->}_{\mf{tr}^{\e}_{\underline{\e}}}[dl]\ar@{->}^{\mf{tr}^{\e}_{\ov{\e}}}[dr] &  \\
& \widehat{\mc{O}}_{\rm osc,\underline{\e}}\ar@{->}_{\mf{tr}^{\underline{\e}}_{\underline{\e}'}}[dl] & &  \widehat{\mc{O}}_{\rm osc,\ov{\e}}\\
\widehat{\mc{O}}_{\rm osc,\underline{\e}'} & & &
}
\end{equation}
where $\widehat{\mc{O}}_{\rm osc}=\widehat{\mc{O}}_{\rm osc,\underline{\e}'}$ and $\widehat{\mc{O}}_{\rm osc,\ov{\e}}$ is the category of finite-dimensional representations of $U_{\td{q}}(\widehat{\gl}_{a+b})$ (see Remark \ref{rem:odd homogeneous}). 
Let $({\mb l},{\mb c})\in \mc{P}^+_\ell$ be given such that $\W({\mb l},{\mb c})=\W_{\underline{\e}'}({\mb l},{\mb c})\neq 0$. 
We have $\W_{\e}({\mb l},{\mb c})\neq 0$ and $\mf{tr}^\e_{\ov{\e}}(\W_{\e}({\mb l},{\mb c}))=\W_{\ov{\e}}({\mb l},{\mb c})$ for all sufficiently large $a,b\ge 1$.
Furthermore, by \eqref{eq:tr of simples} and Theorem \ref{thm:trunc sends simple to simple}, we also have  
\begin{equation}\label{eq:multiplicity of classical comp}
{\rm Hom}_{U_q(\gl_n)}(\V^{\la},\W({\mb l},{\mb c}))= 
{\rm Hom}_{\mathring{\mc U}(\ov{\e})}(\V^{\la}_{\ov{\e}},\W_{\ov{\e}}({\mb l},{\mb c})),
\end{equation}
when both of them are non-zero.

We expect that there is an equivalence between $\widehat{\mc{O}}_{\rm osc,\e}$ and $\widehat{\mc{O}}_{\rm osc,{\mb \e}}$ for ${\mb \e}\in \{\,\underline{\e}, \ov{\e}\,\}$, which yields the correspondence \eqref{eq:tr of simples} (after taking a suitable limit of $\e$ and ${\mb \e}$).
This would be an affine and quantum analogue of the super duality between the integrable highest weight representations and oscillator representations of $U(\gl_{\infty})$ (see \cite{CW}).
} 
\end{rem}

\subsection{Kirillov-Reshetikhin type modules}
We define an analogue of Kirillov-Reshetikhin type modules in $\widehat{\mc{O}}_{\rm osc}$ as follows. 
For $l\in\Z$, $s\in \mathbb{N}$ and $c\in \Bbbk^\times$, let 
\begin{equation*}
 \mathcal{W}^{l,s}(c)=\W({\mb l},{\mb c}),
\end{equation*}
where $\boldsymbol{l}=(\underbrace{l,\dots,l}_{s\ \text{times}})$ and $\boldsymbol{c}=c(q^{2-2s},\dots,q^{-2},1)$.

\begin{thm}
We have $ \mathcal{W}^{l,s}(c)\neq 0$ if and only if $(l^s)\in \mc{P}(GL_s)_{(r,n-r)}$. In this case, $\mathcal{W}^{l,s}(c)$ is classically irreducible, that is,  we have as a $U_q(\gl_n)$-module
\[
\mathcal{W}^{l,s}(c)\cong \mathcal{V}^{(l^{s})}.
\]
\end{thm}
\pf The first assertion follows from Theorem \ref{thm:fusion construction}.
So it is enough to show the second one.
Let ${\e}$ and $\ov{\e}$ denote the sequences given in the proof of Theorem \ref{thm:fusion construction}. Suppose that $a,b$ are sufficiently large.

Note that ${\rm Im}R_{{\mb l},\ov{\e}}({\mb c})$ is the Kirillov-Reshetikhin module of type $A_{a+b-1}^{(1)}$, and ${\rm Im}R_{{\mb l},\ov{\e}}({\mb c})\cong \mc{V}_{\ov{\e}}^{(l^s)}$ as a $\mathring{\mc{U}}(\ov{\e})$-module.
Consider ${\rm Im}R_{{\mb l},{\e}}({\mb c})$. If there is a submodule isomorphic to $\mc{V}_{\e}^\la$ for $\la\neq (l^s)$, then we have $\mf{tr}^\e_{\ov{\e}}(\mc{V}_{\e}^\la)=\mc{V}_{\ov{\e}}^\la \neq 0$ for all sufficiently large $a,b$ by Theorem \ref{thm:trunc sends simple to simple} (cf.~\eqref{eq:GL e}).
This is a contradiction.
Hence ${\rm Im}R_{{\mb l},{\e}}({\mb c})\cong \mc{V}_{\e}^{(l^s)}$ as a $\mathring{\mc{U}}({\e})$-module. This implies our assertion since 
$\mathcal{W}^{l,s}(c)=\mf{tr}^\e_{\underline{\e}'}\left({\rm Im}R_{{\mb l},{\e}}({\mb c})\right)$.
\qed

\begin{rem}\label{rem:non-zero KR}
{\rm
Note that $\mathcal{W}^{0,s}(c)\neq 0$ for all $s\ge 1$.
On the other hand, when $l\neq 0$, 
we have $\mathcal{W}^{l,s}(c)$ is not zero if and only if 
\begin{equation*}
\begin{cases}
 s\le n-r & \text{if $l> 0$},\\
 s\le r & \text{if $l< 0$}.
\end{cases}
\end{equation*}
} 
\end{rem}

\begin{cor}
If $(l^s)\in \mc{P}(GL_s)_{(r,n-r)}$, then the character of $\mathcal{W}^{l,s}(c)$ as a $U_q(\gl_n)$-module is given by
\begin{equation*}
 {\rm ch}\mathcal{W}^{l,s}(c)=t^s\sum_{\mu,\nu}c^{(l^s)}_{\mu\nu^*}s_\mu(x)s_\nu(y).
\end{equation*}
In particular, we have for all $s\ge 1$
\begin{equation*}
\begin{split}
 {\rm ch}\mathcal{W}^{0,s}(c)
 &=t^s\sum_{\ell(\mu)\le \min\{r,n-r,s\}}s_\mu(x)s_\mu(y)  \\
 &=\dfrac{t^s}{\prod_{i,j}(1-x_iy_j)} \quad\quad \text{if $s\ge \max\{r,n-r\}$}. 
\end{split}
\end{equation*}
\end{cor}
\pf It follows from Proposition \ref{prop:ch of classical osc} and Remark \ref{rem:non-zero KR}.
\qed

\begin{rem}
{\rm 
Let $({\mb l},{\mb c})\in \mc{P}^+_\ell$ be such that $\W({\mb l},{\mb c})$ is non-zero and classically irreducible. Assume further that ${\mb l}=[\mu,\nu]_\ell$ as in \eqref{eq:generalized partition}. 
Since $\V^{\mb l}\subset \W({\mb l},{\mb c})$ by Remark \ref{rem:LR for q osc}, we have $\V^{\mb l}= \W({\mb l},{\mb c})$.
If $\ell\ge n$, then we have by Corollary \ref{cor:ch of classical osc} 
\begin{equation*}
{\rm ch}\W({\mb l},{\mb c})={\rm ch}\V^{\mb l}=t^\ell s_\mu(x)s_{\nu}(y)\dfrac{1}{\prod_{i,j}(1-x_iy_j)}.
\end{equation*}
Indeed, a classically irreducible $\W({\mb l},{\mb c})$ naturally corresponds to the one in the category of finite-dimensional representations in the sense of \eqref{eq:triple of truncation osc}.
We have shown in the proof of Theorem \ref{thm:fusion} that $\W_{{\ov \e}}({\mb l},{\mb c})$ is not zero if $a,b$ are sufficiently large. Then it is not difficult to see from \eqref{eq:multiplicity of classical comp} that if $\W_{{\ov \e}}({\mb l},{\mb c})$ is classically irreducible, then so is $\W({\mb l},{\mb c})$.
}
\end{rem}

\section{$q$-oscillator representations of quantum affine superalgebras}\label{sec:spectral decomposition}

\subsection{Generalized quantum group $\mc{U}(\e)$}\label{subsec:GQG}
Let us first recall the notion of the generalized quantum groups (of affine type $A$) \cite{KOS}, which is also closely related to the quantum affine superalgebras  (cf.~\cite{KL,Ya99}).

We assume the following notations in this section:
\begin{itemize}

\item[$\bullet$] $\e=(\e_1,\dots,\e_n)$ : a sequence with $\e_i\in \{0,1\}$ ($i=1,\dots, n$),

\item[$\bullet$] $I=\{\,0,1,\ldots,n-1\,\}$,

\item[$\bullet$] $r\in I\setminus\{0,1,n-1\}$,
 
\item[$\bullet$] $\mathbb{I}= \{1<2<\cdots <n\}$ : a linearly ordered set with $\mathbb{Z}_2$-grading $\mathbb{I}=\mathbb{I}_0\sqcup\mathbb{I}_1$ such that $$\mathbb{I}_0=\{\,i\,|\,\e_i=0\,\},\quad \mathbb{I}_1=\{\,i\,|\,\e_i=1\,\},$$

\item[$\bullet$] $\I^-=\{\,1,\dots,r\,\}$ and $\I^+=\{\,r+1,\dots,n\,\}$

\item[$\bullet$] $P = \Z\La_{r,\e} \oplus \Z\de_1\oplus\cdots\oplus \Z\de_n$ with a symmetric bilinear form $(\,\cdot\,|\,\cdot\,)$ such that 
\begin{equation*}\label{eq:extended bilinear form}
(\de_i|\de_j)=(-1)^{\e_i}\delta_{ij}, \quad 
(\La_{r,\e}|\La_{r,\e})=0, \quad 
(\de_i|\La_{r,\e})=
\begin{cases}
0 & (i\in \I^-),\\
1 & (i\in \I^+),
\end{cases}
\end{equation*}

\item[$\bullet$] $\alpha_i=\de_i-\de_{i+1} \in P$ $(i\in I)$,

\item[$\bullet$] $I_{\rm even}=\{\,i\in I\,|\,(\alpha_i|\alpha_i)=\pm 2\,\}$, 
$I_{\rm odd}=\{\,i\in I\,|\,(\alpha_i|\alpha_i)=0\,\}$,

\item[$\bullet$] $q_i=\si q^{\si}$ $(i\in \mathbb{I})$, that is,
\begin{equation*}
q_i=
\begin{cases}
q & \text{if $\e_i=0$},\\
-q^{-1} & \text{if $\e_i=1$},\\
\end{cases} \quad (i\in \mathbb{I}),
\end{equation*}

\item[$\bullet$] ${\bq}(\,\cdot\,,\,\cdot\,)$,  $\hat{\bq}(\,\cdot\,,\,\cdot\,)$ : symmetric biadditive functions from $P\times P$ to $\Bbbk^{\times}$ given by
\begin{equation*}
\begin{split}
&\bq(\mu,\nu) = \prod_{i\in \mathbb{I}}q_i^{\mu_i \nu_i},\quad \hat{\bq}(\mu,\nu) = q^{\sum_{j\in \I^+}(\ell'\mu_j+\ell\nu_j)}\bq(\mu,\nu) \quad (\mu,\nu\in P),
\end{split}
\end{equation*}
where $\mu=\ell\La_{r,\e}+\sum_{i\in \I}\mu_i\de_i$ and $\nu=\ell'\La_{r,\e}+\sum_{i\in \I}\nu_i\de_i$.

\end{itemize}
Then we define ${\mathcal{U}}(\e)$ to be the associative $\Bbbk$-algebra with $1$ 
generated by $k_\mu, e_i, f_i$ for $\mu\in P$ and $i\in I$ 
satisfying
{\allowdisplaybreaks
\begin{gather*}
k_0=1, \quad k_{\mu +\mu'}=k_{\mu}k_{\mu'} \quad (\mu, \mu' \in P),\label{eq:Weyl-rel-1-e} \\ 
k_\mu e_i k_{-\mu}=\bq(\mu,\alpha_i)e_i,\quad 
k_\mu f_i k_{-\mu}=\bq(\mu,\alpha_i)^{-1}f_i, \\ 
e_if_j - f_je_i =\delta_{ij}\frac{k_{\alpha_i} - k_{-\alpha_i}}{q-q^{-1}}, \\
e_i^2= f_i^2 =0 \quad (i\in I_{\rm odd}),
\end{gather*}
\begin{gather*}\label{eq:Serre-rel-1}
e_i e_j -  e_j e_i = f_i f_j -  f_j f_i =0
 \quad \text{($i-j\not\equiv \pm 1\!\!\pmod n$)},\\ 
\begin{array}{ll}
e_i^2 e_j- (-1)^{\e_i}[2] e_i e_j e_i + e_j e_i^2= 0\\ 
f_i^2 f_j- (-1)^{\e_i}[2] f_i f_j f_i+f_j f_i^2= 0
\end{array}
\quad \text{($i\in I_{\rm even}$, $i-j\equiv \pm 1\!\!\pmod n$)}, 
\end{gather*}
\begin{gather*}\label{eq:Serre-rel-2}
\begin{array}{ll}
  e_{i}e_{i-1}e_{i}e_{i+1}  
- e_{i}e_{i+1}e_{i}e_{i-1} 
+ e_{i+1}e_{i}e_{i-1}e_{i} \\  
\hskip 2cm - e_{i-1}e_{i}e_{i+1}e_{i} 
+ (-1)^{\e_i}[2]e_{i}e_{i-1}e_{i+1}e_{i} =0, \\ 
  f_{i}f_{i-1}f_{i}f_{i+1}  
- f_{i}f_{i+1}f_{i}f_{i-1} 
+ f_{i+1}f_{i}f_{i-1}f_{i}  \\  
\hskip 2cm - f_{i-1}f_{i}f_{i+1}f_{i} 
+ (-1)^{\e_i}[2]f_{i}f_{i-1}f_{i+1}f_{i} =0,
\end{array}\quad \text{($i\in I_{\rm odd}$)}.
\end{gather*}}
We call $\mathcal{U}(\e)$ the {\em generalized quantum group of affine type $A$ associated to $\e$}. It has a Hopf algebra structure with the comultiplication $\Delta$ and the antipode $S$
\begin{equation}\label{eq:comult-e}
\begin{split}
\Delta(k_\mu)&=k_\mu\otimes k_\mu, \\ 
\Delta(e_i)&= 1\ot e_i + e_i\ot k_i^{-1}, \\
\Delta(f_i)&= f_i\ot 1 + k_i\ot f_i , \\  
S(q^h)=q^{-h},& \ \ S(e_i)=-e_i k_{\alpha_i}, \ \  S(f_i)=-k_{\alpha_i}^{-1} f_i,\\
\end{split}
\end{equation}
for $\mu\in P$ and  $i\in I$.
Let $\ring{\mathcal{U}}(\e)$ be the $\Bbbk$-subalgebra of $\mathcal{U}(\e)$ generated by $k_\mu$, $e_i$ and $f_i$ for $\mu\in P$ and $i\in I\setminus \{0\}$.

\begin{rem}\label{rem:gqg-epsilon-0}{\rm
Suppose that $\e=(0,0,\dots,0)$. We can identify the weight lattice $P$ of $\mathcal{U}(\epsilon)$ with the dual weight lattice $P^\vee$ of $U_q(\widehat{\mathfrak{gl}}_n)$ by an isomorphism $\phi$ sending $\Lambda_{r,\epsilon}$ to $\de^\vee_{r+1}+\dots+\de^\vee_{n}$ and $\de_i$ to $\de_i^\vee-\sigma(i)c$, where $\sigma(i)=0$ (resp. $\sigma(i)=1$) for $1\le i\le r$ (resp. $r<i\le n$). Under this identification, $\mathcal{U}(\epsilon)$ is isomorphic to $U_q(\agl_n)$.
Moreover, $P$ can also be identified with the weight lattice of $U_q(\widehat{\mathfrak{gl}}_n)$ by $\psi$ sending $\Lambda_{r,\epsilon}$ to $-\Lambda_r$ and $\delta_i$ to $\delta_i$, which naturally relates the bilinear form on $P$ to the canonical pairing by $(\lambda |\mu)=\left<\psi(\lambda),\phi(\mu) \right>$.
Similarly, when $\e=(1,1,\dots,1)$, $\mc{U}(\e)$ is isomorphic to $U_{\td{q}}(\agl_n)$, where $\td{q}=-q^{-1}$.
}
\end{rem}

For a $\mathcal{U}(\e)$-module $V$ and $\la\in P$, we define the $\la$-weight space to be
\begin{equation*}
V_\la 
= \{\,u\in V\,|\,k_{\mu} u= \hat{\bq}(\la,\mu) u \ \ (\mu\in P) \,\},
\end{equation*}
and let
${\rm wt}(V)=\{\,\mu\in P\,|\,V_\mu\neq 0\,\}$.

\subsection{Polynomial representations}\label{subsec:poly repn}\mbox{}
Suppose that $M$ and $N$ denote the numbers of $0$'s and $1$'s in $\e$, respectively.
A partition $\la=(\la_i)_{i\ge 1}$ is called an $(M|N)$-hook partition if $\la_{M+1}\leq N$. We denote the set of all $(M|N)$-hook partitions by $\cP_{M|N}$.
For $\la\in\cP_{M|N}$, let $H_{\la}$ denote the tableau of shape $\la$ defined inductively as follows (see also \cite[(2.37)]{CW}):
\begin{enumerate}
    \item[(i)] If $\epsilon_1 =0$ (resp. $\e_1=1$), then fill the first row (resp. column) of $\lambda$ with $1$.
    \item[(ii)] After filling a subdiagram $\mu$ of $\lambda$ with $1,\dots,k$, fill the first row (resp. column) of the skew Young diagram $\lambda/\mu$ with $k+1$ if $\epsilon_{k+1} =0$ (resp. $\e_{k+1}=1$).
\end{enumerate}
Let $V_\e(\la)$ be the irreducible highest weight $\ring{\U}(\e)$-module with the highest weight $\sum_{i\in \mathbb{I}} m_i\delta_i$, where $m_i$ is the number of $i$'s in the tableau $H_\lambda$.
If we put $V=V_\e((1))$, then $V^{\otimes \ell}$ is semisimple for $\ell\ge 1$ and its irreducible components are $V_\e(\la)$'s for  $(M|N)$-hook partitions $\la$ of $\ell$ (cf.~\cite{BKK,KY}). We denote by $V_\e(-\la)$ the dual representation of $V_\e(\la)$ with respect to \eqref{eq:comult-e}.

%Let $\Z[P]$ be the group ring of the weight lattice with basis $\{\,e^\mu\,|\,\mu\in P\,\}$. Put $x_i=e^\de_i$ ($i\in \I$). Let $x=\{\,x_i\,|\,i\in \I\,\}$ be the set of commuting variables. For $\la\in \cP_{M|N}$, let $hs_\la(x)$ be the hook Schur polynomial, which is the character of $V_\e(\la)$ \cite{BR}.

%Let $r\in I\setminus\{0,1,n-1\}$ be as in Section \ref{sec:notations}. 
We write $\e=(\e_-,\e_+)$, where
\begin{equation}\label{eq:e for parabolic}
\e_-=(\e_1,\dots,\e_{r}),\quad \e_+=(\e_{r+1},\dots,\e_{n}).
\end{equation}
Let $\mathring{\U}(\e_-,\e_+)$ be the subalgebra of $\mathring{\U}(\e)$ generated by $k_\mu$, $e_i$ and $f_i$ for $\mu\in P$ and $i\in I\setminus\{0,r\}$.
Let $M_\pm$ and $N_\pm$ be the numbers of $0$'s and $1$'s in $\e_\pm$, respectively.

For $\mu\in \cP_{M_-|N_-}$ and $\nu\in \cP_{M_+|N_+}$, let 
\begin{equation*}\label{eq:mixed poly}
V_{(\e_-,\e_+)}(\mu,\nu)=V_{\e_-}(-\mu)\ot V_{\e_+}(\nu),
\end{equation*}
which is a representation of $\mathring{\U}(\e_-,\e_+)$.
It is clear that a tensor product $V_{(\e_-,\e_+)}(\alpha,\beta)\ot V_{(\e_-,\e_+)}(\gamma,\delta)$ is semisimple and decomposes into a sum of $V_{(\e_-,\e_+)}(\mu,\nu)$'s.
%Note that the character of $V_{(\e_-,\e_+)}(\mu,\nu)$ is $hs_\mu(x^{-1}_-)hs_\nu(x_+)$, where $x_-^{-1}=\{\,x_i^{-1}\,|\,i\in \I^-\,\}$. 

\subsection{$q$-oscillator representations of $\U(\e)$}\label{subsec:q-osc for e}
Let
\begin{equation*}
\Z^n_+(\e)=\{\,{\bf m}=(m_1,\ldots,m_n)\,|\,\e_i=0 \Rightarrow m_i\in\Z_+,\  \e_i=1 \Rightarrow m_i\in \{0,1\}\,\}.
\end{equation*}
We keep the notations \eqref{eq:notation for state vector}.
Let
\begin{equation*}
\begin{split}
\W_{\e} &= \bigoplus_{{\bf m}\in \Z^n_+(\e)}\Bbbk|{\bf m}\rangle,\\
\W_{l,\e} &= \bigoplus_{l({\bf m})=l}\Bbbk|{\bf m}\rangle \quad (l\in \Z).
\end{split}
\end{equation*}
Then we have the following, which is an analogue of Proposition \ref{prop:osc module-2}.

\begin{prop}\label{prop:osc module-3}
Let $x\in\Bbbk^\times$ be given.
\begin{itemize}
\item[(1)]
The following formula defines a $\U(\e)$-module structure on $\W_\e$, which we denote by $\W_\e(x)$:
{\allowdisplaybreaks
\begin{align*}
&
\begin{array}{l}
k_{\La_{r,\e}} \ket{\bf m}= q^{\sum_{j\in \I^+}m_j}\ket{\bf m}, \\
k_{\de_i}\ket{\bf m}=q_i^{-m_i}\ket{\bf m} \quad (1\le i\le r),\\
k_{\de_j}\ket{\bf m}=q_j^{m_j}q\ket{\bf m} \quad (r+1\le j\le n),\\
\end{array}\\
&
\begin{array}{l}
e_{0}\ket{\bf m} = x\ket{{\bf m}+\be_{1}+\be_{n}},\\
f_{0}\ket{{\bf m}} = -x^{-1}[m_{1}][m_{n}]\ket{{\bf m}-\be_{1}-\be_{n}},\\
\end{array}\\
&
\begin{array}{l}
e_{i}\ket{\bf m} =[m_{i}]\ket{{\bf m}-\be_{i}+\be_{i+1}}\\
f_{i}\ket{\bf m}  =[m_{i+1}]\ket{{\bf m}+\be_{i}-\be_{i+1}}\\
\end{array}\quad (1\leq i < r),\\
&
\begin{array}{l}
e_{r}\ket{\bf m} =-[m_{r}][m_{r+1}]\ket{{\bf m}-\be_{r}-\be_{r+1}},\\
f_{r}\ket{\bf m} =\ket{{\bf m}+\be_{r}+\be_{r+1}},\\
\end{array}\\
&
\begin{array}{l}
e_{j}\ket{\bf m} =[m_{j+1}]\ket{{\bf m}+\be_{j}-\be_{j+1}}\\
f_{j}\ket{\bf m}  =[m_{j}]\ket{{\bf m}-\be_{j}+\be_{j+1}}\\
\end{array}\quad (r < j\leq n-1).
\end{align*}}

\item[(2)] For $l\in \Z$, $\W_{l,\e}$ is an irreducible $\U(\e)$-module, which we denote by $\W_{l,\e}(x)$. Moreover, $\W_{l,\e}(x)$ is an irreducible $\mathring{\U}(\e)$-module.  
\end{itemize}
\end{prop}
\pf It is straightforward to check that the operators above satisfy the defining relations of $\U(\e)$. We leave the details to the reader.
\qed
\vskip 2mm

The weight of $\ket{\bf m}\in \W_\e(x)$ is given by
\begin{equation}\label{eq:weight of m e}
{\rm wt}(\ket{\bf m})= 
\La_{r,\e} - \sum_{i\in \I^-}m_i\delta_i  + \sum_{j\in \I^+}m_j\delta_j,
\end{equation} 

\begin{rem}\label{rem:odd homogeneous}
{\rm
Note that $\W_{l,\e}(x)$ is infinite-dimensional unless $\e_i=1$ for all $i$, in which case, it is isomorphic to the usual fundamental representation of type $A_{n-1}^{(1)}$ with spectral parameter $x$, up to a twist by an automorphism of $U_{\widetilde{q}}(\widehat{\mathfrak{gl}}_n)$:
$k_r \mapsto -k_r$, $e_r \mapsto -e_r$, $f_r \mapsto f_r$,
and fixing other generators. Here $\La_{r,\e}$ is the $r$-th fundamental weight. 

On the other hand, when $\epsilon_i =0$ for all $i$, $\mathcal{W}_{l,\epsilon}(x)$ is the fundamental $q$-oscillator representation of $U_q(\mathfrak{gl}_n)$ under the isomorphism explained in  Remark~\ref{rem:gqg-epsilon-0}. Recall that $\Lambda_{r,\epsilon}$ in \eqref{eq:weight of m e} is compatible with $-\La_r$ in  \eqref{eq:weight of m} by the identification of weight lattices (see Remark~\ref{rem:gqg-epsilon-0}).
} 
\end{rem}

\begin{lem}\label{lem:semisimple over levi}
Let $\e=(\e_-,\e_+)$ be as in \eqref{eq:e for parabolic}.
For $\ell\ge 1$ and $x\in \Bbbk^\times$, $\W_\e(x)^{\ot \ell}$ is semisimple as a $\mathring{\U}(\e_-,\e_+)$-module. An irreducible $\mathring{\U}(\e_-,\e_+)$-submodule of $\W_\e(x)^{\ot \ell}$ is isomorphic to $V_{(\e_-,\e_+)}(\mu,\nu)$ for some $\mu\in \cP_{M_-|N_-}$ and $\nu\in \cP_{M_+|N_+}$. 
\end{lem}
\pf Given $a_\pm\ge 0$, let
\begin{equation*}
\W_{(a_-,a_+),\e}(x) = \bigoplus_{| \mathbf{m}|_\pm=a_\pm}\Bbbk \ket{\bf m}.
\end{equation*}
Then we have
\begin{equation}\label{eq:decomp into mixed poly}
\begin{split}
&\W_{l,\e}(x)=\bigoplus_{a_+-a_-=l}\W_{(a_-,a_+),\e}(x),\\
&\W_{(a_-,a_+),\e}(x)\cong V_{(\e_-,\e_+)}((a_-),(a_+)),
\end{split}
\end{equation}
as a $\mathring{\U}(\e_-,\e_+)$-module.
More precisely, the righthand side should be understood as tensored by the trivial $\mathring{\U}(\e_-,\e_+)$-module with weight $\La_{r,\e}$, say $\Bbbk_{\La_{r,\e}}$.

Hence it follows from \eqref{eq:decomp into mixed poly} that $\W_\e(x)^{\ot \ell}$ is semisimple as a $\mathring{\U}(\e_-,\e_+)$-module and each of its irreducible components is isomorphic to $V_{(\e_-,\e_+)}(\mu,\nu)$, more precisely $V_{(\e_-,\e_+)}(\mu,\nu)\ot \Bbbk_{\ell\La_{r,\e}}$, for some $\mu\in \cP_{M_-|N_-}$ and $\nu\in \cP_{M_+|N_+}$.
\qed\newline

Next, let us consider the semisimplicity of $\W_{\e}(x)^{\ot \ell}$ as a $\mathring{\U}(\e)$-module.
Define a non-degenerate symmetric bilinear form on $\mathcal{W}_{\epsilon}(x)$ by 
\begin{equation}\label{eq:bilinear form}
(\ket{\mathbf{m}},\ket{\mathbf{m}^{\prime}})=\delta_{\mathbf{m},\mathbf{m}^{\prime}}q^{-\sum_{i=1}^{n}\frac{m_{i}(m_{i}-1)}{2}}\prod_{i=1}^{n}[m_{i}]! 
\end{equation}
for $\ket{\mathbf{m}},\ket{\mathbf{m}^{\prime}}\in \mathcal{W}_{\epsilon}(x)$.
Define an anti-involution $\eta$ on $\mathring{\mathcal{U}}(\epsilon)$
by
\begin{equation*}
\begin{split}
\eta(k_{\mu})&=k_{\mu},\\
\eta(e_{i})&=
\begin{cases}
(-q^{2})^{\epsilon_{i}-\epsilon_{i+1}}q_{i}f_{i}k_{\alpha_{i}}^{-1} & \text{if }i<r\\
(-q^{2})^{\epsilon_{r}-1}q_{r}f_{r}k_{\alpha_{r}}^{-1} & \text{if }i=r\\
q_{i}f_{i}k_{\alpha_{i}}^{-1} & \text{if }i>r
\end{cases},\\
\eta(f_{i})&=
\begin{cases}
(-q^{2})^{\epsilon_{i+1}-\epsilon_{i}}q_{i}^{-1}k_{\alpha_{i}}e_{i} & \text{if }i<r\\
(-q^{2})^{1-\epsilon_{r}}q_{r}^{-1}k_{\alpha_{r}}e_{r} & \text{if }i=r\\
q_{i}^{-1}k_{\alpha_{i}}e_{i} & \text{if }i>r
\end{cases}, 
\end{split}
\end{equation*}
for $\mu\in P$ and $i\in I\setminus\{0\}$. 
Then we have $(\eta\otimes\eta)\circ\Delta=\Delta\circ\eta$.

\begin{lem}\label{lem:pol and semisimple}
We have the following.
\begin{itemize}
 \item[(1)] The bilinear form \eqref{eq:bilinear form} is a polarization, that is, $(x v,w)=(v,\eta(x)w)$ 
for $x\in \mathring{\U}(\e)$ and $v,w\in \W_{\e}(x)$.

\item[(2)] For $\ell\ge 1$, $\W_{\e}(x)^{\ot \ell}$ is semisimple as a $\mathring{\U}(\e)$-module. 
\end{itemize}
\end{lem}
\pf (1) It can be checked in a straightforward way.
 (2) Consider an $A_\infty$-lattice $\mathcal{L}_\infty$ of $\mathcal{W}_\epsilon(x)$
 \[
 \mathcal{L}_\infty=\bigoplus_{\mathbf{m}\in\mathbb{Z}_{+}^{n}(\epsilon)} A_\infty \ket{\mathbf{m}},
 \]
 where $A_\infty$ is the subring of $\mathbb{Q}(q)$ consisting of rational functions regular at $q=\infty$. 
 By the definition of the bilinear form, we have $(\mathcal{L}_{\infty},\mathcal{L}_{\infty})\subset A_{\infty}$, and we obtain the $\mathbb{Q}$-valued induced form on $\mathcal{L}_{\infty}/q^{-1}\mathcal{L}_{\infty}$. Note that the induced form is positive-definite because the image of $\{\ket{\mathbf{m}}\,|\,\mathbf{m}\in\mathbb{Z}_{+}^{n}(\epsilon)\}$
in $\mathcal{L}_{\infty}/q^{-1}\mathcal{L}_{\infty}$ is an orthonormal
basis. Now we can follow the same arguments in \cite[Theorem~2.12]{BKK} to conclude that $\W_{\e}(x)^{\ot \ell}$ is semisimple.
\qed\newline

A decomposition of $\W_\e(x)^{\ot \ell}$ as a $\mathring{\U}(\e)$-module will be given in Section \ref{subsec:decomp of W e}. For this, we introduce an analogue of $\mc{V}^\la$.
%For $\gamma\in P$, let $V_{q,\e}(\gamma)$ be the irreducible highest weight $\mathring{\U}(\e)$-module with highest weight $\gamma$. 

For $\la\in \mc{P}(GL_\ell)$,
we define $\La_{\la,\e}$ in a similar way as in \eqref{eq:Lambda highest weight-2} (see also Remark \ref{rem:comb rule for h.w.}), but now taking $\e$ into account as follows:  

\begin{itemize}
\item[(1)]  
 If $\epsilon_{r+1} =0$ (resp. $\e_{r+1}=1$), then fill the first row (resp. column) of $\la^+$ with $r+1$. After filling a subdiagram $\mu$ of $\lambda^+$ with $r+1,\dots,r+k$, fill the first row (resp. column) of $\la^+/\mu$ with $r+k+1$ if $\epsilon_{r+k+1} =0$ (resp. $\e_{r+k+1}=1$). 
\item[(2)] Fill $\la^-$ in the same way as in $\la^+$ with $r, r-1, \dots$.

\item[(3)] Let $m_i$ be the number of occurrences of $i$'s in $\la^\pm$. Define
\begin{equation*}
\La_{\la,\e}
= \ell \La_{r,\e} -\sum_{i\in \I^-} m_i\de_{i} + \sum_{j\in \I^+} m_j\de_{j}.
\end{equation*}
\end{itemize} 
Depending on $\epsilon$, the weight $\La_{\la,\e}$ may not be well-defined.
Indeed, for $\la\in \mc{P}(GL_\ell)$, $\La_{\la,\e}$ is well-defined if and only if $\la^+\in \cP_{M^+|N^+}$ and $\la^- \in \cP_{M^-|N^-}$ (cf.~Section \ref{subsec:poly repn}).
We put
\begin{equation*}\label{eq:hw osc for e}
\begin{split}
\mc{P}(GL_\ell)_{(\e_-,\e_+)}&=\{\,\la\in \mc{P}(GL_\ell)\,|\,\la^\pm\in \cP_{M^\pm|N^\pm}\,\}.
\end{split}
\end{equation*}

For $\la\in \mc{P}(GL_\ell)_{(\e_-,\e_+)}$, we let $\mc{V}^{\la}_\e$ be the irreducible highest weight $\mathring{\U}(\e)$-module with highest weight $\La_{\la,\e}$.
In particular, we have 
\begin{equation}\label{eq:GL e}
\mc{P}(GL_\ell)_{(\e_-,\e_+)}=
\begin{cases}
 \mc{P}(GL_\ell)_{(r,n-r)} & \text{ if $\e_i=0$ for all $i\in \I$},\\
 \{\,\la\in \mc{P}(GL_\ell)\,|\,\ n-r\ge \la_1\ge \dots\ge \la_\ell\ge -r\,\} &  \text{ if $\e_i=1$ for all $i\in \I$}.
\end{cases}
\end{equation}
 
\subsection{Lie superalgebras and classical limits}
To describe the decomposition of $\W_\e(x)^{\ot \ell}$, we use some results from the theory of super duality \cite{CL}. For this, we need to consider the quantum superalgebra associated to a general linear Lie superalgebra \cite{Ya94} in order to have a well-defined classical limit of $\mc{V}_{\e}^\la$.

Let $\mathring{U}(\e)$ be the associative $\Bbbk$-algebra with $1$ 
generated by $K_\mu, E_i, F_i$ for $\mu\in P$ and $i\in I\setminus\{0\}$ 
satisfying
{\allowdisplaybreaks
\begin{gather*}
K_0=1, \quad K_{\mu +\mu'}=K_{\mu}K_{\mu'} \quad  (\mu, \mu' \in P),\\ 
 K_\mu E_i K_\mu^{-1}=q^{(\mu|\alpha_i)}E_i,\quad 
 K_\mu F_i K_\mu^{-1}=q^{-(\mu|\alpha_i)}F_i, \\ 
 E_iF_j - (-1)^{p(i)p(j)} F_jE_i =
{(-1)^{\e_i}}\delta_{ij}\frac{K_{{\alpha}_i} - K_{-{\alpha}_i}}{q-q^{-1}}, \\
 E_i^2 = F_i^2=0 \quad (i\in I_{\rm odd}\setminus\{0\}), \label{eq:Serre-rel-1-Ya} \\
 E_i E_j - (-1)^{p(i)p(j)} E_j E_i = F_i F_j - (-1)^{p(i)p(j)}  F_j F_i =0\quad (\text{$i-j\not\equiv 1\!\!\pmod{n}$}), \label{eq:Serre-rel-2-Ya} \\ 
\!\!\!
\begin{array}{ll}
E_i^2 E_j- [2] E_i E_j E_i + E_j E_i^2= 0\\ F_i^2 F_j-[2] F_i F_j F_i+F_j F_i^2= 0
\end{array}
\quad (\text{$i\in I_{\rm even}\setminus\{0\}$, $i-j\equiv 1\!\!\pmod{n}$}), \label{eq:Serre-rel-3-Ya}  \\
\!\!\!
\begin{array}{ll}
\left[E_i,\left[\left[E_{i-1},E_i\right]_{(-1)^{p(i-1)}q},E_{i+1}\right]_{(-1)^{(p(i-1)+p(i))p(i+1)}q^{-1}}\right]_{(-1)^{p(i-1)+p(i)+p(i+1)}}=0\\
\left[F_i,\left[\left[F_{i-1},F_i\right]_{(-1)^{p(i)}q},F_{i+1}\right]_{(-1)^{(p(i-1)+p(i))p(i+1)}q^{-1}}\right]_{(-1)^{p(i-1)+p(i)+p(i+1)}}=0 \end{array}\\
\hskip 10cm (i\in I_{\rm odd}\setminus\{0\}),\label{eq:Serre-rel-4-Ya} 
\end{gather*}
\noindent where
$p(i)=\e_i+\e_{i+1}$ $(i\in I)$, and $[X,Y]_t = XY-tYX$ for $t\in \Bbbk^\times$.}

Let $\Sigma$ be the bialgebra over $\Bbbk$ generated by $\sigma_j$ for $j\in \mathbb{I}$ such that $\sigma_i\sigma_j=\sigma_j\sigma_i$ and  $\sigma_j^2=1$ for all $i,j\in \I$, where the comultiplication is given by $\Delta(\sigma_j)=\sigma_j\otimes \sigma_j$ for $j\in \mathbb{I}$.
Then $\mathring{U}(\e)$ is a $\Sigma$-module algebra, where $\Sigma$ acts on $\mathring{U}(\e)$ by 
\begin{equation*}\label{eq:sigma-rel}
\begin{split}
&\sigma_j K_\mu =K_\mu,\quad
\sigma_jE_i=(-1)^{\e_j(\delta_j|\alpha_i)}E_i,\quad 
\sigma_jF_i=(-1)^{\e_j(\delta_j|\alpha_i)}F_i,
\end{split}
\end{equation*}
for $j\in \mathbb{I}$, $\mu\in P$ and $i\in I\setminus{0}$.
Let $\mathring{U}(\e)[\sigma]$ be the semidirect product of $\mathring{U}(\e)$ and $\Sigma$.
We define $\mathring{\mathcal{U}}(\e)[\sigma]$ in the same way.

\begin{lem}[{\cite[Theorem 2.7]{KL}}]
There exists an isomorphism of algebras $\tau : \mathring{U}(\e)[\sigma] \longrightarrow \mathring{\U}(\e)[\sigma]$ such that
\begin{equation*}
\begin{split}
 \tau(E_i)&=e_i \sigma_{\le i}(-\varsigma_i)^{i-1},\ \
 \tau(F_i)=f_i \sigma_{\le i}\varsigma_i^{i},\ \
 \tau(K_{\de_j})=k_{\de_j} \sigma_{j},\ \
 \tau(K_{\La_{r,\e}})=k_{\La_{r,\e}},
\end{split}
\end{equation*}
for $i\in I\setminus\{0\}$ and $j\in \I$, where $\sigma_{\le i}=\sigma_1\cdots\sigma_i$ and $\varsigma_i=\sigma_i\sigma_{i+1}$.
\end{lem}  
 
\begin{rem}\label{rem:super duality}
{\rm 
(1) The isomorphism in \cite[Theorem 2.7]{KL} is given between the quantum {\em affine} superalgebra and $\U(\e)$ extended by $\Sigma$, while the subalgebras generated by  $k_\mu, K_\mu$ $(\mu\in P)$ are slightly different from the ones here. In particular, we have $\tau(K_{\alpha_i})=k_{\alpha_i}\varsigma_i$ for $i\in I\setminus\{0\}$.

(2) Let $\C^{\e}$ be the $n$-dimensional superspace with basis $\{\,v_i\,|\,i\in \I\,\}$. We assume that ${\rm deg} v_i=\epsilon_i$ ($i\in\I$).  We may identify $\gl(\e):={\rm End}_\C(\C^{\e})$ with the Lie superalgebra of matrices $(a_{ij})_{i,j\in\I}$, which is spanned by the elementary matrices  $E_{ij}$ with $1$ at the $(i,j)$-position and zero elsewhere.
The algebra $\mathring{U}(\e)$ can be viewed as the quantized enveloping algebra of one-dimensional central extension of $\gl(\e)$. We may regard $\gl(\ov{\e})$ as a subalgebra of $\gl(\e)$ by identifying $\C^{\ov{\e}}$ as a subspace of $\C^\e$, where $\ov{\e}$ is the subsequence of $\e$ consisting of $0$'s.} 
\end{rem}
 
 Let $V$ be a $\mathring{\U}(\e)$-module. 
One can extend the $\mathring{\mathcal{U}}(\epsilon)$-action to $\mathring{\mathcal{U}}(\epsilon)[\sigma]$ by 
$\sigma_{i}u = (-1)^{\epsilon_i \mu_i}u$ 
for $i\in\mathbb{I}$ and $u\in V_\mu$ with $\mu=\ell\La_{r,\e}+\sum_{i\in \I}\mu_i\de_i$. 
Let $V^\tau=\{\,v^\tau\,|\,v\in V\,\}$ denote the $\mathring{U}(\e)$-module obtained from $V$ by applying $\tau$. 
Suppose that $V$ has weight space decomposition, that is, $V=\bigoplus_{\nu\in P}V_\nu$.

\begin{lem}\label{lem:tau pullback of wt space}
Under the above hypothesis, $V^\tau$ also has a weight space decomposition $V^\tau=\bigoplus_{\nu\in P}\left(V^\tau\right)_\nu$ with
\[
\left(V^\tau\right)_\nu =  \left\{u\in V^\tau \,|\, K_\mu u= q^{(\mu|\nu)}u \ \  (\mu\in P)\right\}.
\]
Indeed, we have
$\left(V^\tau\right)_\nu =\{\,v^\tau\,|\,v\in V_\nu\,\}$.
\end{lem}
\pf Let $v\in V_\nu$ be given with $\nu=\ell'\La_{r,\e}+\sum_{i\in \I}\nu_i\de_i=\ell'\La_{r,\e}+\nu^\circ$. 
For $\mu=\ell\La_{r,\e}+\sum_{i\in \I}\mu_i\de_i=\ell\La_{r,\e}+\mu^\circ$, we have
\begin{equation*}
\begin{split}
 K_\mu v^\tau= (\tau(K_\mu) v)^\tau &
 =(\tau(K_{\mu^\circ})k_{\ell\La_{r,\e}}v)^\tau
 =q^{(\ell\La_{r,\e}|\nu^\circ)}(\tau(K_{\mu^\circ})\,v)^\tau\\
 &=q^{(\ell\La_{r,\e}|\nu^\circ)}\left(k_{\mu^\circ}\prod_{i\in \I}\sigma_i^{\mu_i}\, v\right)^\tau\\
 &=q^{(\ell\La_{r,\e}|\nu^\circ)}q^{(\ell'\La_{r,\e}|\mu^\circ)} \bq(\mu,\nu)\prod_{i\in \I}(-1)^{\e_i\mu_i\nu_i}v^\tau\\
 &=q^{(\ell\La_{r,\e}|\nu^\circ)+(\ell'\La_{r,\e}|\mu^\circ)+(\mu^\circ|\nu^\circ)} v^\tau=q^{(\mu|\nu)}v^\tau.
\end{split}
\end{equation*}
\qed\newline  

For $X=(X_{ij})_{i,j\in \I}\in \gl(\e)$, let ${\rm str}(X)=\sum_{i\in \mathbb{I}}(-1)^{\e_i}X_{ii}$ be the supertrace of $X$. 
Let $\beta$ be a 2-cocycle on $\gl(\e)$ given by $\beta(X,Y)={\rm str}([J,X]Y)$, where $J=\sum_{j\in \I^+}E_{jj}$ (cf.~\cite{CK08}). 
Let $\gl(\e)^e=\gl(\e)\oplus \mathbb{C} c$ be a central extension with respect to $\beta$ such that $c$ is a central element and
\begin{equation*}
 [X,Y]=XY-(-1)^{p(X)p(Y)}YX +\beta(X,Y)c,
\end{equation*}
for homogeneous elements $X, Y$ with parities $p(X)$, $p(Y)$ respectively.
In particular, we have
\begin{equation*}\label{eq:Weyl relation for Lie super}
 [E_{i,i+1}, E_{i+1,i}] = E_{i,i}-(-1)^{\e_i+\e_{i+1}}E_{i+1,i+1} + \de_{ir}c=:\alpha_i^\vee\quad (i\in I\setminus\{0\}).
\end{equation*}

Let $P^\vee=\Z c \oplus \bigoplus_{i\in \I}\Z E_{ii}$. 
If we define a pairing on $P\times P^\vee$ such that 
\begin{equation*}
 \langle \de_i,E_{jj}\rangle=\de_{ij},\quad 
\langle \de_i, c \rangle= \langle \La_{r,\e}, E_{ii} \rangle=0,\quad 
 \langle \La_{r,\e},c \rangle=1\quad (i,j\in \I),
\end{equation*}
then we have an isomorphism of abelian groups $\phi : P \longrightarrow P^\vee$ sending
\begin{equation*}
\begin{split}
 \phi(\de_i)&= D_i:=(-1)^{\e_i}E_{ii} + \sigma(i)c\quad (i\in \I),\\
 \phi(\La_{r,\e})&= D^+:=\sum_{j\in \I^+}E_{jj},
\end{split}
\end{equation*}
where $\sigma(j)=0$ (resp. $1$) for $i\in \I^-$ (resp. $i\in \I^+$), such that 
\begin{equation*}
 (\la|\mu)=\langle \la,\phi(\mu)\rangle\quad (\la,\mu\in P).
\end{equation*} 
Note that $(\la|\alpha_i)=(-1)^{\e_i}\langle \la,\alpha_i^\vee \rangle$ for $i\in I\setminus\{0\}$.

Suppose that $V$ is a highest weight $\mathring{\U}(\e)$-module with highest weight $\La_{\la,\e}$ for some $\la\in \mc{P}(GL_\ell)_{(\e_-,\e_+)}$.
Consider the $\mathring{U}(\e)$-module $V^\tau$ and let
$\overline{V^\tau}$ be its classical limit defined as in \eqref{eq:classical limit}. 
Let ${\rm E}_i$, ${\rm F}_i$, ${\rm D}_j$, and ${\rm D}_{\La_{r,\e}}$ be the operators induced from $E_i$, $F_i$, $\{K_{\de_j}\}$, and $\{K_{\La_{r,\e}}\}$ in $\mathring{U}(\e)$, respectively for $i\in I\setminus\{0\}$ and $j\in \I$.
Let $\ov{U}$ be the subalgebra of ${\rm End}_{\mathbb{C}}\left(\ov{V^\tau}\right)$ generated by ${\rm E}_i$, ${\rm F}_i$, ${\rm D}_j$, and ${\rm D}_{\La_{r,\e}}$. 

\begin{lem}\label{lem:classical limit to gl}
Under the above hypothesis, there exists a homomorphism of $\C$-algebras from $U(\gl(\e)^e)$ to $\ov{U}$ sending 
\begin{equation*}\label{eq:map on classical limit}
\begin{split}
&E_{i,i+1}\ \longmapsto\ {\rm E}_i,\quad
E_{i+1,i}\ \longmapsto\ {\rm F}_i \quad (i\in I\setminus\{0\}),\\
&D^+ \longmapsto {\rm D}_{\La_{r,\e}},\quad
D_j \ \longmapsto\ {\rm D}_j \quad (j\in \I). 
\end{split}
\end{equation*}
\end{lem}
\pf The defining relations among $E_{i,i+1}$'s (resp. $E_{i+1,i}$'s) follow from \cite[Theorem 10.5.8]{Ya94}.
So it remains to show that the following relations in $U(\gl(\e)^e)$ are preserved under the above correspondence:
\begin{gather*}
 [D,D']=0,\ \
 [E_{i,i+1},E_{j+1,j}]= \de_{ij}\alpha_i^\vee,\\
 [D,E_{i,i+1}]=\langle \alpha_i,D \rangle E_{i,i+1},\quad
 [D,E_{i+1,i}]=-\langle \alpha_i,D \rangle E_{i+1,i},
\end{gather*}
for $i,j\in I\setminus\{0\}$ and $D, D'\in \C D^++\sum_{i\in \mathbb{I}}\C D_i$.

Let us check the second relation for example.
For $u\in (V^\tau)_\nu$, let $\ov{u}$ denote its image in $\ov{V^\tau}$. 
Then by Lemma \ref{lem:tau pullback of wt space}
\begin{equation*}
\begin{split}
 [{\rm E}_i,{\rm F}_i]\cdot \ov{u} 
 &= \ov{[E_i,F_i]\cdot u} = \ov{(-1)^{\e_i}\{K_{\alpha_i}\}\cdot u}
 = \ov{(-1)^{\e_i}[(\nu|\alpha_i)]u} \\
 &= (-1)^{\e_i}(\nu|\alpha_i)\ov{u}
 = \langle \nu, \alpha_i^\vee\rangle u
 = \alpha_i^\vee\cdot \ov{u}.
\end{split}
\end{equation*}
We can also check the other relations directly.
\qed

\begin{cor}\label{cor:classical limit of h.w. mod}
$\ov{V^\tau}$ is a highest weight $U(\gl(\e)^e)$-module with highest weight $\La_{\la,\e}$.
\end{cor}
\pf By Lemma \ref{lem:classical limit to gl}, $\ov{V^\tau}$ is a highest weight $U(\gl(\e)^e)$-module. Let $u$ be a weight vector of $\ov{V^\tau}$ with weight $\nu$. For $\mu\in P$, we have
\begin{equation*}
 \phi(\mu)\cdot\ov{u}
 =\ov{\{K_{\mu}\}\cdot u}
 =\ov{[(\mu|\nu)] u}
 =(\mu|\nu)\ov{u}
 =\langle \nu,\phi(\mu) \rangle \ov{u}.
\end{equation*}
This implies that $\ov{u}$ is of weight $\nu$ with respect to the action of the Cartan subalgeba of $\gl(\e)^e$. In particular, the highest weight of $\ov{V^\tau}$ is $\La_{\la,\e}$.
\qed

\subsection{Reduction homomorphism}\label{subsec:reduction homo}
Suppose that $n\ge 5$.
Let $\e'=(\e'_1,\dots,\e'_{n-1})$ be the sequence obtained from $\e$ by removing $\e_i$ for some $i\in \mathbb{I}$.
Let $\mathbb{I}'=\{\,1,\dots,n-1\,\}$ with the $\mathbb{Z}_2$-grading induced from $\mathbb{I}$. Let $\mathbb{I}_\pm'$ be the subsequences of $\mathbb{I}'$ induced from $\mathbb{I}_\pm$ and let $r'=|\mathbb{I}'_-|$.

We assume that the weight lattice for $\mathcal{U}(\e')$ is $P'=\Z\La_{r,\epsilon^\prime}\oplus\bigoplus_{l\in \mathbb{I}'}\mathbb{Z}\delta'_l$.
Put $I'=\{0,1,\cdots,n-2\}$.  
Denote by $k'_{\mu}$, $e'_j$, and $f'_j$ the generators of $\mathcal{U}(\e')$ for $\mu\in P'$ and $j\in I'$. 
Then we define $\hat{k}_{\delta'_l}, \hat{e}_j, \hat{f}_j\in \mathcal{U}(\e)$ for $l\in \mathbb{I}'$ and $j \in I'$ as follows: 
\begin{equation*}
\hat{k}_{\Lambda_{r',\epsilon^\prime}}=k_{\Lambda_{r,\epsilon}}, \quad \hat{k}_{\delta'_l}=
\begin{cases}
k_{\delta_l} & \text{for $1\leq l\leq i-1$}, \\
k_{\delta_{l+1}} & \text{for $i \leq l\leq n-1$,}
\end{cases}
\end{equation*}

{\allowdisplaybreaks
{\em Case 1}. If $2\leq i\leq r$, then
\begin{gather*}\label{eq:generators for e'}
(\hat{e}_j, \hat{f}_j) =
\begin{cases}
\left(e_j,f_j\right) & \text{for $j\leq i-2$},\\
\left([e_{i},e_{i-1}]_{\bq_{i-1,i}^{-1}}, [f_{i-1},f_{i}]_{\bq_{i-1,i}}\right) & \text{for $j= i-1$},\\
\left(e_{j+1},f_{j+1}\right) & \text{for $j\geq i$}.
\end{cases}
\end{gather*}

{\em Case 2}. If $r+1\leq i \leq n-1$, then
\begin{gather*}\label{eq:generators for e'-2}
(\hat{e}_j, \hat{f}_j) =
\begin{cases}
\left(e_j,f_j\right) & \text{for $j\leq i-2$},\\
\left([e_{i-1},e_{i}]_{\bq_{i-1,i}}, [f_{i},f_{i-1}]_{\bq_{i-1,i}^{-1}}\right) & \text{for $j= i-1$},\\
\left(e_{j+1},f_{j+1}\right) & \text{for $j\geq i$}.
\end{cases}
\end{gather*}

{\em Case 3}. If $i=n$, then  
\begin{gather*}\label{eq:generators for e'-3}
(\hat{e}_j, \hat{f}_j) =
\begin{cases}
(e_j,f_j) & \text{for $j\neq 0$},\\
\left([e_{n-1},e_{0}]_{\bq_{n-1,0}},[f_{0},f_{n-1}]_{\bq_{n-1,0}^{-1}}\right) & \text{for $j=0$}.
\end{cases}
\end{gather*}

{\em Case 4}. If $i=1$, then
\begin{gather*}\label{eq:generators for e'-4}
(\hat{e}_j, \hat{f}_j) =
\begin{cases}
\left([e_{1},e_{0}]_{\bq^{-1}_{0,1}},[f_{0},f_{1}]_{\bq_{0,1}}\right) & \text{for $j=0$},\\
(e_{j+1},f_{j+1}) & \text{for $j\neq 0$}.
\end{cases}
\end{gather*}}
Here $\bq_{a,b}=\bq(\alpha_a,\alpha_b)$ for $a,b\in I$, and $[X,Y]_t = XY-tYX$ for $X,\,Y\in\mathcal{U}(\epsilon)$, $t\in\Bbbk^\times$.

\begin{thm}[{\cite[Theorem 4.3]{KY}}]\label{thm:folding homomorphism}
Under the above hypothesis, there exists a homomorphism of $\Bbbk$-algebras 
\begin{equation*}
\xymatrixcolsep{2pc}\xymatrixrowsep{3pc}
\xymatrix{
\phi^\e_{\e'} : \mathcal{U}(\e') \ \ar@{->}[r] & \ \mathcal{U}(\e),
}
\end{equation*}
such that for $l\in \mathbb{I}'$ and $j\in I'$
\begin{equation*}
\phi^\e_{\e'}(k'_{\La_{r',\e'}})=\hat{k}_{\La_{r',\e'}},\quad
\phi^\e_{\e'}(k'_{\delta'_l})=\hat{k}_{\delta'_l},\quad 
\phi^\e_{\e'}(e'_j)=\hat{e}_j,\quad \phi^\e_{\e'}(f'_j)=\hat{f}_j.
\end{equation*}
\end{thm}

More generally, suppose that the sequence $\e'=(\e'_1,\dots,\e'_{n-t})$ with $1\le t\le n-3$ is obtained from $\e$ by removing $\e_{i_1},\dots,\e_{i_t}$ for some $i_1<\dots<i_t$ .
For $0\leq k\leq t$, let $\e^{(k)}$ be a sequence such that
\begin{itemize}
\item[(1)] $\e^{(0)}=\e$, $\e^{(t)}=\e'$, 

\item[(2)] $\e^{(k)}$ is obtained from $\e^{(k-1)}$ by removing $\e_{i_k}$ for $1\leq k\leq t$.
\end{itemize}
We define a homomorphism of $\Bbbk$-algebras
\begin{equation}\label{eq:reduction homo}
\xymatrixcolsep{2pc}\xymatrixrowsep{3pc}
\xymatrix{
\phi^\e_{\e'} : \mathcal{U}(\e') \ \ar@{->}[r] & \ \mathcal{U}(\e),
}
\end{equation}
by $\phi^{\e}_{\e'}=\phi^{\e^{(0)}}_{\e^{(1)}}\circ \phi^{\e^{(1)}}_{\e^{(2)}}\circ \dots \circ\phi^{\e^{(t-1)}}_{\e^{(t)}}$.

\subsection{Truncation functor}\label{subsec:truncation functor}
Let 
\begin{equation*}
P_{\ge 0}=\Z\La_{r,\e} - \sum_{i\in \I^-}\Z_+\de_i + \sum_{j\in \I^+}\Z_+\de_j.
\end{equation*}
Let 

\begin{itemize}
 \item[$\bullet$] $\mc{C}(\e)$ : the category of $\mathcal{U}(\e)$-modules $V$ such that $V=\bigoplus_{\mu\in P} V_\lambda$ and ${\rm wt}(V)\subset P_{\ge 0}$,
 
 \item[$\bullet$] $\mc{C}^\ell(\e)$ ($\ell\in \mathbb{Z}$) : the full subcategory of $\mc{C}(\e)$ consisting of $V$ such that $\la \in \ell\La_{r,\e}- \sum_{i\in \I^-}\Z_+\de_i + \sum_{j\in \I^+}\Z_+\de_j$ for $\la\in {\rm wt}(V)$.
 
\end{itemize}
Then we have $\mc{C}(\e)=\bigoplus_{\ell\in \mathbb{Z}}\mc{C}^\ell(\e)$. 
Note that $\mc{C}(\e)$ is closed under taking submodules, quotients and tensor products, while $\mc{C}^\ell(\e)$ is closed under taking submodules and quotients.
We define $\mathring{\mc{C}}(\e)$ for $\mathring{\U}(\e)$-modules and $\mathring{\mc C}(\e_-,\e_+)$ for $\mathring{\U}(\e_-,\e_+)$-modules in the same way and $\mathring{\mc{C}}^{\ell}(\e)=\mathring{\mc{C}}(\e)\cap \mc{C}^\ell(\e)$.
We have $\W_\e(x)^{\ot \ell}\in \mc{C}^{\ell}(\e)$ for $\ell\ge 1$ by \eqref{eq:weight of m e}.

Now, let $\e=(\e_1,\dots,\e_n)$ and $\e'=(\e'_1,\dots,\e'_{n'})$ be as in \eqref{eq:reduction homo}, where $n'=n-t$. Let $\e=(\e_-,\e_+)$ be as in \eqref{eq:e for parabolic}, and let $\e'=(\e'_-,\e'_+)$ be induced from $(\e_-,\e_+)$ with $\e'_-=(\e'_1,\dots,\e'_{r'})$ and $\e'_+=(\e'_{r'+1},\dots,\e'_{n'})$ satisfying $r'\ge 2$ and $n'-r'\ge 2$. We define $\mc{C}(\e')$, $\mc{C}^\ell(\e')$, $\mathring{\mc{C}}(\e')$, and $\mathring{\mc{C}}^{\ell}(\e')$ in the same way.

For a $\mathcal{U}(\e)$-module $V$ in $\mc{C}(\e)$, define
\begin{equation}\label{eq:truncation-1}
\mf{tr}^\e_{\e'}(V) = 
\bigoplus_{\substack{\mu\in {\rm wt}(V) \\ ({\rm pr}_{\I}(\mu)|\de_{i_1})=\dots=({\rm pr}_{\I}(\mu)|\de_{i_t})=0}}V_\mu,
\end{equation}
where ${\rm pr}_{\I}: P_{\ge 0} \longrightarrow \bigoplus_{i\in \I}\Z\de_i$ denotes the canonical projection.
We denote by $\pi^\e_{\e'} : V \longrightarrow \mf{tr}^\e_{\e'}(V)$ the natural projection. For any $\mathcal{U}(\e)$-modules $V, W$ in $\mc{C}(\e)$ and $f \in {\rm Hom}_{\mathcal{U}(\e)}(V,W)$, let
\begin{equation*}
\xymatrixcolsep{2pc}\xymatrixrowsep{3pc}
\xymatrix{
\mf{tr}^\e_{\e'}(f) : \mf{tr}^\e_{\e'}(V) \ \ar@{->}[r] &  \mf{tr}^\e_{\e'}(W)
}
\end{equation*}
be the $\Bbbk$-linear map given by $\mf{tr}^\e_{\e'}(f)(v)=f(v)$ for $v\in \mf{tr}^\e_{\e'}(V)$. Then we have the following commutative diagram of $\Bbbk$-vector spaces:
\begin{equation*}\label{eq:truncation commuting diagram}
\xymatrixcolsep{4pc}\xymatrixrowsep{3pc}\xymatrix{
V \ar@{->}[d]_{\pi^\e_{\e'}}\ar@{->}[r]^{f} &  W \ar@{->}[d]^{\pi^\e_{\e'}} \\
\mf{tr}^\e_{\e'}(V) \ar@{->}[r]^{\mf{tr}^\e_{\e'}(f)} & \mf{tr}^\e_{\e'}(W)
}
\end{equation*}
\begin{prop}[{\cite[Propositions 4.4]{KY}}]\label{prop:truncation}
Under the above hypothesis,
\begin{itemize}
\item[(1)] the space $\mf{tr}^\e_{\e'}(V)$ is invariant under the action of $\mathcal{U}(\e')$ via $\phi^{\e}_{\e'}$, and hence a $\mathcal{U}(\e')$-module in $\mc{C}(\e')$,

\item[(2)] the map $\mf{tr}^\e_{\e'}(f) : \mf{tr}^\e_{\e'}(V) \longrightarrow \mf{tr}^\e_{\e'}(W)$ is $\mathcal{U}(\e')$-linear,

\item[(3)] the space $\mf{tr}^\e_{\e'}(V\ot W)$ is naturally isomorphic to the tensor product of $\mf{tr}^\e_{\e'}(V)$ and $\mf{tr}^\e_{\e'}(W)$ as a $\mathcal{U}(\e')$-module.

\end{itemize}
\end{prop}

Hence we have a functor 
\begin{equation*}
\xymatrixcolsep{2pc}\xymatrixrowsep{3pc}
\xymatrix{
\mf{tr}^\e_{\e'} : \mc{C}(\e)  \ar@{->}[r] &  \mc{C}(\e'),}
\end{equation*}
which we call {\em truncation}. Note that $\mf{tr}^\e_{\e'}$ is exact by its definition and monoidal by Proposition \ref{prop:truncation} in the sense of \cite[Appendix A.1]{KKK}.
We may also define $\mf{tr}^\e_{\e'} : \ring{\mc{C}}(\e) \longrightarrow \ring{\mc{C}}(\e')$ and $\mf{tr}^\e_{\e'} : \ring{\mc{C}}(\epsilon_-,\epsilon_+) \longrightarrow \ring{\mc{C}}(\epsilon'_-,\epsilon'_+)$ in the same way.

Let $M'$ and $N'$ be the numbers of $0$'s and $1$'s in $\e'$, and let $M'_\pm$ and $N'_\pm$ be the numbers of $0$'s and $1$'s in $\e'_\pm$.
Then we can check the following.
\begin{lem}\label{lem:truncation of fund}
For $l\in\mathbb{Z}$ and $x\in \Bbbk^\times$, we have as a $\U(\e')$-module
\begin{equation*}
\mf{tr}^\e_{\e'}(\W_{l,\e}(x))\cong
\begin{cases}
\W_{l,\e'}(x) & \text{if $(l)\in \cP_{M'|N'}$},\\
0 & \text{otherwise}.
\end{cases}
\end{equation*}
In particular, we have $\mf{tr}^\e_{\e'}\left(\W_{\e}(x)^{\ot\ell}\right)\cong \W_{\e'}(x)^{\ot\ell}$ for $\ell\ge 1$.
\end{lem}

We remark that the action of $\U(\e')$ on $\W_{\e'}(x)^{\ot\ell}$ does not depend on the choice of the sequence $\{\e^{(t)}\}_{1\le k\le t}$ in \eqref{eq:reduction homo}.

\begin{lem}\label{lem:truncation of poly}
For $\mu\in \cP_{M_-|N_-}$ and $\nu\in \cP_{M_+|N_+}$, we have as a $\mathring{\U}(\e'_-,\e'_+)$-module
\begin{equation*}
\mf{tr}^\e_{\e'}\left(V_{(\e_-,\e_+)}(\mu,\nu)\right)\cong
\begin{cases}
V_{(\e'_-,\e'_+)}(\mu,\nu)& \text{if $\mu\in \cP_{M'_-|N'_-}$ and $\nu\in \cP_{M'_+|N'_+}$},\\
0 & \text{otherwise}.
\end{cases}
\end{equation*}
\end{lem}
\pf By \cite[Proposition 4.5]{KY}, we have 
$$\mf{tr}^{\e_\pm}_{\e'_\pm}(V_{\e_\pm}(\pm\la))\cong V_{\e'_\pm}(\pm\la),$$
for $\la\in\cP_{M_\pm|N_\pm}$, where we assume that the righthand side is zero unless $\la\in\cP_{M'_\pm|N'_\pm}$.
This proves the assertion.
\qed\newline

\subsection{Decomposition of $\W_{\e}(x)^{\ot \ell}$}\label{subsec:decomp of W e}
In this subsection, we consider a special sequence $\e$ and its homogeneous subsequences.
For $a\ge 1$, let 
\begin{gather*}
\e^{(a)}=(\underbrace{0,1,0,1,\dots,0,1}_{2a}),\\
\underline{\e}^{(a)}=(\underbrace{0,\dots,0}_{a}),\quad 
\overline{\e}^{(a)}=(\underbrace{1,\dots,1}_{a}).
\end{gather*}
Suppose that $n=2a+2b+1$ for some $a,b\ge 1$ with $r=2a$ and let 
\begin{equation}\label{eq:e(a,b)}
 \e=\e^{(a,b)}:=(\e^{(a)},\e^{(b)},0)
\end{equation}
be the concatenation of $\e^{(a)}$, $\e^{(b)}$ and $0$ such that $\e_-=\e^{(a)}$ and $\e_+=(\e^{(b)},0)$.
Let $\underline{\e}$ and $\overline{\e}$ be the sequences of length $n'=a+b+1$ and $a+b$ respectively, with $r'=r/2$ given by
\begin{equation}\label{eq:e'(a,b)}
\underline{\e}=(\underline{\e}^{(a)},\underline{\e}^{(b)},0),\quad
\overline{\e}=(\overline{\e}^{(a)},\overline{\e}^{(b)}).  
\end{equation}
Note that $\U(\underline{\e})\cong U_q(\widehat{\mf{gl}}_{a+b+1})$ and $\U(\overline{\e})\cong U_{\td{q}}(\widehat{\mf{gl}}_{a+b})$ where $\td{q}=-q^{-1}$.

For $k\ge 1$, let us identify $\e^{(a,b)}$ with the subsequence of $\e^{(a+k,b+k)}$ by removing the first and the last $2k$ entries in $\e^{(a+k)}$ and $(\e^{(b+k)},0)$ respectively. Let $\mf{tr}^{(a+k,b+k)}_{(a,b)}$ denote the truncation functor associated to them.

\begin{prop}\label{prop:hwsub-W^l-super}
Let $\ell\ge 1$ and $x\in \Bbbk^\times$ be given.
Any highest weight $\mathring{\U}(\e)$-submodule in $\W_\e(x)^{\ot \ell}$ is isomorphic to $\mc{V}^{\la}_\e$ for some $\la\in \mc{P}(GL_\ell)_{(\e_-,\e_+)}$.
\end{prop}
\pf 
Suppose that $V$ is a highest weight $\mathring{\U}(\e)$-submodule of $\mathcal{W}_\epsilon(x)^{\otimes\ell}$ and hence irreducible by Lemma \ref{lem:pol and semisimple}(2). Let $v$ be a highest weight vector of $V$ with highest weight $\gamma$. Note that $\gamma=\ell\La_{r,\e}+\sum_{i\in \mathbb{I}}\gamma_i\de_i\in P_{\ge 0}$ for some $\gamma_i\in\Z$. 

We may identify $\W_{\e}(x)$ as a subspace of $\W_{\td{\e}}(x)$ with the entries in $\ket{{\bf m}}$ outside $\e^{(a,b)}$ being zero, where $\td{\e}=\e^{(a+k,b+k)}$ for $k\ge 1$.
Then we may assume that $v$ is a vector in $\W_{\td{\e}}(x)^{\ot \ell}$, for all $k\ge 1$. Since $\W_{\td{\e}}(x)^{\ot \ell}\in \mathring{\mc{C}}(\td{\e})$, it is not difficult to see from Proposition \ref{prop:osc module-3}(1) that $v$ is also a highest weight vector with respect to $\mathring{\U}\left(\e^{(a+k,b+k)}\right)$ for $k\ge 1$. 
Let $\widetilde{V}$ be the irreducible highest weight $\mathring{\U}(\td{\e})$-module generated by $v$ in  $\W_{\td{\e}}(x)^{\ot \ell}$.
So we may assume that $a$ and $b$ are large enough.

Consider $\mf{tr}^\e_{\ov{\e}}(V)$. 
By Lemma \ref{lem:semisimple over levi}, we have $v\in V_{(\e_-,\e_+)}(\mu,\nu)$ for some $\mu\in \cP_{a|a}$ and $\nu\in \cP_{b+1|b}$.
Note that for any other $\mathring{\U}(\e_-,\e_+)$-component $V_{(\e_-,\e_+)}(\zeta,\xi)$ of $V$, we have $|\mu|\le |\zeta|$ and $|\nu|\le |\xi|$, where $|\la|$ denotes the size of a partition $\la$ since the action of $f_r$ maps $V_{(\e_-,\e_+)}(\alpha,\beta)$ to $V_{(\e_-,\e_+)}(\gamma,\delta)$ with $|\gamma|=|\alpha|+1$ and $|\delta|=|\beta|+1$.

Since we assume that $a$ and $b$ are large enough, we have by Lemma \ref{lem:truncation of poly}
\begin{equation}\label{eq:tr preserves poly}
\mf{tr}^\e_{\ov{\e}}(V_{(\e_-,\e_+)}(\mu,\nu))=V_{(\ov{\e}_-,\ov{\e}_+)}(\mu,\nu)\neq 0.
\end{equation}
Let $\ov{w}$ be the $\mathring{\U}(\ov{\e}_-,\ov{\e}_+)$-highest weight vector in $V_{(\ov{\e}_-,\ov{\e}_+)}(\mu,\nu)$. By minimality of $|\mu|$ and $|\nu|$, $\ov{w}$ is a $\mathring{\U}(\ov{\e})$-highest weight vector with highest weight, say $\ov{\gamma}$.
Since $\mf{tr}^\e_{\ov{\e}}(\W_\e(x)^{\ot \ell})=\W_{\ov{\e}}(x)^{\ot \ell}$ is 
finite-dimensional (see Remark \ref{rem:odd homogeneous}), we have 
\begin{equation*}
\ov{\gamma} = \La_{\la,\ov{\e}},
\end{equation*}
for some $\la\in \mc{P}(GL_\ell)_{(\overline{\epsilon}_{-},\overline{\epsilon}_{+})}$. This implies that $\gamma = \La_{\la,\e}$ by \eqref{eq:tr preserves poly}.
\qed 

\begin{lem}\label{lem:truncation by rank}%\comments{It has been added (2022.2.27).}\todo{Revised for the pf of Thm6.17 (Lee)}
Let $\td{\e}=\e^{(a+k,b+k)}$ for $k\ge 1$. 
\begin{itemize}
 \item[(1)] For $\la\in \mc{P}(GL_\ell)_{(\e_-,\e_+)}$ such that $\V^\la_\e\subset \W_\e(x)^{\ot \ell}$, we have $\V^\la_{\td{\e}}\subset \W_{\td{\e}}(x)^{\ot \ell}$.
 
 \item[(2)] For $\la\in \mc{P}(GL_\ell)_{(\td{\e}_-,\td{\e}_+)}$ such that $\V^\la_{\td{\e}}\subset \W_{\td{\e}}(x)^{\ot \ell}$, we have 
\[
 \mf{tr}^{\td{\e}}_\e\left(\V^\la_{\td{\e}}\right)=
 \begin{cases}
  \mathcal{V}^{\lambda}_\epsilon & \text{if } \la \in \mc{P}(GL_\ell)_{(\e_-,\e_+)}, \\
  0 & \text{if } \la\not\in \mc{P}(GL_\ell)_{(\e_-,\e_+)}.
 \end{cases}
\]
\end{itemize}
\end{lem}

\pf (1) Let $v$ be a $\mathring{\U}(\e)$-highest weight vector of ${\mc{V}_{\e}^\la}$ in $\W_\e(x)^{\ot \ell}$. Then we may regard $v$ as an element in $\W_{\td{\e}}(x)^{\ot \ell}$. As seen in the proof of Proposition \ref{prop:hwsub-W^l-super}, $v$ is also a $\mathring{\U}(\td{\e})$-highest weight vector. By Proposition \ref{prop:hwsub-W^l-super}, the $\mathring{\U}(\td{\e})$-submodule generated by $v$ is isomorphic to $\V^\la_{\td{\e}}$. So we have $\V^\la_{\td{\e}}\subset \W_{\td{\e}}(x)^{\ot \ell}$.

(2) First assume $\la \in \mc{P}(GL_\ell)_{(\e_-,\e_+)}$. If we let $v$ be a $\mathring{\U}(\tilde{\e})$-highest weight vector of ${\mc{V}_{\tilde\e}^\la}$, then by assumption, $v$ belongs to $\mathfrak{tr}^{\tilde \epsilon}_{\epsilon} \left(\mathcal{V}^{\lambda}_{\tilde\epsilon}\right)$ and is of highest weight over $\mathring{\U}(\e)$ as well. Since $\mathfrak{tr}^{\tilde \epsilon}_{\epsilon} \left(\mathcal{V}^{\lambda}_{\tilde\epsilon}\right) \subset \mathfrak{tr}^{\tilde \epsilon}_{\epsilon} \left(\mathcal{W}_{\tilde\epsilon}(x)^{\otimes \ell}\right) = \mathcal{W}_\epsilon (x)^{\otimes \ell}$, we obtain $\mathcal{V}^{\lambda}_\epsilon\subset \mathfrak{tr}^{\tilde \epsilon}_{\epsilon} \left(\mathcal{V}^{\lambda}_{\tilde\epsilon}\right)$.  Suppose that $\V^\la_{\e}\subsetneq \mf{tr}^{\tilde\e}_\e\left(\V^\la_{\td{\e}}\right)$. Since $\mf{tr}^{\tilde\e}_\e\left(\V^\la_{\td{\e}}\right)$ is semisimple over $\mathring{\mathcal{U}}(\epsilon)$, there exists a $\mathring{\U}({\e})$-highest weight vector $w\in \mf{tr}^{\tilde\e}_\e\left(\V^\la_{\td{\e}}\right) \setminus \Bbbk v$ such that ${\rm wt}(w)\in \La_{\la,\e}-\sum_{i\in I\setminus\{0\}}\Z_+\alpha_i$. But then $w$ is also a $\mathring{\U}(\td{\e})$-highest weight vector in $\V^\la_{\td{\e}}$ as in (1), which is a contradiction.

For the other one,  %Let $v$ be a $\mathring{\U}(\td{\e})$-highest weight vector of ${\mc{V}_{\td{\e}}^\la}$ in $\W_{\td{\e}}(x)^{\ot \ell}$. 
Let $\mathbb{I}_{\td{\e}}=\mathbb{I}_{\td{\e}}^-\cup\mathbb{I}_{\td{\e}}^+$ denote the index set $\mathbb{I}$ for $\U(\td{\e})$ where $\mathbb{I}_{\td{\e}}^-=\{\,1,2,\dots,2a+2k\,\}$ and $\mathbb{I}_{\td{\e}}^+=\{\,2a+2k+1,\dots,2a+2b+4k+1\,\}$.
We identify the index set $\mathbb{I}$ for $\U(\e)$ with $\mathbb{I}_\e=\mathbb{I}_{\e}^-\cup\mathbb{I}_{\e}^+$ where $\mathbb{I}_{\e}^-=\{\,2k+1,2k+2,\dots,2a+2k\,\}$ and $\mathbb{I}_{\e}^+=\{\,2a+2k+1,\dots,2a+2b+2k+1\,\}$.

If $\la\not\in \mc{P}(GL_\ell)_{(\e_-,\e_+)}$, then we have $({\rm pr}_{\I_{\td{\e}}}(\La_{\la,\td{\e}})|\de_i)\neq 0$ for some $i\in \mathbb{I}_{\td{\e}}\setminus \mathbb{I}_{\e}$. Since $\W_{\td{\e}}(x)^{\ot \ell}\in \mathring{\mc{C}}(\td{\e})$ and $\V^\la_{\td{\e}}$ is a highest weight module, any weight $\mu$ of $\V^\la_{\td{\e}}$ satisfies $({\rm pr}_{\I_{\td{\e}}}(\mu)|\de_i)\neq 0$ for some $i\in \mathbb{I}_{\td{\e}}\setminus \mathbb{I}_{\e}$. 
Hence $\mf{tr}^{\td{\e}}_\e\left(\V^\la_{\td{\e}}\right)=0$.
\qed

\begin{thm}\label{thm:trunc sends simple to simple}
For $\la\in \mc{P}(GL_\ell)_{(\e_-,\e_+)}$, we have
\begin{equation*}
\mf{tr}^\e_{\e'}\left(\mc{V}^{\la}_\e\right)\cong
\begin{cases}
\mc{V}^{\la}_{\e'}  & \text{if $\la\in \mc{P}(GL_\ell)_{(\e'_-,\e'_+)}$},\\
0 & \text{otherwise},
\end{cases}
\end{equation*}
where $\e'=\underline{\e}$ or $\overline{\e}$ \eqref{eq:e'(a,b)}.
\end{thm}
\pf Recall the following well-known decomposition (cf.~ Remark \ref{rem:odd homogeneous}):
\begin{equation*}\label{eq:type A 1}
\W_{\ov{\e}}(x)^{\ot \ell} \cong 
\bigoplus_{\la\in \mc{P}(GL_\ell)_{(\ov{\e}_-,\ov{\e}_+)}}
{\mc{V}_{\ov{\e}}^\la}^{\oplus d^\la},
\end{equation*}
where $d^\la$ is the dimension of $V_{GL_\ell}(\la)$. 
By Proposition \ref{prop:hwsub-W^l-super}, we have
\begin{equation}\label{eq:type A 01}
\W_{\e}(x)^{\ot \ell} \cong \bigoplus_{\la\in S}\ {\mc{V}_{\e}^\la}^{\oplus d^{\la}_{\e}},
\end{equation}
for some $S\subset \mc{P}(GL_\ell)_{(\e_-,\e_+)}$ and $d_{\e}^\la\in \mathbb{N}$ ($\la\in S$). \vskip 2mm

{\em Step 1}.
Suppose that $\e'=\ov{\e}$. Let $\la\in S$ be given.
We claim that 
\begin{equation}\label{eq:step 1}
\mf{tr}^\e_{\ov{\e}}\left(\mc{V}^{\la}_\e\right)\cong
\begin{cases}
\mc{V}^{\la}_{\ov{\e}}  & \text{if $\la\in \mc{P}(GL_\ell)_{(\ov{\e}_-,\ov{\e}_+)}$},\\
0 & \text{otherwise},
\end{cases}
\end{equation}

Let $V=\mc{V}_{\e}^\la$ with a highest weight vector $v$. By Lemma \ref{lem:semisimple over levi}, we have $v\in W:=V_{(\e_-,\e_+)}(\la^-,\la^+)$ for $\la^-\in \cP_{a|a}$ and $\la^+\in \cP_{b+1|b}$.

Suppose first that $\la\in \mc{P}(GL_\ell)_{(\ov{\e}_-,\ov{\e}_+)}$.
Then $\mf{tr}^\e_{\overline \e}\left(W\right)$ is non-zero by Lemma \ref{lem:truncation of poly}, and hence $\mf{tr}^\e_{\overline \e}(V)$ is non-zero. Moreover, as in the proof of Proposition~\ref{prop:hwsub-W^l-super}, $\mf{tr}^\e_{\overline \e}(V)$ contains a $\mathring{\mathcal{U}}(\overline{\epsilon})$-highest weight vector of highest weight $\Lambda_{\lambda,\overline{\epsilon}}$, so that $\mathcal{V}_{\overline{\epsilon}}^{\lambda}\subset \mf{tr}^\e_{\overline \e}(V)$.

Consider the $\mathring{U}(\e)$-module $V^\tau$. 
By Corollary \ref{cor:classical limit of h.w. mod},  
its classical limit $\overline{V^\tau}$ is a highest weight $U(\gl(\e)^e)$-module with highest weight $\La_{\la,\e}$. 
We consider $\mf{tr}^\e_{\ov{\e}}\left(\ov{V^\tau}\right)$, which can be defined in the same way as in \eqref{eq:truncation-1}, and which is a representation of $U(\gl(\ov{\e})^e)$ (cf.~Remark \ref{rem:super duality}(2)).
Now we can apply the same method in the proof of \cite[Lemma 3.5]{CL} to conclude that $\mf{tr}^\e_{\ov{\e}}\left(\ov{V^\tau}\right)$ is also a highest weight $U(\gl(\ov{\e})^e)$-module with highest weight $\La_{\la,\ov{\e}}$, and hence irreducible since it is integrable. 
This implies that the character of $\mf{tr}^\e_{\ov{\e}}\left(\ov{V^\tau}\right)$ is the equal to that of $\mc{V}_{\ov{\e}}^\la$.
On the other hand, the character of $\mf{tr}^\e_{\ov{\e}}\left(\ov{V^\tau}\right)$ is equal to that of $\mf{tr}^\e_{\ov{\e}}\left(V\right)$. Therefore, $\mf{tr}^\e_{\ov{\e}}\left(V\right)= \V^\la_{\ov{\e}}$.

Next, suppose that $\la\not\in \mc{P}(GL_\ell)_{(\ov{\e}_-,\ov{\e}_+)}$.
Choose $\td{\e}=\e^{(a+k,b+k)}$ for some $k\ge 1$ such that $\la\in \mc{P}(GL_\ell)_{(\ov{\td{\e}}_-,\ov{\td{\e}}_+)}$. 
By Lemma \ref{lem:truncation by rank}, we have $\V^\la_{\td{\e}}\subset \W_{\td{\e}}(x)^{\ot \ell}$ and $\mf{tr}^{\td{\e}}_\e\left(\V^\la_{\td{\e}}\right)=V$. 
Moreover, by {\em Step 1}, we have 
$\mf{tr}^{\td{\e}}_{\ov{\td{\e}}}\left(\V^\la_{\td{\e}}\right)\cong \V^\la_{\ov{\td{\e}}}$, which is finite-dimensional.
By considering the character of $\V^\la_{\ov{\td{\e}}}$, it is not difficult to see that 
$\mf{tr}^{\ov{\td{\e}}}_{\ov{\e}}(\V^\la_{\ov{\td{\e}}})=0$. Hence
\begin{equation*}
0= \mf{tr}^{\ov{\td{\e}}}_{\ov{\e}}(\V^\la_{\ov{\td{\e}}})=\mf{tr}^{\ov{\td{\e}}}_{\ov{\e}}\circ \mf{tr}^{\td{\e}}_{\ov{\td{\e}}}\left(\V^\la_{\td{\e}}\right)
 =\mf{tr}^{\e}_{\ov{\e}}\circ \mf{tr}^{\td{\e}}_{\e}\left(\V^\la_{\td{\e}}\right)
 =\mf{tr}^{\e}_{\ov{\e}}\left(V\right).
\end{equation*}
This proves the claim.\vskip 2mm

{\em Step 2}.
We claim $S=\mc{P}(GL_\ell)_{(\e_-,\e_+)}$ and $d^\la_\e=d^\la$ in \eqref{eq:type A 01}.

Suppose that $\la\in \mc{P}(GL_\ell)_{(\e_-,\e_+)}$ is given. Consider $\td{\e}=\e^{(a+k,b+k)}$ for a sufficiently large $k\ge 1$ such that $\la\in \mc{P}(GL_\ell)_{(\ov{\tilde{\e}}_-,\ov{\tilde{\e}}_+)}$.
By {\em Step 1} \eqref{eq:step 1} and Lemma \ref{lem:truncation of fund}, we have $\V^\la_{\tilde{\e}}\subset \W_{\tilde{\e}}(x)^{\ot \ell}$.
By Lemma \ref{lem:truncation by rank}(2), we have $\V^\la_{\e} = \mf{tr}^{\td{\e}}_\e\left(\V^\la_{\td{\e}}\right)\subset \W_\e(x)^{\ot \ell}$, that is, $\la\in S$. Therefore we have $S=\mc{P}(GL_\ell)_{(\e_-,\e_+)}$ and by \eqref{eq:step 1}, $d^\la_{\e}=d^\la$.
\vskip 2mm 

{\em Step 3}. Suppose that $\e'=\underline{\e}$.
Let $\la\in \mc{P}(GL_\ell)_{(\underline{\e}_-,\underline{\e}_+)}$ be given. Note that $\mc{P}(GL_\ell)_{(\underline{\e}_-,\underline{\e}_+)}\subset \mc{P}(GL_\ell)_{(\e_-,\e_+)}$.
Let $v$ be a highest weight vector in ${\mc{V}_{\e}^\la}$.
By {\em Step 2}, we have $\V^\la_\e\subset \W_{\e}(x)^{\ot \ell}$.
So we may assume that $v\in V_{(\e_-,\e_+)}(\la^-,\la^+)$ by Lemma \ref{lem:semisimple over levi}.
Let $w$ be a $\mathring{\U}(\underline{\e}_-,\underline{\e}_+)$-highest weight vector in $V_{(\underline{\e}_-,\underline{\e}_+)}(\la^-,\la^+)=\mf{tr}^\e_{\underline{\e}}\left(V_{(\e_-,\e_+)}(\la^-,\la^+) \right)$. 
By minimality of $|\la^-|+|\la^+|$ in ${\mc{V}_{\e}^\la}$, it follows that $w$ is a $\mathring{\U}(\underline{\e})$-highest weight vector.
By Proposition \ref{prop:hwsub-W^l}, the highest weight $\mathring{\U}(\underline{\e})$-module generated by $w$ is isomorphic to $\mc{V}^\la_{\underline{\e}}$ and
\begin{equation}\label{eq:tr to 0}
 \mc{V}^\la_{\underline{\e}} \subset 
 \mf{tr}^\e_{\underline{\e}}\left( \mc{V}^\la_{\e}\right).
\end{equation}
By Theorem \ref{thm:Howe duality} and Corollary \ref{cor:classical limit of V lambda}, we also have
\begin{equation}\label{eq:type A 0}
\W_{\underline{\e}}(x)^{\ot \ell} \cong \bigoplus_{\la\in \mc{P}(GL_\ell)_{(\underline{\e}_-,\underline{\e}_+)}}{\mc{V}_{\underline{\e}}^\la}^{\oplus d^\la}.
\end{equation}
By Lemma \ref{lem:truncation of fund}, we conclude from {\em Step 2} and \eqref{eq:type A 0} that the equality holds in \eqref{eq:tr to 0} and 
$\mf{tr}^\e_{\underline{\e}}\left( \mc{V}^\la_{\e}\right)=0$ for $\la\in \mc{P}(GL_\ell)_{(\e_-,\e_+)}\setminus \mc{P}(GL_\ell)_{(\underline{\e}_-,\underline{\e}_+)}$.
\qed

\begin{cor}\label{cor:decomp of W for e}
We have the following decomposition:
\begin{equation*} 
\W_{\e}(x)^{\ot \ell} \cong \bigoplus_{\la\in \mc{P}(GL_\ell)_{(\e_-,\e_+)}}\ {\mc{V}_{\e}^\la}^{\oplus d^{\la}}.
\end{equation*}
\end{cor}

We have an analogue of Definition \ref{def:cat osc} for $\mathring{\U}(\e)$, say 
$\mc{O}_{{\rm osc},\e}$, where we replace a $U_q(\gl_n)$-module $\V^\la$ for $\la\in \mc{P}(GL_\ell)_{(r,n-r)}$ with a $\mathring{\U}(\e)$-module $\V_{\e}^\la$ for $\la\in \mc{P}(GL_\ell)_{(\e_-,\e_+)}$.
It is closed under tensor product by Corollary \ref{cor:decomp of W for e}, and hence it is a semisimple tensor category. 
By Theorem \ref{thm:trunc sends simple to simple}, the decomposition of the tensor product $\mc{V}^\la_{\e}\ot \mc{V}^\mu_{\e}$ is the same as in the case of $\mc{V}^\la\ot \mc{V}^\mu$ up to irreducible components.

\subsection{Spectral decomposition}\label{subsec:spectral decomp for e}

Let $\e$ and $\underline{\e}, \ov{\e}$ be as in \eqref{eq:e(a,b)} and \eqref{eq:e'(a,b)}.

The notions in Sections \ref{subsec:aff osc cat type 0} and \ref{subsec:Irreducibility} can be naturally defined for $\U(\e)$. For example, we denote by $\widehat{\mc O}_{{\rm osc},\e}$ the one corresponding to $\widehat{\mc O}_{\rm osc}$.  
Note that $\widehat{\mc O}_{{\rm osc},\e}$ is a full subcategory of $\mc{C}(\e)$ (cf.~Section \ref{subsec:truncation functor}).

For $l,m\in \Z$, we regard $\W_{l,\e}=\W_{l,\e}(1)$ and $\W_{m,\e}=\W_{m,\e}(1)$ as $\U(\e)$-modules. 
By \eqref{eq:LR for q fundamentals} and Theorem \ref{thm:trunc sends simple to simple}, we have the following decomposition as a $\mathring{\U}(\e)$-module:
\begin{equation}\label{eq:LR for q fundamentals e}
\W_{l,\e}\ot \W_{m,\e}\cong \W_{m,\e}\ot \W_{l,\e} 
=\bigoplus_{t\in\Z_+} \mc{V}^{(l_1+t,l_2-t)}_{\e}
=\bigoplus_{i\ge -L} \mc{V}^{(l_1+L+i,l_2-L-i)}_{\e},
\end{equation}
where $l_1=\max\{l,m\}$ and $l_2=\min\{l,m\}$ and $L$ is given in \eqref{eq:L}.

\begin{thm}\label{thm:irr of tensor product of fund e}
The tensor product $\W_{l,\e}\ot \W_{m,\e}$ is an irreducible $\U(\e)$-module in $\widehat{\mc O}_{{\rm osc},\e}$.
\end{thm}
\pf We may assume that 
\begin{equation*}\label{eq:embedding of 0 into e}
 \W_{l,\underline{\e}}\ot \W_{m,\underline{\e}} \subset \W_{l,\e}\ot \W_{m,\e},
\end{equation*}
as a $\Bbbk$-space. 
Let $v\in \W_{l,\e}\ot \W_{m,\e}$ be given. By our choice of $\e$, we may send $v$ to a non-zero vector $v'$ in $\W_{l,\underline{\e}}\ot \W_{m,\underline{\e}}$ by applying finitely many $E_{i}$, $F_{i-1}$ (if $i\leq r$) or $E_{i-1}$, $F_{i}$ (if $i> r$) with $\e_i=1$. 
Since $\W_{l,\underline{\e}}\ot \W_{m,\underline{\e}}$ is irreducible by Theorem \ref{thm:irr of tensor product of fund 0}, the $\U(\underline{\e})$-submodule generated by $v'$ contains all $\mc{V}^{(l_1+t,l_2-t)}_{\underline{\e}}$ ($t\in\Z_+$).
Since $\mc{V}^{(l_1+t,l_2-t)}_{\underline{\e}}\subset \mc{V}^{(l_1+t,l_2-t)}_{\e}$
and $\mc{V}^{(l_1+t,l_2-t)}_{\e}$ is an irreducible $\mathring{\U}(\e)$-module, the $\U(\e)$-submodule generated by $v'$ is $\W_{l,\e}\ot \W_{m,\e}$ by \eqref{eq:LR for q fundamentals e}. Hence $\W_{l,\e}\ot \W_{m,\e}$ is irreducible.
\qed

\begin{cor}\label{cor:irreducibility e}
The tensor product $(\mathcal{W}_{l,\e})^{{}^{\wedge}}_{\rm aff} \ot (\mathcal{W}_{m,\e})^{{}^{\wedge}}_{\rm aff}$ is an irreducible $\Bbbk(z_1,z_2)\ot \U(\e)$-module.
\end{cor}

The universal $R$ matrix for $\U(\e)$ is well-defined by \cite[Theorem 3.7]{KL} since the numbers of $0$'s and $1$'s are not equal by our choice of $\e$.

Let us keep the notations in Section \ref{subsec:normalized R matrix}.
We assume that the vector $v_l$ in \eqref{eq:v_l} belongs to $ \W_{l,\e}$ by identifying $v_l=\ket{-l\be_{r-1}}$ if $l<0$ and $v_l=\ket{l\be_{r+1}}$ if $l>0$.
Then the universal $R$ matrix for $\U(\e)$ induces a map
\begin{equation*}
\xymatrixcolsep{4pc}\xymatrixrowsep{3pc}\xymatrix{
(\mathcal{W}_{l,\e})_{\rm aff}\, \widehat{\ot}\,(\mathcal{W}_{m,\e})_{\rm aff} \ \ar@{->}[r] &\ (\mathcal{W}_{m,\e})_{\rm aff}\, \widetilde{\ot}\, (\mathcal{W}_{l,\e})_{\rm aff}}, 
\end{equation*}
such that
\begin{equation*}
\mc{R}^{\rm univ}_{(l,m)}(v_l\ot v_m)= a_\e\left(\tfrac{z_1}{z_2}\right)\ot(v_m\ot v_l),
\end{equation*}
for some invertible $a_\e(\frac{z_1}{z_2})\in \Bbbk\llbracket \frac{z_1}{z_2}\rrbracket$. 
Then we define the normalized $R$ matrix by
\begin{equation*}
\mc{R}^{\rm norm}_{(l,m),\e}= c(z)a_\e(z)^{-1}\mc{R}^{\rm univ}_{(l,m),\e},
\end{equation*}
where $z={z_1}/{z_2}$ and $c(z)$ is given in \eqref{eq:c(z)}.
We have a unique $\Bbbk[z_1^{\pm 1},z_2^{\pm 1}]\ot \U(\e)$-linear map
\begin{equation}\label{eq:R matrix on W_l tensor W_m for e}
\xymatrixcolsep{3pc}\xymatrixrowsep{3pc}\xymatrix{
(\mathcal{W}_{l,\e})_{\rm aff} \ot (\mathcal{W}_{m,\e})_{\rm aff}\ \ar@{->}^{\hskip -2cm \mc{R}^{\rm norm}_{(l,m),\e}}[r] &  \
\Bbbk(z_1,z_2)\ot_{\Bbbk[z_1^{\pm 1},z_2^{\pm 1}]}\left((\mathcal{W}_{m,\e})_{\rm aff}\,{\ot}\, (\mathcal{W}_{l,\e})_{\rm aff}\right),
}
\end{equation}
satisfying 
\begin{equation}\label{eq:normalization e}
\mc{R}^{\rm norm}_{(l,m),\e}\big\vert_{\V_\e^{(l_1+L,l_2-L)}}
=c(z){\rm id}_{\V_\e^{(l_1+L,l_2-L)}}. 
\end{equation}

\begin{lem}\label{lem:truncation from e to (0)}%\comments{This lemma has been added for clarification (K.2022.01.29)}
Under the above hypothesis, we have
$$\mf{tr}^{\e}_{\underline{\e}}\left( \mc{R}^{\rm norm}_{(l,m),{\e}} \right)=
 \mc{R}^{\rm norm}_{(l,m),\underline{\e}}.$$
\end{lem}
\pf By Lemma \ref{lem:truncation of fund} and \eqref{eq:R matrix on W_l tensor W_m for e}, $\mf{tr}^{\e}_{\underline{\e}}\left( \mc{R}^{\rm norm}_{(l,m),{\e}} \right)$ induces a $\mc{U}(\underline{\e})$-linear map 
\begin{equation*}
\xymatrixcolsep{3pc}\xymatrixrowsep{3pc}\xymatrix{
(\mathcal{W}_{l})_{\rm aff} \ot (\mathcal{W}_{m})_{\rm aff}\ \ar@{->}^{\hskip -2cm \mc{R}^{\rm norm}_{(l,m)}}[r] &  \
\Bbbk(z_1,z_2)\ot_{\Bbbk[z_1^{\pm 1},z_2^{\pm 1}]}\left((\mathcal{W}_{m})_{\rm aff}\,{\ot}\, (\mathcal{W}_{l})_{\rm aff}\right),
}
\end{equation*}
which satisfies  $\mf{tr}^{\e}_{\underline{\e}}\left( \mc{R}^{\rm norm}_{(l,m),{\e}} \right)(v_l\ot v_m)=c(z)v_m\ot v_l$ by \eqref{eq:normalization e}.
Hence it coincides with $\mc{R}^{\rm norm}_{(l,m),\underline{\e}}$ in \eqref{eq:R matrix on W_l tensor W_m} by its uniqueness.
\qed\newline

For $t\in\Z_+$, we define a $\mathring{\U}(\e)$-linear map $\mc P^{l,m}_{t} : \W_{l,\e} \otimes \W_{m,\e} \longrightarrow \W_{m,\e} \otimes \W_{l,\e}$ as follows:

\begin{itemize}
 \item[(1)] Let $\td{\e}=\e^{(a+k,b+k)}$ for some $k>0$ such that $\V^{(l_1+t,l_2-t)}_{\ov{\td{\e}}}\neq 0$.
 
 \item[(2)] Let $v_0(l,m,t)$ and $v_0(m,l,t)$ be the $\mathring{\U}(\ov{\td{\e}}\,)$-highest weight vectors of $\mc{V}^{(l_1+t,l_2-t)}_{\ov{\td{\e}}}$ in $\W_{l,\td{\e}} \otimes \W_{m,\td{\e}}$ and $\W_{m,\td{\e}} \otimes \W_{l,\td{\e}}$ respectively, which correspond to the elements in their crystal bases at $\td{q}=0$.
 
 \item[(3)] Let $\td{\mc{P}}^{l,m}_t : \W_{l,\td{\e}} \otimes \W_{m,\td{\e}} \longrightarrow \W_{m,\td{\e}} \otimes \W_{l,\td{\e}}$ be a $\mathring{\U}(\td{\e})$-linear map given by
\begin{equation*}
\begin{split}
&\td{\mc{P}}^{l,m}_t(v_0(l,m,t))=v_0(m,l,t),\ \
\td{\mc{P}}^{l,m}_t\left(\V_{\td{\e}}^{(l_1+t',l_2-t')}\right)=0 \quad (t'\neq t).
\end{split}
\end{equation*}
 
 \item[(4)] Define 
\begin{equation}\label{eq:normalized projector}
 \mc P^{l,m}_{t}=\mf{tr}^{\td{\e}}_{\e}\left(\td{\mc{P}}^{l,m}_t\right).
\end{equation}
\end{itemize}
Note that the map $\mc P^{l,m}_{t}$ is well-defined since the highest weight vectors $v_0(l,m,t)$ and $v_0(m,l,t)$ do not depend on the choice of $\td{\e}$ or $k$. 
Let $\kappa\in\Bbbk^\times$ be such that 
\begin{equation}\label{eq:kappa}
 \mc{P}^{l,m}_L(v_l\ot v_m) = \kappa v_m\ot v_l.
\end{equation}
Then we have the following spectral decomposition of $\mc{R}^{\rm norm}_{(l,m),\e}$.

\begin{thm}\label{thm:spectral decomposition for e}
For $l, m\in\mathbb{Z}$, we have
\begin{equation*}
\begin{split}
\kappa\mc{R}^{\rm norm}_{(l,m),\e} 
&= \mc P^{l,m}_0 + \sum_{t=1}^\infty\prod_{i=1}^{t}\dfrac{1-q^{|l-m|+2i}z}{z-q^{|l-m|+2i}} \mc P^{l,m}_t\\
&= \mc P^{l,m}_0 
+ \dfrac{1-q^{|l-m|+2}z}{z-q^{|l-m|+2}} \mc P^{l,m}_1 
+ \dfrac{1-q^{|l-m|+2}z}{z-q^{|l-m|+2}}\dfrac{1-q^{|l-m|+4}z}{z-q^{|l-m|+4}} \mc P^{l,m}_2 
+ \dots .
\end{split}
\end{equation*}
\end{thm}
\pf 
Let us write
\begin{equation*}
 \mc{R}^{\rm norm}_{(l,m),\e} 
= \sum_{t=0}^\infty\rho_{t,\e}(z) \mc P^{l,m}_t,
\end{equation*}
where $\rho_{L,\e}(z)=c(z)\kappa^{-1}$.
Let $\td{\e}=\e^{(a+k,b+k)}$ for $k\ge 1$. 
By the uniqueness of the normalized $R$ matrix, we have
\begin{equation*}
 \mf{tr}^{\td{\e}}_{\e}\left( \mc{R}^{\rm norm}_{(l,m),\td{\e}} \right)=
 \mc{R}^{\rm norm}_{(l,m),\e}.
\end{equation*}
So we may assume that $a, b$ in $\e=\e^{(a,b)}$ are large enough to compute $\rho_{t,\e}(z)$.
Suppose that $0<a+l,a+m<a+b$.

Note that $\U(\ov{\e})$ is isomorphic to $U_{\td{q}}(\agl_{a+b})$, where $\td{q}=-q^{-1}$. 
Then $\W_{l,\ov{\e}}(x)$ is isomorphic to the $(a+l)$-th fundamental representation with spectral parameter $x$, and 
$\mc{V}^{(l_1+t,l_2-t)}_{\ov{\e}}=\mc{V}^{(l+m-l_2+t,l_2-t)}_{\ov{\e}}$ is isomorphic to the finite-dimensional irreducible $U_{\td{q}}(\gl_{a+b})$-module corresponding to the partition with two columns of lengths $l'_1\ge l'_2$ given by   
\begin{equation*}
(l'_1,l'_2)=((a+l)+(a+m)-(a+l_2)+t,(a+l_2)-t).
\end{equation*}

Let $\mc{R}^{\rm norm}_{(l,m),\ov{\e}}$ be the unique $\U(\ov{\e})$-linear map such that 
\begin{equation*}
\xymatrixcolsep{3pc}\xymatrixrowsep{1pc}\xymatrix{
(\mathcal{W}_{l,\ov{\e}})_{\rm aff} \ot (\mathcal{W}_{m,\ov{\e}})_{\rm aff}\ \ar@{->}^{\hskip -2cm \mc{R}^{\rm norm}_{(l,m),\ov{\e}}}[r] &  \
\Bbbk(z_1,z_2)\ot_{\Bbbk[z_1^{\pm 1},z_2^{\pm 1}]}\left((\mathcal{W}_{m,\ov{\e}})_{\rm aff}\,{\ot}\, (\mathcal{W}_{l,\ov{\e}})_{\rm aff}\right) \\
v_0(l,m,0) \ \ar@{|->}[r] & \ v_0(m,l,0)
}.
\end{equation*}
By the uniqueness of the normalized $R$ matrix \eqref{eq:R matrix on W_l tensor W_m for e}, we have
$
\mf{tr}^{\e}_{\ov{\e}}\left( \mc{R}^{\rm norm}_{(l,m),{\e}} \right)=
 c\mc{R}^{\rm norm}_{(l,m),\ov{\e}},
$
for some $c\in \Bbbk^\times$, 
which implies
\begin{equation*}
 \mc{R}^{\rm norm}_{(l,m),\ov{\e}} 
= \sum_{t=0}^{a+l_2}c^{-1}\rho_{t,\e}(z) \mf{tr}^{\e}_{\ov{\e}}\left(\mc P^{l,m}_t\right).
\end{equation*}
By construction of $\mc{P}^{l,m}_t$ in \eqref{eq:normalized projector}, the map $\mf{tr}^{\e}_{\ov{\e}}\left(\mc P^{l,m}_t\right) : \W_{l,\ov{\e}} \otimes \W_{m,\ov{\e}} \longrightarrow \W_{m,\ov{\e}} \otimes \W_{l,\ov{\e}}$ is a $\mathring{\U}(\ov{\e})$-linear map sending $v_0(l,m,t)$ to $v_0(m,l,t)$ on $\V^{(l_1+t,l_2-t)}$ and the other components $\V^{(l_1+t',l_2-t')}$ ($t'\neq t$) to $0$.

Therefore it follows from \cite{DO} that 
putting $t'=a+l_2-t$ ($0\le t\le a+l_2$), the coefficient $c^{-1}\rho_{t,\e}(z)$ of $\mf{tr}^{\e}_{\ov{\e}}\left(\mc P^{l,m}_{t}\right)$ is given by 
\begin{equation}\label{eq:DO}
\prod_{i=t'+1}^{a+l_2}\frac{1-q^{2a+l+m-2i+2}z}{z-q^{2a+l+m-2i+2}}
=\prod_{i=1}^{t}\frac{1-q^{|l-m|+2i}z}{z-q^{|l-m|+2i}},
\end{equation}
where it is assumed to be $1$ for $t=0$, and the second equality follows from 
\begin{equation*}
 \begin{split}
  &2a+l+m-2(a+l_2)+2=|l-m|+2,\\
  &2a+l+m-2t'=2a+l+m-2a-2l_2+2t=|l-m|+2t.
 \end{split}
\end{equation*}
(Here we follow the convention in \cite[Theorem 3.15]{KL} for \eqref{eq:DO}.)
In particular, if $t'=a+l_2-L$, that is, $t=L$, 
then we have by \eqref{eq:DO}
\begin{equation*}
 c^{-1}\rho_{L,\e}(z)=c(z).
\end{equation*}
Since $\rho_{L,\e}(z)=c(z)\kappa^{-1}$, we have $c=\kappa^{-1}$, and therefore  
\begin{equation*}\label{eq:spec coeff rho}
 \rho_{t,\e}(z)
 =\kappa^{-1}\prod_{i=1}^{t}\frac{1-q^{|l-m|+2i}z}{z-q^{|l-m|+2i}}.
\end{equation*}
This completes the proof.
\qed

\begin{cor}\label{cor:truncation from e to (1)}%\comments{This corollary has been added for clarification (K.2022.01.29)}
 Under the above hypothesis, we have
 $$\mf{tr}^{\e}_{\ov{\e}}\left(\kappa \mc{R}^{\rm norm}_{(l,m),{\e}} \right)=
 \mc{R}^{\rm norm}_{(l,m),\ov{\e}}.$$
\end{cor}

Finally, the following corollary implies Theorem \ref{thm:spectral decomposition}.

\begin{cor}\label{cor:spectral decomp}
For $l,m\in \Z$, we have
\begin{equation*}
\begin{split}
\kappa\mc{R}^{\rm norm}_{(l,m),\underline{\e}} 
&= \mc P^{l,m}_0 + \sum_{t=1}^\infty\prod_{i=1}^{t}\dfrac{1-q^{|l-m|+2i}z}{z-q^{|l-m|+2i}} \mc P^{l,m}_t.
\end{split}
\end{equation*}
\end{cor}
\pf We apply Lemma \ref{lem:truncation from e to (0)} to Theorem \ref{thm:spectral decomposition for e}.
\qed


{\small
\begin{thebibliography}{HK}

\bibitem{BCP}
J. Beck, V. Chari, A. Pressley, {\em An algebraic characterization of the affine canonical basis}, Duke Math. J. {\bf 99} (1999) 455--487.

\bibitem{BKK}
G. Benkart, S.-J. Kang, M. Kashiwara, {\em Crystal bases for the
quantum superalgebra $U_q(\mathfrak{gl}(m,n))$}, J. Amer. Math. Soc. \textbf{13} (2000) 295-331.

%\bibitem{BR} 
%A.~Berele, A.~Regev, {\em Hook Young diagrams with applications to combinatorics and representations of Lie superalgebras}, Adv.~Math.~{\bf 64} (1987), 118--175.

\bibitem{BGKNR16}
H. Boos, F. G\"ohmann, A. Kl\"umper, K. S. Nirov, A. V. Razumov, {\em Oscillator versus prefundamental representations}, J. Math. Phys. {\bf 57} (2016) 111702.

\bibitem{BGKNR17}
H. Boos, F. G\"ohmann, A. Kl\"umper, K. S. Nirov, A. V. Razumov, {\em Oscillator versus prefundamental representations II. Arbitrary higher ranks}, J. Math. Phys. {\bf 58} (2017) 093504.


\bibitem{CP}
V. Chari, A. Pressley, {\em Quantum affine algebras and their representations}, in Representations of groups (Banff, AB, 1994), CMS Conference Proceedings, Vol. 16, American Mathematical Society, Providence, RI (1995) 59--78.

%\bibitem{CP}
%V. Chari, A. Pressley, {\em Quantum affine algebras and affine Hecke algebras}, Pacific J. Math. \textbf{174} (2) (1996) 295--326.



\bibitem{CK08} 
S.-J.~Cheng, J.-H.~Kwon, {\em Howe duality and
Kostant's homology formula for infinite dimensional Lie
superalgebras}, Int. Math. Res. Not. IMRN 2008, Art. ID rnn 085, 52pp.

\bibitem{CL}
S.-J. Cheng, N. Lam, {\em Irreducible characters of general linear superalgebra and super duality}, Comm. Math. Phys. \textbf{298} (2010) 645--672.

%\bibitem{CLW}
%S.-J. Cheng, N. Lam, W. Wang, {\em Super duality and irreducible characters of ortho-symplectic Lie superalgebras}, Invent. Math. \textbf{183} (2011) 189--224.

\bibitem{CLZ}
S.-J. Cheng, N. Lam, R. B. Zhang, {\em Character formula for infinite-dimensional unitarizable modules of the general linear superalgebra},  
J. Algebra \textbf{273} (2004) no. 2, 780–-805. 

\bibitem{CW}
S.-J. Cheng, W. Wang, {\em Dualities and Representations of Lie Superalgebras}, Graduate Studies in Mathematics 144, Amer. Math. Soc., 2013.

\bibitem{DO}
E. Date, M. Okado, {\em Calculation of excitation spectra of the spin model related with the vector representation of the quantized affine algebra of type $A_n^{(1)}$}, Internat. J. Modern Phys. A  \textbf{9} (1994) 399--417.

\bibitem{EHW}
T. Enright, R. Howe, N. Wallach, {\em Classification of unitary highest weight modules}, in Representation Theory of Reductive Groups (P. Trombi, Ed.), pp. 97 – 144, Birkh\"{a}user, Boston, 1983.

\bibitem{FH}
E. Frenkel, D. Hernandez, {\em Baxter’s relations and spectra of quantum integrable models}, Duke Math. J. {\bf 164} (2015) 2407--2460.

\bibitem{H05}
D. Hernandez, {\em Representations of quantum affinizations and fusion product}, Transform. Groups \textbf{10} (2005) 163--200.

\bibitem{HJ}
D. Hernandez, M. Jimbo, {\em Asymptotic representations and Drinfeld rational fractions}, Compositio Math. {\bf 148} (2012) 1593--1623.

\bibitem{HL10}
D. Hernandez, B. Leclerc, {\em Cluster algebras and quantum affine algebras}, Duke Math. J. \textbf{154} (2010) 265--341.

%\bibitem{HL16}
%D. Hernandez, B. Leclerc, {\em Cluster algebras and category $\mc{O}$ for representations of Borel subalgebras of quantum affine algebras}, Algebra Number Theory \textbf{10} (2016) 2015--2052.

\bibitem{H} 
R.~Howe, {\em Remarks on classical invariant theory}, Trans. Amer. Math. Soc. {\bf 313} (1989) 539--570.

\bibitem{HTW}
R. Howe, E.-C. Tan, J. Willenbring, {\em Stable branching rules for classical symmetric pairs}, Trans. Amer. Math. Soc. \textbf{357}  (2004) 1601--1626.

\bibitem{JKP}
I.-S. Jang, J.-H. Kwon, E. Park, {\em Unipotent quantum coordinate ring and prefundamental representations for types $A_n^{(1)}$ and $D_n^{(1)}$}, Int. Math. Res. Not. \textbf{2} (2023) 1119--1172.

\bibitem{Ja}
J.C. Jantzen, {\em Lectures on Quantum Groups}, Amer. Math. Soc. Graduate Studies in Math., vol. 6, Providence, 1996. 

\bibitem{KKK}
S.-J. Kang, M. Kashiwara, M. Kim, {\em Symmetric quiver Hecke algebras and $R$-matrices of quantum affine algebras}, Invent. math. \textbf{211} (2018) 591–685.

\bibitem{KKK2}
S.-J. Kang, M. Kashiwara, M. Kim, {\em Symmetric quiver Hecke algebras and $R$-matrices of quantum affine algebras, II}, Duke. Math. J. \textbf{164} (2015) 1549--1602.

\bibitem{KMN}
S.-J. Kang, M. Kashiwara, K. C. Misra, T. Miwa, T. Nakashima and A. Nakayashiki, 
{\em Affine crystals and vertex models}, Internat. J. Modern Phys. A \textbf{7} (suppl. 1A) (1992) 449--484.

\bibitem{Kas02}
M. Kashiwara, {\em On level zero representations of quantum affine algebras}, Duke Math. J. \textbf{112} (2002) 117--175.

\bibitem{KK}
M. Kashiwara, M. Kim, {\em Laurent phenomenon and simple modules of quiver Hecke algebras}, Compos. Math. \textbf{155.12} (2019) 2263--2295.

\bibitem{KKKO15}
S.-J. Kang, M. Kashiwara, M.-H. Kim, S.-J. Oh, {\em Simplicity of heads and socles of tensor products}, Compos. Math. \textbf{151} (2015) 377--396.

\bibitem{KKOP1}
M. Kashiwara, M. Kim, S. Oh, E. Park, {\em Monoidal categorification and quantum affine algebras}, Compos. Math. \textbf{156} (2020) 1039--1077.

\bibitem{KKOP2}
M. Kashiwara, M. Kim, S. Oh, E. Park, {\em Monoidal categorification and quantum affine algebras II} preprint (2021) arXiv:2103.10067.

\bibitem{KV}
M. Kashiwara, M. Vergne, {\em On the Segal-Shale-Weil
representations and harmonic polynomials}, Invent. Math. \textbf{44}
(1978) 1--47.

\bibitem{Ku18}
A. Kuniba, {\em Tetrahedron equation and quantum R matrices for q-oscillator representations mixing particles and holes}, SIGMA Symmetry Integrability Geom. Methods Appl. 14 (2018), Paper No. 067

\bibitem{KuO}
A. Kuniba, M. Okado,
{\em Tetrahedron equation and quantum R matrices for infinite-dimensional modules of $U_q(A^{(1)}_1)$ and $U_q(A^{(2)}_2)$}, J. Phys. A \textbf{46} (2013) 485203.

\bibitem{KOS}
A. Kuniba, M. Okado, S. Sergeev, {\em Tetrahedron equation and generalized quantum groups}, 
J. Phys. A \textbf{48} (2015) 304001 (38p).

\bibitem{K08}
J.-H. Kwon, {\em Rational semistandard tableaux and character
formula for the Lie superalgebra $\widehat{\gl}_{\infty|\infty}$},  Adv. Math. \textbf{217} (2008) 713--739.

\bibitem{KL}
J.-H. Kwon, S.-M. Lee, {\em Super duality for quantum affine algebras of type $A$}, Int. Math. Res. Not. \textbf{23} (2022) 18446--18525.

\bibitem{KO21}
J.-H. Kwon, M. Okado, {\em Higher level $q$-oscillator representations for $U_q(C^{(1)}_n)$, $U_q(C^{(2)}(n+1))$ and $U_q(B^{(1)}(0,n))$}, Comm. Math. Phys. \textbf{385} (2021) 1041--1082. 

\bibitem{KY}
J.-H. Kwon, J. Yu, {\em $R$-matrix for generalized quantum groups of type $A$}, J. Algebra, \textbf{566} (2021) 309--341.

\bibitem{MY}
E. Mukhin, C. A. S. Young, {\em Affinization of category $\mathcal{O}$ for quantum groups}, Trans. Amer. Math. Soc. {\bf 366} (2014) 4815--4847.

\bibitem{SS}
B. E. Sagan, R. Stanley, {\em Robinson-Schensted algorithms for skew
tableaux}, J. Combin. Theory Ser. A \textbf{55} (1990), 161--193.

\bibitem{Ya94}
H. Yamane, {\em Quantized enveloping algebras associated to simple Lie superalgebras and universal R-matrices},  Publ. Res. Inst. Math. Sci. \textbf{30} (1994) 15--84.

\bibitem{Ya99}
H. Yamane, {\em On defining relations of affine Lie superalgebras and affine quantized universal enveloping superalgebras}, Publ. Res. Inst. Math. Sci. \textbf{35} (1999) 321--390.

\end{thebibliography}}
\end{document}